# PROLONGEMENT ANALYTIQUE DE FONCTIONS $\zeta$ ET DE FONCTIONS $L$

### par **Pierre Colmez**

## 1. INTRODUCTION

### 1.1. La conjecture de Hasse–Weil

La fonction zêta de Riemann, définie par la formule $\zeta(s) = \sum_{n\geq 1} n^{-s}$ sur le demi-plan $\mathrm{Re}(s) > 1$, admet un prolongement méromorphe à tout le plan complexe, holomorphe en dehors d'un pôle simple en $s = 1$, et une équation fonctionnelle[1]

$$\Gamma_{\mathbf{R}}(s)\zeta(s) = \Gamma_{\mathbf{R}}(1-s)\zeta(1-s)$$

Elle peut aussi se définir comme le produit eulérien $\prod_p (1 - p^{-s})^{-1}$, où le produit porte sur les nombres premiers $p$ ou, de manière équivalente, sur les idéaux maximaux de $\mathbf{Z}$.

Ce dernier point de vue se généralise pour donner naissance aux fonctions zêta de Hasse–Weil : si $\Lambda$ est un anneau de type fini sur $\mathbf{Z}$ (i.e. $\Lambda$ de la forme $\mathbf{Z}[X_1, \ldots, X_n]/I$), on définit sa fonction zêta par la formule $\zeta_\Lambda(s) = \prod_{\mathfrak{m}} (1 - |\Lambda/\mathfrak{m}|^{-s})^{-1}$ où le produit porte sur les idéaux maximaux de $\Lambda$ et $|\Lambda/\mathfrak{m}|$ est le cardinal du corps fini $\Lambda/\mathfrak{m}$. Pour $\Lambda = \mathbf{Z}$, on retrouve la fonction zêta de Riemann.

*Remarque 1.1.* — Si $\mathfrak{m}$ est un idéal maximal de $\Lambda$, alors $\mathfrak{m}$ contient $p$ pour un unique nombre premier $p$ ; il s'ensuit que l'on a une factorisation $\zeta_\Lambda(s) = \prod_p \zeta_{\Lambda/p}(s)$ en produit de facteurs d'Euler.

On montre facilement que, si $\Lambda$ est un quotient de $\mathbf{Z}[X_1, \ldots, X_n]$, le produit converge pour $\mathrm{Re}(s) > n + 1$ et définit une fonction holomorphe sur ce demi-plan.

CONJECTURE 1.2. — (HASSE–WEIL). *La fonction $\zeta_\Lambda$ possède un prolongement méromorphe à tout le plan complexe.*

Voici ce qu'écrit WEIL dans son commentaire sur [218] dans lequel il prouve le premier cas significatif de cette conjecture :

> Peu avant la guerre, si mes souvenirs sont exacts, G. de Rham me raconta qu'un
> de ses étudiants, Pierre Humbert, était allé à Göttingen avec l'intention d'y travailler sous la direction de Hasse, et que celui-ci lui avait proposé un problème sur lequel de Rham désirait mon avis. Une courbe elliptique $C$ étant donnée sur le corps des rationnels, il s'agissait principalement, il me semble, d'étudier le produit infini des fonctions zêta des courbes $C_p$ obtenues en réduisant $C$ modulo $p$ pour tout nombre premier pour lequel $C_p$ est de genre 1 ; plus précisément, il fallait rechercher si ce produit possède un prolongement analytique et une équation fonctionnelle. [...]

---

1. $\Gamma_{\mathbf{R}}(s) = \pi^{-s/2}\Gamma(s/2)$, où $\Gamma$ est la fonction $\Gamma$ d'Euler ; on a $\Gamma_{\mathbf{R}}(s)\Gamma_{\mathbf{R}}(s+1) = \Gamma_{\mathbf{C}}(s) := 2(2\pi)^{-s}\Gamma(s)$.

J'avoue avoir pensé que Hasse avait été trop optimiste. Je dis à de Rham, non seulement que le problème me semblait trop difficile pour un étudiant, mais même que je ne voyais aucune raison pour que le produit en question eût les propriétés que Hasse lui supposait. [...]

Ce qui finalement me donna confiance dans cette idée, ce fut une fois de plus l'analogie entre corps de nombres et corps de fonctions ; j'observai (v. [217, p. 99]) que la conjecture de Hasse pour les courbes sur un corps de nombres correspond exactement à mes conjectures de 1949 pour les surfaces sur un corps fini, au sujet desquelles il ne me restait plus aucun doute.

Contrairement à ce que le cas $\Lambda = \mathbf{Z}$ pourrait laisser croire, cette conjecture est remarquablement difficile à prouver, même à $\Lambda$ fixé. Un des problèmes est qu'il n'y a, en général, aucun moyen de décrire ce que doivent être les pôles du prolongement.

*Remarque 1.3.* — Si $\mathscr{X} = \operatorname{Spec} \Lambda$, on a aussi $\zeta_\Lambda(s) = \prod_x (1 - |\kappa(x)|^{-s})^{-1}$, où $x$ parcourt l'ensemble des points fermés de $\mathscr{X}$, et $\kappa(x)$ est le corps résiduel en $x$. Cette formule a un sens pour tout schéma $\mathscr{X}$ de type fini sur $\mathbf{Z}$, ce qui permet de définir la fonction zêta de Hasse–Weil $\zeta_\mathscr{X}$ d'un tel schéma. Comme on peut partitionner $\mathscr{X}$ en schémas affines $\operatorname{Spec} \Lambda_i$, on a $\zeta_\mathscr{X} = \prod_i \zeta_{\Lambda_i}$ et la conjecture de Hasse–Weil pour $\mathscr{X}$ est donc conséquence de celle pour les $\Lambda_i$.

*1.1.1. Le cas de la dimension* 0. — Il y a quand même un certain nombre de cas où cette conjecture est connue depuis belle lurette.

• $\Lambda = \mathbf{Z}[\sqrt{-1}]$ où $\zeta_\Lambda(s) = \zeta(s) L(\chi, s)$, où $L(\chi, s)$ est la fonction $L$ de Dirichlet $\sum_{n \geq 1} \chi(n) n^{-s} = \prod_p (1 - \chi(p) p^{-s})^{-1}$, où $\chi : (\mathbf{Z}/4\mathbf{Z})^* \to \{\pm 1\}$ est le caractère non trivial (i.e. $\chi(1) = 1$ et $\chi(-1) = -1$), dont le prolongement analytique à tout le plan complexe est élémentaire.

• Plus généralement, $\Lambda = \mathbf{Z}[\boldsymbol{\mu}_N]$ où $\zeta_\Lambda$ est le produit de $\zeta$ et de fonctions $L$ de Dirichlet.

• Encore plus généralement, $\Lambda$ anneau des entiers d'un corps de nombres $F$ ; alors $\zeta_\Lambda$ est la fonction zêta de Dedekind $\zeta_F$ de $F$ ; son prolongement méromorphe à tout le plan complexe a été établi par HECKE [131] ; cette fonction est holomorphe en dehors d'un pôle simple en $s = 1$.

*1.1.2. Le cas où $p = 0$ dans $\Lambda$.* — Comme $\zeta_{\mathscr{U} \cup \mathscr{V}} = \zeta_\mathscr{U} \zeta_\mathscr{V} \zeta_{\mathscr{U} \cap \mathscr{V}}^{-1}$, si $\mathscr{U}$ et $\mathscr{V}$ sont des fermés de Zariski de $\mathscr{X}$, on peut se ramener au cas où $\mathscr{X}$ est une variété algébrique irréductible définie sur $\mathbf{F}_q$. On a

$$-\log \zeta_{\mathscr{X}}(s) = \sum_x \sum_{k\geq 1} \tfrac{1}{k}|\kappa(x)|^{-ks} = \sum_{k\geq 1} \tfrac{1}{k}|\mathscr{X}(\mathbf{F}_{q^k})|\, q^{-ks}$$

(Si $x$ est un point fermé de $\mathscr{X}$ et $|\kappa(x)| = q^k$, alors $x$ correspond à une orbite de cardinal $k$ de $\mathrm{Gal}(\overline{\mathbf{F}}_q/\mathbf{F}_q)$ agissant sur $\mathscr{X}(\mathbf{F}_{q^k})$.)

Une des célèbres conjectures de Weil sur le nombre de points des variétés sur les corps finis [216] affirme que $\zeta_{\mathscr{X}}$ est une fraction rationnelle en $q^{-s}$ (et donc aussi en $p^{-s}$); ceci a été démontré par DWORK [86]. Ceci prouve la conjecture de Hasse–Weil en caractéristique $p$.

*1.1.3. Variétés de dimension $\geq 1$ sur un corps de nombres.* — Pour énoncer plus facilement les résultats suivants, disons que *deux séries de Dirichlet* $F = 1 + a_2 2^{-s} + a_3 3^{-s} + \cdots$ *et* $G = 1 + b_2 2^{-s} + b_3 3^{-s} + \cdots$ *sont équivalentes* (et notons $F \sim G$) s'il existe un ensemble fini $S$ de nombres premiers et, pour $p \in S$, une fonction rationnelle $R_p$ en $p^{-s}$, tels que $G = F \cdot \prod_{p \in S} R_p$.

Par exemple, le théorème de DWORK implique que $\zeta_\Lambda \sim 1$ si $p = 0$ dans $\Lambda$ (c'est le cas, plus généralement, s'il existe un entier $N \geq 1$ tel que $N = 0$ dans $\Lambda$). L'existence d'un prolongement méromorphe à tout le plan complexe pour une série de Dirichlet ne dépend que de sa classe modulo cette relation d'équivalence.

*Remarque 1.4.* — Si $X$ est une variété projective (ou propre) définie sur $\mathbf{Q}$ (ou sur un corps de nombres $F$), on peut en choisir un modèle $\mathscr{X}$ sur $\mathbf{Z}$ (ou sur l'anneau des entiers $\mathscr{O}_F$ de $F$); la classe de $\zeta_{\mathscr{X}}$ modulo $\sim$ ne dépend pas de $\mathscr{X}$ grâce à la rationalité des fonctions zêtas des variétés sur les corps finis. On note $\zeta_X$ cette classe; on dit que $X$ vérifie la conjecture de Hasse–Weil si c'est le cas pour une des $\zeta_{\mathscr{X}}$.

Les exemples ci-dessus concernaient la dimension 0 ou la caractéristique $p$. À partir de maintenant, nous allons nous intéresser uniquement à des variétés définies sur un corps de nombres $F$.

• Si $X$ est une courbe de genre 0 définie sur $F$, alors $\zeta_X(s) \sim \zeta_F(s)\zeta_F(s-1)$, et la conjecture est une conséquence du résultat de HECKE mentionné ci-dessus.

• WEIL [218] a traité le cas d'une courbe $X$ d'équation $y^d = ax^e + b$, avec $a,b \in F^*$ et $e,d$ entiers $\geq 2$; si $F$ contient $\boldsymbol{\mu}_d$ et $\boldsymbol{\mu}_e$, cette courbe a beaucoup d'automorphismes, ce qui permet de factoriser la fonction $\zeta_X(s)$ sous la forme $\frac{\zeta_F(s)\zeta_F(s-1)}{L(\chi_1,s)\cdots L(\chi_{2g},s)}$ où $g$ est le genre de $X$ et les $\chi_i$ sont des caractères de Hecke de $F$ pour lesquels le prolongement analytique de la fonction $L$ est connu depuis les travaux de HECKE [132].

• Si $d = 2$ et $e = 3, 4$ dans l'exemple précédent, on tombe sur une courbe elliptique avec un automorphisme d'ordre 3 ou 4; cette courbe a donc de la multiplication complexe par $\mathbf{Z}[\frac{-1+\sqrt{-3}}{2}]$ ou $\mathbf{Z}[\sqrt{-1}]$. DEURING [80] a étendu le résultat à toutes les courbes elliptiques à multiplication complexe, et SHIMURA, TANIYAMA [200] aux variétés abéliennes à multiplication complexe de dimension quelconque.

• EICHLER [87] et SHIMURA [199] ont prouvé la conjecture pour les quotients de

courbes modulaires. Si $X$ est un tel quotient, ils ont prouvé que $\zeta(X,s) \sim \frac{\zeta(s)\zeta(s-1)}{L(f_1,s)\cdots L(f_g,s)}$ où $g$ est le genre de $X$ et $f_1,\ldots,f_g$ sont des formes modulaires de poids 2 (et il est élémentaire que $L(f,s)$ admet un prolongement holomorphe à $\mathbf{C}$ tout entier, si $f$ est une forme modulaire).

Les résultats ci-dessus datent des années 1950. Les second et troisième exemples sont « de type $\mathbf{GL}_1$ »; le dernier a les formes modulaires intégrées dans la construction. Si on part d'une situation qui ne relève *a priori* d'aucun de ces deux cas, on n'a pendant longtemps eu que des confirmations numériques partielles pour des courbes elliptiques particulières, définies sur $\mathbf{Q}$; la méthode de Faltings–Serre ([196], lettre du 26/10/1984) a permis de transformer ces calculs en des preuves[2] mais, pour un résultat général, il a fallu attendre les années 1990 et les travaux de WILES [222] menant à la preuve du théorème de Fermat.

THÉORÈME 1.5. — *Si $X$ est une courbe de genre 1 définie sur $\mathbf{Q}$, alors $X$ vérifie la conjecture de Hasse–Weil.*

*Remarque 1.6.* — (i) WILES [222] a prouvé le résultat pour les courbes semi-stables (incluant celles d'équation affine $y^2 = x(x-a)(x+b)$, avec $a,b \in \mathbf{Z} \setminus \{0\}$, premiers entre eux), le cas général a été démontré par BREUIL, CONRAD, DIAMOND, TAYLOR [33]. Ce que démontrent les auteurs ci-dessus est la conjecture de Taniyama–Weil selon laquelle $\zeta_X(s) \sim \frac{\zeta(s)\zeta(s-1)}{L(f,s)}$ où $f$ est une forme modulaire (i.e., $X$ est modulaire).

(ii) On déduit des résultats de TAYLOR [206, 208] certains cas de la conjecture pour des courbes de genre 2 définies sur $\mathbf{Q}$ (celles dont la jacobienne a un anneau d'endomorphismes plus gros que $\mathbf{Z}$). Le nouveau point de vue introduit dans [206, 208] a été le point de départ d'une flopée de beaux résultats dont on trouvera les plus emblématiques au n° 1.3.4, mais qui sont moins précis que [33] car ils ne prouvent qu'une version potentielle de l'analogue de la conjecture de Taniyama–Weil (i.e., $X$ est potentiellement modulaire).

Récemment, dans un petit article de 350 pages, BOXER, CALEGARI, GEE et PILLONI [28] ont démontré l'analogue en genre 2.

THÉORÈME 1.7. — (i) *Si $X$ est une courbe de genre 2 définie sur $\mathbf{Q}$ ou, plus généralement, sur un corps totalement réel, alors $X$ vérifie la conjecture de Hasse–Weil.*

(ii) *Si $X$ est une variété abélienne de dimension 2 définie sur $\mathbf{Q}$ ou, plus généralement, sur un corps totalement réel, alors $X$ vérifie la conjecture de Hasse–Weil.*

---

2. Par exemple, cela a permis à MESTRE [166] de prouver que la courbe elliptique $E$ d'équation $Y^2 + Y = X^3 - 7X + 6$ vérifie la conjecture de Taniyama-Weil; cette courbe est la courbe de conducteur minimal (son conducteur est 5077) avec $\mathrm{rg}(E(\mathbf{Q})) = 3$. La conjecture de Birch et Swinnerton-Dyer prédit que $L(E,s)$ a un zéro d'ordre 3 en $s = 1$ et GROSS, ZAGIER [123] venaient juste de donner un moyen de le vérifier, ce qui permettait d'utiliser un résultat antérieur de GOLDFELD [122] pour obtenir une minoration effective du nombre de classes des corps quadratiques imaginaires; l'utilisation de la courbe $E$ ci-dessus permettait par exemple de donner une liste *inconditionnelle* des corps quadratiques de nombre de classe 3 (cf. [175] pour les détails).

*Remarque 1.8.* — (i) Ce que démontrent les auteurs est un résultat de modularité potentielle (cf. th. 1.25) pour les variétés abéliennes de dimension 2 sur un corps totalement réel (même si on ne s'intéresse qu'aux courbes définies sur **Q**, on est forcé de traiter le cas d'un corps totalement réel général car la technique de modularité potentielle demande de faire une extension des scalaires à un corps totalement réel que l'on ne maîtrise pas). Le théorème s'en déduit par des méthodes qui sont standard depuis [206].

(ii) Il y a des raisons profondes pour lesquelles le cas des courbes de genre 2 est plus délicat que celui des courbes de genre 1 (cf. n° 1.3.5) ; ces mêmes raisons font que plus le genre augmente et plus le résultat est délicat à prouver (au moins avec les techniques ayant permis de prouver le cas des courbes de genres 1 et 2), et il semble que la preuve de la conjecture de Hasse–Weil pour les courbes de genre $\geq 3$ va demander des idées radicalement nouvelles.

Le résultat de modularité potentielle ci-dessus n'est pas la seule avancée importante de ces dernières années sur ce genre de questions (citons par exemple [2, 35, 27, 45, 171, 172, 173]), et nous allons essayer de faire un tour d'horizon de ces avancées.

Dans le reste de l'introduction, on commence (§ 1.2.1) par factoriser les fonctions zêta de Hasse–Weil en termes de fonctions $L$ (motiviques) et on énonce les conjectures concernant le comportement de ces fonctions $L$ (conj. 1.10 et 1.18). Le § 1.3 (principalement, les n°s 1.3.4, 1.3.5 et 1.3.6) est consacré aux énoncés des théorèmes de modularité mentionnés ci-dessus, et le § 1.4 donne des applications arithmétiques des théorèmes de modularité (conjectures de Sato–Tate et de Ramanujan, et équation de Fermat).

Le corps du texte est consacré à la présentation d'un certain nombre d'outils utilisés dans les preuves des théorèmes de modularité : le chap. 2 tourne autour des questions liées aux représentations galoisiennes d'un corps de nombres, le chap. 3 survole une partie de la théorie des représentations automorphes, et le chap. 4 est consacré aux représentations galoisiennes d'un corps local.

Ce texte reste à un niveau assez superficiel (nous nous sommes contenté d'énoncer les résultats et de présenter les outils intervenant dans les preuves — pas tous... — mais nous ne disons rien des preuves elles-mêmes), et nous conseillons au lecteur désireux d'approfondir le sujet de lire les articles de survol récents écrits par des gens plus compétents (par exemple ceux de Calegari [41], de[(3)] Gee [119] et de Caraiani, Shin [47] ou ceux, plus anciens, de Carayol [48] et Taylor [207]).

Je remercie Sophie Morel, Frank Calegari et Michael Harris pour leurs commentaires sur la version parue dans la brochure du séminaire.

**1.2. Fonctions $L$ motiviques**

*1.2.1. Factorisation des fonctions zêta de Hasse–Weil.* — En sus de la rationalité de la fonction zêta des variétés sur les corps finis, Weil [216] avait aussi conjecturé

---

3. Ou encore [92, § 3.1] pour une un point de vue sur le patching à la Taylor-Wiles-Kisin..., malheureusement absent de ce texte.

une factorisation des numérateur et dénominateur de cette fraction rationnelle. Ces conjectures sont aussi des théorèmes :

THÉORÈME 1.9. — *Si $X$ est une variété projective lisse connexe, de dimension $d$, définie sur un corps fini $\kappa$ de cardinal $q$, où $q = p^f$, alors :*

(i) (GROTHENDIECK [124]) $\zeta_X(s) = \frac{P_{X,1}(q^{-s})\cdots P_{X,2d-1}(q^{-s})}{P_{X,0}(q^{-s})\cdots P_{X,2d}(q^{-s})}$ *où $P_{X,i}(T) \in 1 + T\mathbf{Q}[T]$ est le déterminant de*[(4)] $1 - T\mathrm{Frob}_q^{-1}$ *agissant sur $H^i_{\mathrm{\acute{e}t}}(X_{\overline{\kappa}}, \mathbf{Q}_\ell)$ (pour n'importe quel $\ell \neq p$) ; en particulier, $P_{X,0}(T) = 1 - T$ et $P_{X,2d}(T) = 1 - q^d T$.*

(ii) (DELIGNE [74]) *Les racines de $P_{X,i}$ sont toutes de valeur absolue $q^{-i/2}$* (Hypothèse de Riemann sur les corps finis).

Soit maintenant $X$ une variété propre et lisse définie sur un corps de nombres $F$. On déduit du (i) du th. 1.9 une factorisation de la fonction zêta de Hasse–Weil de $X$ sous la forme
$$\zeta_X(s) = \prod_{i=0}^{2\dim X} L(h^i(X), s)^{(-1)^i}$$
où $h^i(X)$ est la famille des $H^i_{\mathrm{\acute{e}t}}(X_{\overline{\mathbf{Q}}}, \mathbf{Q}_p)$, où $p$ décrit l'ensemble des nombres premiers. Cette famille forme un système compatible (voir ci-dessous) de représentations galoisiennes (i.e. de représentations du groupe de Galois absolu $\mathrm{Gal}_F$ de $F$). On sait associer une fonction $L$ dite *motivique* à un tel système compatible et les fonctions $L$ motiviques ont des propriétés conjecturales bien meilleures que les fonctions zêta de Hasse–Weil (cf. conj. 1.18). Par exemple, si $X$ est une courbe projective lisse, alors $L(h^0(X), s) = \zeta_F(s)$, $L(h^2(X), s) = \zeta_F(s-1)$ et on conjecture que $L(h^1(X), s)$ admet un prolongement holomorphe à tout le plan complexe et une équation fonctionnelle reliant $L(h^1(X), s)$ et $L(h^1(X), 2-s)$.

*1.2.2. Corps de nombres.* — Introduisons un certain nombre de notations standard concernant les corps de nombres. Un corps de nombres $F$ est une extension finie de $\mathbf{Q}$. L'anneau des entiers $\mathscr{O}_F$ de $F$ est l'ensemble des $x \in F$ dont le polynôme minimal sur $\mathbf{Q}$ est à coefficients dans $\mathbf{Z}$.

Une *place* de $F$ est une classe d'équivalence de normes non triviales sur $F$. Si $v$ est une place de $F$, on note $F_v$ le complété de $F$ pour $v$.

Les places $v$ de $\mathbf{Q}$ ont des représentants naturels $|\ |_v$, où $|\ |_\infty$ est la valeur absolue et $|\ |_p$ est la norme $p$-adique normalisée par $|p|_p = p^{-1}$ si $p$ est un nombre premier. Si $w$ est une place de $F$, et $\|\ \|_w$ est une norme dans la classe de $w$, la restriction de $\|\ \|_w$ à $\mathbf{Q}$ définit une place $v$ de $\mathbf{Q}$, et on dit que *$w$ est au-dessus de $v$*, ce que l'on note $w \mid v$. On dit que *$w$ est finie* si $v$ est un nombre premier. On définit $|\ |_w$ sur $F_w$ par $|x|_w = |\mathrm{N}_{F_w/\mathbf{Q}_v}|_v$ : c'est une norme sur $F_w$ sauf si $F_w = \mathbf{C}$ (où c'est le carré de la valeur absolue), mais la normalisation est choisie pour que l'on ait *la formule du produit* : $\prod_v |x|_v = 1$, si $x \in F^*$.

---

4. $H^i_{\mathrm{\acute{e}t}}(X_{\overline{\kappa}}, \mathbf{Q}_\ell)$ est muni d'une action de $\mathrm{Gal}(\overline{\kappa}/\kappa)$ et $\mathrm{Frob}_q$ est l'élément de $\mathrm{Gal}(\overline{\kappa}/\kappa)$ agissant par $x \mapsto x^q$ (i.e. c'est le frobenius arithmétique).

Si $v \mid \infty$, on dit que $v$ est *réelle* si $F_v = \mathbf{R}$ et *complexe* si $F_v = \mathbf{C}$. On note $r_1$ (resp. $r_2$) le nombre de places réelles (resp. complexes) ; on a $r_1 + 2r_2 = [F : \mathbf{Q}]$. On dit que $F$ est *totalement réel* si $r_2 = 0$ ; on dit que $F$ est *un corps CM* si $F$ est totalement réel ou si c'est une extension quadratique totalement imaginaire (i.e. $r_1 = 0$) d'un corps totalement réel (auquel cas on dit que $F$ est *un corps CM imaginaire*). Si $\mathbf{Q}^{\text{CM}}$ (resp. $\mathbf{Q}^{\text{t-r}}$) désigne le composé de tous les corps CM (resp. totalement réels), alors $[\mathbf{Q}^{\text{CM}} : \mathbf{Q}^{\text{t-r}}] = 2$ (et donc $\mathbf{Q}^{\text{CM}} = \mathbf{Q}^{\text{t-r}}(\sqrt{-1})$) et l'élément non trivial de $\operatorname{Gal}(\mathbf{Q}^{\text{CM}}/\mathbf{Q}^{\text{t-r}})$ est la conjugaison complexe.

Si $v$ est finie, on note $\mathscr{O}_v$ l'anneau des entiers de $F_v$ et $\kappa_v$ son corps résiduel ; si $v \mid p$, $\kappa_v$ est une extension finie de $\mathbf{F}_p$ dont on note $f_v$ le degré et $|v|$ le cardinal (et donc $|p| = p$ et $|v| = p^{f_v}$). On note $I_v \subset \operatorname{Gal}_{F_v}$ le *sous-groupe d'inertie* (i.e. l'ensemble des $\sigma$ tels que $|\sigma(x) - x|_v < |x|_v$ pour tout $x \in \overline{F}_v^*$) ; on a alors une suite exacte

$$1 \mapsto I_v \to \operatorname{Gal}_{F_v} \to \operatorname{Gal}_{\kappa_v} \to 1$$

On note $\operatorname{Frob}_v \in \operatorname{Gal}_{\kappa_v}$ le *frobenius arithmétique* $x \mapsto x^{|v|}$.

On note $\mathbf{A}_F$ l'*anneau des adèles* de $F$ : c'est le produit restreint des $F_v$ relativement aux $\mathscr{O}_v$, i.e. l'ensemble des suites $(x_v)_v$ avec $x_v \in F_v$ pour toute $v$ et $x_v \in \mathscr{O}_v$ pour presque [5] toute $v$ finie (de manière plus compacte, $\mathbf{A}_F = (\mathbf{R} \otimes_{\mathbf{Q}} F) \times (\widehat{\mathbf{Z}} \otimes_{\mathbf{Z}} F)$ où $\widehat{\mathbf{Z}} = \prod_p \mathbf{Z}_p$ est le complété profini de $\mathbf{Z}$). On note $\mathbf{A}_F^{]\infty[}$ l'anneau des adèles finies (i.e. le produit restreint des $F_v$, pour $v$ finie) ; on a alors $\mathbf{A}_F = (\prod_{v \mid \infty} F_v) \times \mathbf{A}_F^{]\infty[}$, ce qui permet d'écrire $x = (x_v)_v \in \mathbf{A}_F$ sous la forme $(x_\infty, x^{]\infty[})$.

*1.2.3. Fonctions L d'Artin*. — Soit $F$ un corps de nombres. Si $\rho \colon \operatorname{Gal}_F \to \mathbf{GL}(V)$, où $V$ est un $\mathbf{C}$-espace vectoriel de dimension finie, est une représentation continue, alors, pour des raisons topologiques [6], $\rho$ se factorise à travers un quotient fini. Il s'ensuit que la restriction de $\rho$ à $\operatorname{Gal}_{F_v}$ est non ramifiée (i.e. le sous-groupe d'inertie $I_v$ agit trivialement : $V^{I_v} = V$) pour tout $v$ en dehors d'un ensemble fini $S$.

On associe à $\rho$ une fonction $L$ (dite *fonction L d'Artin*) par la formule suivante dans laquelle $E_v[T]$ est le déterminant de $1 - T\operatorname{Frob}_v^{-1}$ agissant sur $V^{I_v}$ :

$$L(\rho, s) = \prod_v L_v(\rho, s), \quad \text{où } L_v(\rho, v) = E_v(|v|^{-s})^{-1}$$

Par construction, $E_v \in 1 + T\mathbf{C}[T]$, et $E_v$ est de degré $\leq n$ avec égalité pour tout $v \notin S$. Le produit converge absolument dans le demi-plan $\operatorname{Re}(s) > 1$ et y définit une fonction holomorphe. La conjecture suivante, vieille d'un siècle, est le prototype de toutes les conjectures sur les fonctions $L$.

---

5. I.e. à l'exception d'un nombre fini de telles $v$.

6. Les groupes de Galois sont des groupes profinis (et munis de la topologie de groupe profini), et un voisinage $U$ assez petit de $1 \in \mathbf{GL}_d(V)$ ne contient pas de sous-groupe autre que $\{1\}$ ; or l'image inverse de $U$ par $\rho$ contient un voisinage ouvert de $1 \in \operatorname{Gal}_F$, et donc contient un sous-groupe d'indice fini qui s'envoie sur $\{1\}$ d'après ce qui précède.

Conjecture 1.10. — (Artin [4]). *La fonction $L(\rho,s)$ possède un prolongement analytique à $\mathbf{C}$, holomorphe en dehors d'un pôle éventuel en $s=1$ d'ordre $\dim_{\mathbf{C}} V^{\mathrm{Gal}_F}$, et ne s'annule pas sur la droite $\mathrm{Re}(s)=1$.*

*De plus, il existe un entier $N$ divisible uniquement par les nombres premiers au-dessus des éléments de $S$, et une constante $\varepsilon \in \mathbf{C}$, tels que l'on ait l'équation fonctionnelle*

$$\Bigl(\prod_{v|\infty}\Gamma(\rho_v,s)\Bigr)L(\rho,s) = \varepsilon N^{-s}\Bigl(\prod_{v|\infty}\Gamma(\rho_v,1-s)\Bigr)L(\check\rho,1-s)$$

*où $\check\rho$ est la représentation duale, $\Gamma(\rho_v,s) = \Gamma_{\mathbf{C}}(s)^n$ si $v$ est une place complexe*[7]*, et $\Gamma(\rho_v,s) = \Gamma_{\mathbf{R}}(s)^{n_v^+}\Gamma_{\mathbf{R}}(s+1)^{n_v^-}$ avec $n_v^\pm = \dim V^{\mathrm{Frob}_v=\pm 1}$ (et $\mathrm{Frob}_v$ est la conjugaison complexe) si $v$ est une place réelle.*

• Si $F=\mathbf{Q}$ et $\dim V=1$, il résulte du théorème de Kronecker–Weber (selon lequel toute extension galoisienne de $\mathbf{Q}$, de groupe de Galois abélien, est incluse dans l'extension cyclotomique) que $L(\rho,s)$ est aussi une fonction $L$ associée à un caractère de Dirichlet ; on en déduit la conjecture dans ce cas.

• Si $F$ est quelconque et $\dim V=1$, la conjecture résulte de la théorie globale du corps de classes qui fournit une identité $L(\rho,s) = L(\chi,s)$ où $\chi$ est un caractère localement constant de $\mathbf{A}_F^*/F^*$, et l'analogue pour $L(\chi,s)$ de la conj. 1.10 a été démontré par Hecke [132] (la preuve de Hecke a été adélisée par Tate [203]).

• Artin avait prouvé qu'une puissance de $L(\rho,s)$ est méromorphe et admet une équation fonctionnelle comme ci-dessus. Brauer a prouvé que $L(\rho,s)$ admet un prolongement méromorphe et une équation fonctionnelle en utilisant le résultat suivant qui fournit une identité du type $L(\rho,s) = \prod_{i\in I} L(\chi_i,s)^{n_i}$, où les $L(\chi_i,s)$ sont holomorphes grâce à la théorie du corps de classe et aux résultats de Hecke [132] déjà mentionnés.

Théorème 1.11. — (Brauer [30]). *Soient $G$ un groupe fini et $\rho$ une représentation de $G$ de dimension finie sur un corps algébriquement clos de caractéristique $0$. Alors on peut écrire $\rho$, vue comme représentation virtuelle, sous la forme $\rho = \sum_{i\in I} n_i \mathrm{Ind}_{H_i}^G \chi_i$ où, si $i\in I$, $H_i$ est un sous-groupe résoluble de $G$, $\chi_i$ est une représentation de dimension $1$ de $H_i$, et $n_i \in \mathbf{Z}$.*

• On ne peut pas, en général, faire en sorte que les $n_i$ soient $\geq 0$ (ce qui prouverait l'holomorphie) et, depuis l'avancée de Brauer, les progrès ont été très limités (ils sont résumés dans le premier point du n° 1.3.6).

*1.2.4. Généralités sur les fonctions L*

*Définition.* — Soit $F$ un corps de nombres. Une fonction $[L]$ de degré $n$ sur $F$ est une classe d'équivalence de produit eulériens $\prod_v L_v(s)$, où $v$ décrit les places finies de $F$, $L_v(s) = E_v(|v|^{-s})^{-1}$ et $E_v \in 1+T\mathbf{C}[T]$ est de degré $n$ pour toute $v$ sauf un nombre fini, deux tels produits étant dits équivalents si $L_v = L_v'$ pour toute $v$ sauf un nombre fini.

Si $E_v$ est de degré $n$, on le factorise sous la forme $\prod_{i=1}^n (1-\alpha_{v,i}T)$. On dit que $[L]$ est *pure de poids $w$* si $|\alpha_{v,i}| = |v|^{-w/2}$ pour tous $v$ et $i$ sauf pour un nombre fini.

---

7. $\Gamma_{\mathbf{R}}$ et $\Gamma_{\mathbf{C}}$ sont définies dans la note 1.

Les fonctions [$L$] intéressantes ont en général un représentant privilégié mais la définition des facteurs d'Euler manquants est souvent subtile (voire conjecturale) : un tel représentant est une fonction $L$. Par exemple, les fonctions $L$ d'Artin sont les représentants privilégiés de leur classe d'équivalence.

*Opérations sur les fonctions* [$L$]. — Les opérations naturelles sur les représentations d'Artin admettent des traductions en termes de fonctions [$L$] :

• Si [$L^{(1)}$] et [$L^{(2)}$] sont de degrés $n$ et $m$, on définit [$L^{(1)} \oplus L^{(2)}$] (de degré $n+m$) et [$L^{(1)} \otimes L^{(2)}$] (de degré $nm$) par

$$[L^{(1)} \oplus L^{(2)}] = [L^{(1)}L^{(2)}] \quad \text{et} \quad (L^{(1)} \otimes L^{(2)})_v = \prod_{i=1}^n \prod_{j=1}^m (1 - \alpha_{v,i}^{(1)}\alpha_{v,j}^{(2)}|v|^{-s})^{-1}$$

• On définit la duale [$\check{L}$] de [$L$] par $\check{L}_v = \prod_{i=1}^n (1 - \alpha_{v,i}^{-1}|v|^{-s})^{-1}$

• Si [$L$] est sur $K$ (resp. sur $F$) de degré $n$, on définit son induite [$\mathrm{I}_K^F L$] à $F$ (resp. sa restriction [$\mathrm{R}_F^K L$] à $K$), qui est de degré $[K:F]n$ (resp. $n$), par

$$(\mathrm{I}_K^F L)_v = \prod_{w|v} L_w \text{ (donc } [\mathrm{I}_K^F L] = [L]) \quad \text{et} \quad (\mathrm{R}_F^K L)_w = \prod_{i=1}^n (1 - \alpha_{v,i}^{[K_w:F_v]} \mathrm{N}w^{-s})^{-1}$$

• Si [$L$] est de degré $n$, et si $r \colon \mathbf{GL}_n(\mathbf{C}) \to \mathbf{GL}_N(\mathbf{C})$ est un morphisme de groupes de Lie (un tel morphisme est automatiquement induit par un morphisme de groupes algébriques), on peut définir [$r \circ L$] (de degré $N$) : par exemple, si $n = 2$, on dispose de $\mathrm{Sym}^k \colon \mathbf{GL}_2(\mathbf{C}) \to \mathbf{GL}_{k+1}(\mathbf{C})$ (induite par l'action de $\mathbf{GL}_2(\mathbf{C})$ sur la puissance symétrique $k$-ième de l'action standard de $\mathbf{GL}_2(\mathbf{C})$ sur $\mathbf{C}^2$), et on a

$$(\mathrm{Sym}^k \circ L)_v = \prod_{i=0}^k (1 - \alpha_{v,1}^i \alpha_{v,2}^{k-i}|v|^{-s})^{-1}$$

*Prolongement analytique et équation fonctionnelle.* — On dit qu'une fonction [$L$], de degré $n$, admet :

• un prolongement analytique si un (ou, de manière équivalente, tout) représentant peut se prolonger en une fonction méromorphe sur $\mathbf{C}$ tout entier,

• une équation fonctionnelle si [$L$] a un représentant privilégié $L$ tel qu'il existe $w, a_1, \ldots, a_n \in \mathbf{C}$, une constante $\varepsilon$, et un entier $N$, tels que, pour tout $s \in \mathbf{C}$,

$$\Big(\prod_i \Gamma_{\mathbf{R}}(s - a_i)\Big)L(s) = \varepsilon N^s \Big(\prod_i \Gamma_{\mathbf{R}}(w + 1 - s - a_i)\Big)\check{L}(w + 1 - s)$$

Le th. 1.11 a comme conséquence l'énoncé suivant, très utile pour prouver la méromorphie de fonctions [$L$] potentiellement modulaires.

PROPOSITION 1.12. — *Soit* [$L$] *sur* $F$. *On suppose qu'il existe une extension galoisienne* $K/F$, *finie, telle que, pour tout* $F \subset K' \subset K$ *avec* $\mathrm{Gal}(K/K')$ *résoluble, et toute représentation* $\chi$ *de dimension* 1 *de* $\mathrm{Gal}(K/K')$, *la fonction* [$(\mathrm{R}_F^{K'} L) \otimes \chi$] *admette un prolongement analytique et une équation fonctionnelle. Alors* [$L$] *admet un prolongement analytique et une équation fonctionnelle.*

*1.2.5. Systèmes compatibles et conjecture de Fontaine–Mazur.* — Il est commode de considérer des représentations de $\mathrm{Gal}_F$ sur des $\overline{\mathbf{Q}}_p$-espaces vectoriels plutôt que des $\mathbf{Q}_p$-espaces vectoriels [8]. Cela conduit, pour définir la notion de système compatible, à fixer des plongements

$$\iota_\infty \colon \overline{\mathbf{Q}} \hookrightarrow \mathbf{C} \quad \text{et} \quad \iota_p \colon \overline{\mathbf{Q}} \hookrightarrow \overline{\mathbf{Q}}_p, \text{ pour tout } p.$$

*Définition 1.13*. — Un *F-système faiblement compatible* de dimension $d$ est une famille $(M_p)_p$, où $p$ décrit l'ensemble des nombres premiers, et :

• $M_p$ est un $\overline{\mathbf{Q}}_p$-module de rang $d$, muni d'une action linéaire continue de $\mathrm{Gal}_F$,

• Il existe $S$ fini, tel que $M_p$ soit non ramifiée en dehors de $S \cup S_p$, où $S_p = \{v,\ v \mid p\}$.

• Si $v \notin S \cup S_p$, alors le déterminant $E_v(T)$ de $1 - T\mathrm{Frob}_v^{-1}$ agissant sur $M_p$ est à coefficients dans $\overline{\mathbf{Q}}$ et ne dépend pas de $p$.

On dit que $(M_p)_p$ est :

• *irréductible* si chaque $M_p$ est irréductible,

• *pur de poids $w$*, si les racines $\alpha$ de $E_v$ vérifient $|\sigma(\alpha)|_\infty = |v|^{-w/2}$ et $|\sigma(\alpha)|_\ell = 1$ pour tous $\sigma \in \mathrm{Gal}_{\mathbf{Q}}$ et $\ell$ avec $v \nmid \ell$.

• *compatible* si :

 ⋄ pour tout $p$ et tout $v \mid p$, la restriction de $M_p$ à $\mathrm{Gal}_{F_v}$ est [9] de Rham, et si, pour tout $\tau \in \mathrm{Hom}(F, \overline{\mathbf{Q}})$, les $\tau$-poids de Hodge–Tate (cf. rem. 4.7) de $M_p$ forment un multi-ensemble $H_\tau$ ne dépendant pas de $p$ (et on dit que $(M_p)_p$ est *régulier* si les éléments de $H_\tau$ ont tous multiplicité 1, pour tout $\tau$),

 ⋄ pour tout $p \notin S$ et tout $v \mid p$, la restriction de $M_p$ à $\mathrm{Gal}_{F_v}$ est cristalline et le déterminant de $1 - T\varphi^{f_v}$ agissant sur $\mathbf{D}_{\mathrm{cris}}(M_p)$ est $E_v(T)$.

La plupart du temps nous allons fixer une base de $M_p$ sur $\overline{\mathbf{Q}}_p$, ce qui transforme $M_p$ en un morphisme continu $\rho_p \colon \mathrm{Gal}_F \to \mathbf{GL}_d(\overline{\mathbf{Q}}_p)$, et nous parlerons d'un système (faiblement) compatible $(\rho_p)_p$ plutôt que $(M_p)_p$.

*Exemple 1.14*. — Soit $\rho \colon \mathrm{Gal}_F \to \mathbf{GL}_d(\mathbf{C})$ une représentation continue (et donc d'image finie puisque $\mathrm{Gal}_F$ est profini). Alors $\mathrm{Tr}(\rho(\sigma)) \in \overline{\mathbf{Q}}$ pour tout $\sigma \in \mathrm{Gal}_F$ et on peut conjuguer $\rho$ de telle sorte que $\rho$ soit à valeurs dans $\mathbf{GL}_d(\overline{\mathbf{Q}})$, ce qui permet de transformer $\rho$ en un système faiblement compatible. Un tel système, dit *système d'Artin*, est compatible (et même fortement compatible, voir ci-dessous), pur de poids 0, et de poids de Hodge–Tate tous nuls.

---

8. Notons que, si $\rho \colon \mathrm{Gal}_F \to \mathbf{GL}_d(\overline{\mathbf{Q}}_p)$ est un morphisme continu de groupes topologiques, il résulte du théorème de Baire qu'il existe une extension finie $L$ de $\mathbf{Q}_p$ telle que l'image de $\rho$ soit incluse dans $\mathbf{GL}_d(L)$, et donc qu'une $\overline{\mathbf{Q}}_p$-représentation de $\mathrm{Gal}_F$ s'obtient par extension des scalaires à partir d'une représentation sur une extension finie de $\mathbf{Q}_p$.

9. Voir §4.2, en particulier n° 4.2.2 et le début de n° 4.2.3 pour la définition des termes « de Rham », « cristalline », et les notations $\mathbf{D}_{\mathrm{cris}}$, $\mathbf{D}_{\mathrm{dR}}$, $\varphi$.

*Remarque 1.15*. — (i) Si $X$ est une variété propre et lisse définie sur $F$, les $H^i(X_{\overline{\mathbf{Q}}}, \overline{\mathbf{Q}}_p)$ forment un système compatible, pur de poids $i$, et la multiplicité de $j$ comme $\tau$-poids de Hodge–Tate est dim $H^{i-j}(X, \Omega^j)$ pour tout $\tau \in \mathrm{Hom}(F, \overline{\mathbf{Q}})$. (La faible compatibilité et la pureté résultent du th. 1.9, la compatibilité et la multiplicité des poids de Hodge–Tate résulte des théorèmes de comparaison $p$-adiques (cf. th. 4.1).)

(ii) Dans ce cas, si $v \in S$, on a aussi une compatibilité (conjecturale, en général) entre ce qui se passe pour la restriction de $M_p$ à $\mathrm{Gal}_{F_v}$ pour $v \mid p$ et pour $v \nmid p$, et on peut introduire la notion de *forte compatibilité* pour tenir compte de cette compatibilité supplémentaire.

(iii) On voit mal comment construire des systèmes faiblement compatibles qui ne soient pas d'origine géométrique ; en particulier *il est naturel de penser qu'un système faiblement compatible est compatible (et même fortement compatible)*. La conjecture ci-dessous est nettement plus surprenante mais s'est révélée extrèmement fructueuse.

CONJECTURE 1.16. — (FONTAINE–MAZUR [101, 109]) *Soit $M_p$ une $\mathbf{Q}_p$-représentation irréductible de $\mathrm{Gal}_F$, non ramifiée en dehors de $S$ fini, dont la restriction à $\mathrm{Gal}_{F_v}$ est de Rham pour tout $v \mid p$. Alors $M_p$ provient de la géométrie, i.e. $M_p$ est un sous-quotient*[10] *d'un $H^i_{\mathrm{\acute{e}t}}(X_{\overline{\mathbf{Q}}}, \mathbf{Q}_p(j))$ où $X$ est une variété propre et lisse définie sur $F$, $i \geq 0$ et $j \in \mathbf{Z}$.*

*Remarque 1.17*. — (i) Une représentation vérifiant les conditions de la conjecture est dite *géométrique*. Une conséquence de la conjecture est qu'une représentation géométrique peut s'étendre en un système compatible (un analogue arithmétique de la conjecture des compagnons de DELIGNE [75]).

(ii) Un cas particulier de la conjecture est celui où $M_p$ est de Rham à poids de Hodge–Tate tous nuls. Dans ce cas, la conjecture implique que $\mathrm{Gal}_F$ agit à travers un quotient fini, ce qui évoque une conjecture géométrique de DE JONG [72].

*1.2.6. Fonctions L motiviques*. — Si $\rho = (\rho_p)_p$ est un système compatible, non ramifié en dehors de $S$, soit $[L(\rho, s)]$ la classe de $L^{(S)}(\rho, s) := \prod_{v \notin S} E_v(|v|^{-s})^{-1}$.

Si $(\rho_p)_p$ est pur de poids $w$, le produit converge sur le demi-plan $\mathrm{Re}(s) > 1 + \frac{w}{2}$ et y définit une fonction holomorphe.

CONJECTURE 1.18. — *Si $\rho = (\rho_p)_p$ est un système compatible, pur de poids $w$, alors $[L(\rho, s)]$ contient un représentant privilégié ayant un prolongement holomorphe*[11] *et une équation fonctionnelle ; de plus ce prolongement ne s'annule pas sur $\mathrm{Re}(s) = 1 + \frac{w}{2}$.*

*Remarque 1.19*. — (i) SERRE [194] a donné des formules précises pour les facteurs d'Euler $E_v(|v|^{-s})$ manquants et les $a_i$, $N$ et $\varepsilon$ de l'équation fonctionnelle. Les $a_i$ dépendent des poids de Hodge–Tate et de l'action de la conjugaison complexe, les $E_v$ sont définis par $E_v(T) = \det(1 - T\mathrm{Frob}_v^{-1})$ agissant sur le sous-espace de $\rho_p$ fixé par $I_v$. Cela donne *a priori* un polynôme à coefficient dans $\overline{\mathbf{Q}}_p$ et il n'est pas clair que ce polynôme soit indépendant de $p$ (pour les systèmes de la forme $h^i(X)$ la question de cette indépendance se pose depuis 50 ans ; on sait la prouver pour les $h^1(X)$).

---

10. $\mathbf{Q}_p(j)$ désigne l'espace $\mathbf{Q}_p$ sur lequel $\mathrm{Gal}_F$ agit par la puissance $j$-ième du caractère cyclotomique.
11. On se permet un nombre fini de pôles d'ordre fini pour ne pas avoir à traiter différemment les facteurs de la forme $\zeta(s-n)$.

(ii) DENINGER [77] a proposé un autre formalisme pour définir ces facteurs locaux dans le cas des $h^i(X)$.

**1.3. Modularité de fonctions $L$**

Il est possible que l'on dispose un jour d'une cohomologie arithmétique expliquant les bonnes propriétés des fonctions $L$ motiviques (comme envisagé dans [78, 79]) mais, en l'état actuel des choses, la preuve de l'existence de prolongements analytiques de fonctions $L$ motiviques demande de relier ces fonctions $L$ à des fonctions $L$ modulaires (ou au moins automorphes) dont les fonctions $L$ de formes modulaires sont le prototype.

*1.3.1. Modularité et modularité potentielle.* — Une fonction $[L]$ sur $F$ est dite :

• *modulaire* (cf. n° 3.1.3) si c'est la classe de $L(\pi, s)$ où $\pi$ est une représentation automorphe [12] de $\mathbf{GL}_n(\mathbf{A}_F)$,

• *potentiellement modulaire* s'il existe $K$, galoisienne sur $\mathbf{Q}$, telle que $[\mathrm{R}_F^K L]$ soit modulaire,

• *automorphe* (cf. n° 3.3.2) si c'est la classe de $L(\pi, r, s)$ où $\pi$ est une représentation automorphe d'un $F$-groupe réductif $\mathbb{G}$ et $r$ est une représentation de son dual de Langlands ${}^L\mathbb{G}$.

Si $X$ est un objet possédant une fonction $[L]$ (notée $[L(X, s)]$), on dit que $X$ est (potentiellement) modulaire si $[L(X, s)]$ l'est. Pour une fonction $[L]$ motivique, on a aussi une notion de *modularité forte* (cf. n° 1.3.2 ci-dessous).

*Remarque 1.20.* — (i) Les fonctions $L$ motiviques de degré 1 sont modulaires (traduction de la théorie globale du corps de classes) et même fortement modulaires ; celles associées aux formes modulaires aussi.

(ii) On *conjecture que les fonctions $[L]$ motiviques et automorphes sont modulaires*. Compte-tenu des remarquables propriétés des fonctions $[L]$ modulaires (th. 3.3), la modularité d'une fonction $[L]$-motivique implique la conj. 1.18 pour cette fonction $[L]$.

(iii) La modularité potentielle implique [13] aussi la conj. 1.18 à ceci près que l'on

---

12. Si $L_v = E_v(|v|^{-s})^{-1}$ avec $E_v$ de degré $n$, on sait associer à $L_v$ une représentation $\Pi_v$ de $\mathbf{GL}_n(F_v)$, non ramifiée (i.e. possédant un vecteur $e_v \neq 0$ fixe par $\mathbf{GL}_n(\mathscr{O}_v)$) ; $[L]$ est modulaire si et seulement si il existe une représentation automorphe $\pi = \otimes'_v \pi_v$ de $\mathbf{GL}_n(\mathbf{A}_F)$ telle que l'on ait $\pi_v = \Pi_v$ pour presque tout $v$ (le symbole $\otimes'$ désigne le produit tensoriel restreint : on considère des tenseurs du type $\otimes_v \phi_v$ avec $\phi_v = e_v$ pour presque tout $v$ ; le groupe $\mathbf{GL}_n(\mathbf{A}_F)$ est le produit restreint des $\mathbf{GL}_n(F_v)$ relativement aux $\mathbf{GL}_n(\mathscr{O}_v)$, et il agit continûment sur le produit tensoriel restreint car $e_v$ est fixe par $\mathbf{GL}_d(\mathscr{O}_v)$) ; cette représentation $\pi$, si elle existe, est unique d'après le théorème de multiplicité 1 (cf. th. 3.5).

13. Ce n'est pas évident mais cela découle de la prop. 1.12 et des travaux d'ARTHUR, CLOZEL [5] sur le changement de base mentionnés dans le n° 3.3.5 : si $[\mathrm{R}_F^K L]$ est modulaire, alors $[\mathrm{R}_F^{K'} L \otimes \chi]$ est modulaire pour tout $K' \subset K$ avec $\mathrm{Gal}(K/K')$ résoluble et toute représentation $\chi$ de dimension 1 de $\mathrm{Gal}(K'/K)$.

n'obtient qu'un prolongement méromorphe au lieu d'un prolongement holomorphe ; comme le montre la conjecture d'Artin, le fossé entre la méromorphie et l'holomorphie est le plus souvent infranchissable mais la méromorphie est déjà un énoncé remarquable.

(iv) Si $[L^{(1)}]$ et $[L^{(2)}]$ sont modulaires, alors $[L^{(1)} \otimes L^{(2)}]$ est automorphe mais on ne sait pas, en général, prouver que $[L^{(1)} \otimes L^{(2)}]$ est modulaire ; cependant on sait prouver que $[L^{(1)} \otimes L^{(2)}]$ a un prolongement holomorphe et une équation fonctionnelle (cf. th. 3.6). Ceci permet de prouver la conj. 1.18 pour certaines fonctions $[L]$ motiviques sans prouver leur modularité (mais cela utilise quand même leur automorphie).

La décennie 1995-2005 a vu l'éclosion d'un corpus de résultats (dont [5, 61, 62, 63, 64, 147, 148, 149]) du type « automorphe ⇒ modulaire » mais on est loin d'avoir résolu tous les problèmes. Ces résultats sont utilisés pleinement dans la preuve de l'implication « potentiellement modulaire ⇒ prolongement méromorphe et équation fonctionnelle ».

Il y a eu une moisson récente de résultats du type « motivique ⇒ potentiellement modulaire » ou même « motivique ⇒ modulaire » ; on déduit de certains de ces résultats des implications « automorphe ⇒ modulaire » inaccessibles pour le moment par des méthodes purement automorphes.

*1.3.2. Modularité forte.* — L'énoncé optimal est « motivique ⇒ *fortement modulaire* » qui inclut la compatibilité avec la correspondance de Langlands locale (HARRIS, TAYLOR [129], HENNIART [135], SCHOLZE [189] pour $\mathbf{GL}_n$, et FARGUES, SCHOLZE [98] pour un groupe général) y compris aux mauvaises places (pour lesquelles le facteur d'Euler de la fonction $L$ est très loin de décrire les objets locaux correspondants). Plus précisément :

• Si $\rho_p \colon \mathrm{Gal}_F \to \mathbf{GL}_n(\overline{\mathbf{Q}}_p)$ est géométrique (au sens du (i) de la rem. 1.17), on sait associer à $\rho_p$, pour tout $v$, une représentation $\mathrm{WD}_v(\rho_p)$ de dimension $n$ du groupe de Weil–Deligne de $F_v$ (cf. (iii) de la rem. 4.6).

• La correspondance de Langlands locale associe une représentation $\pi_v(\mathrm{WD}_v)$ de $\mathbf{GL}_n(F_v)$ à une représentation $\mathrm{WD}_v$ de dimension $n$ du groupe de Weil–Deligne[14] de $F_v$.

• On dit que $\rho_p$ est *fortement modulaire* s'il existe une représentation automorphe[15] $\pi(\rho) = \otimes'_v \pi_v(\rho)$ de $\mathbf{GL}_n(\mathbf{A}_F)$ telle que $\pi_v(\rho) = \pi_v(\mathrm{WD}_v(\rho))$ pour toute[16] $v$.

*Remarque 1.21.* — (i) Si $\rho = (\rho_p)_p$ est fortement compatible, la représentation $\mathrm{WD}_v(\rho_p)$ ne dépend pas de $p$, et donc « $\rho_p$ fortement modulaire pour au moins un $p$ » ⇔ « $\rho_p$ fortement modulaire pour tout $p$ ». Si cette condition est vérifiée on dit que $\rho$ est fortement modulaire, et on note $\pi(\rho)$ la représentation $\pi(\rho_p)$ pour n'importe quel $p$.

---

14. Si $v$ est une place finie, on note $\mathrm{W}_{F_v} \subset \mathrm{Gal}_{F_v}$, le groupe de Weil de $F_v$ ($\sigma \in \mathrm{W}_{F_v}$ si et seulement si son image dans $\mathrm{Gal}_{\kappa_v}$ est $x \mapsto x^{p^{\deg \sigma}}$, avec $\deg \sigma \in \mathbf{Z}$) ; une représentation du groupe de Weil–Deligne de $F_v$ est un couple $(\rho, N)$ où $\rho$ une représentation lisse de $\mathrm{W}_K$ et $N$ un opérateur vérifiant $N\rho(\sigma) = p^{-\deg \sigma} \rho(\sigma) N$, pour tout $\sigma \in \mathrm{W}_{F_v}$.

15. Voir la note 12 pour la notation $\otimes'$.

16. Si $v \mid \infty$, cela se traduit par le fait que $\pi_v(\rho)$ encode les $\tau$-poids de Hodge–Tate de $\rho_p$, pour les $\tau \in \mathrm{Hom}(F, \overline{\mathbf{Q}})$ induisant $v$, plus l'action de la conjugaison complexe si $v$ est réelle.

(ii) On a « $\rho$ fortement modulaire $\Rightarrow \rho$ modulaire » : en effet, $L_v(\pi(\rho_p), s) = L_v(\rho_p, s)$ pour tout $v \nmid p$ tel que $\rho_p$ est non ramifiée en $v$, et donc $[L(\pi(\rho), s)] = [L(\rho, s)]$.

Conjecture 1.22. — (Fontaine–Mazur + Langlands)
*Si $\rho \colon \mathrm{Gal}_F \to \mathbf{GL}_n(\overline{\mathbf{Q}}_p)$ est géométrique, alors $\rho$ est fortement modulaire.*

La différence avec la conj. 1.16 est que l'on ne prétend pas que $\rho$ provient de la géométrie, et il y a pas mal de cas où on sait prouver la conj. 1.22 (sous la forme plus faible « géométrique $\Rightarrow$ modulaire », suffisante pour les applications aux fonctions $L$) mais pas la conj. 1.16.

Les conj. 1.22 et 1.16 sont connues pour $n = 1$ (théorie de la multiplication complexe), et pour les représentations impaires de dimension 2 de $\mathrm{Gal}_{\mathbf{Q}}$ (combinaison de travaux d'Emerton [89], Kisin [154] et Pan [176] pour la conj. 1.22, auxquels il faut rajouter la contribution de Scholl [188] pour la conj. 1.16). Tous les théorèmes de relèvement modulaire du n° 2.3.2 démontrent des cas particuliers de la conj. 1.22 (sous la forme faible, en général).

*Remarque 1.23.* — Soit $\rho = (\rho_p)_p$ un système compatible, modulaire, de dimension $n$, et soit $\pi(\rho) = \otimes'_v \pi_v(\rho)$ la représentation automorphe de $\mathbf{GL}_n(\mathbf{A}_F)$ associée. Par définition, $L(\rho, s) \sim L(\pi(\rho), s)$ ce qui se traduit par $\pi_v(\rho) = \pi_v(\mathrm{WD}_v(\rho))$ pour tout $v \notin S(\rho)$, où $S(\rho)$ est fini. Or $\rho$ est complètement déterminée par les $\pi_v(\mathrm{WD}_v(\rho))$ pour $v \notin S(\rho)$ d'après le théorème de Čebotarev, et $\pi(\rho)$ est complètement déterminée par les $\pi_v(\rho)$, pour $v \notin S(\rho)$ d'après le théorème de multiplicité 1 (cf. th. 3.5). Il est donc envisageable que l'on puisse prouver un résultat du type « modulaire $\Rightarrow$ fortement modulaire » pour un système compatible [17].

*1.3.3. Une caricature des preuves de modularité.* — Les preuves de résultats de modularité semblent toute coulées dans un même moule mais les résultats de ces dernières années (pour $\ell_0 \geq 1$, voir ci-dessous) relevaient de l'utopie il y a 10 ans.

- *Théorèmes de relèvement modulaire.* — On s'intéresse à l'implication « motivique $\Rightarrow$ modulaire » mais la seule stratégie dont on dispose pour prouver un tel énoncé est de disposer d'une construction inverse « modulaire + conditions [18] $\Rightarrow$ motivique », et de vérifier que le système compatible $\rho = (\rho_p)_p$ dont on veut prouver la modularité est dans l'image.

Les conditions mises sur $\pi$ fournissent un sous-groupe algébrique minimal $\mathbb{G}$ de $\mathbf{GL}_n$ « défini sur $\mathscr{O}_F$ », tel que $\rho_p(\mathrm{Gal}_F) \subset \mathbb{G}(\overline{\mathbf{Z}}_p)$ pour tout $p$ (typiquement, si on fixe le caractère central de $\pi$, cela fixe le déterminant de $\rho_\pi$ et si on impose une autodualité à $\pi$, cela force une autodualité de $\rho_\pi$ qui se traduit par le fait que $\mathbb{G} \subset \mathbf{GSp}_n$ ou $\mathbf{GO}_n$).

Pour prouver que $\rho$ est dans l'image il suffit de prouver que $\rho_p$ est dans l'image,

---

17. Un tel résultat serait l'analogue de la « conjecture $\varepsilon$ » de Serre [195] pour les représentations modulo $p$ de $\mathrm{Gal}_{\mathbf{Q}}$.

18. Une condition minimum à imposer est que $\pi$ soit algébrique ; il vaut mieux fixer le caractère central et on peut aussi vouloir que $\pi$ soit autoduale, etc.

pour un $p$ bien choisi. Pour cela, la stratégie est, depuis Wiles [222], de prouver que la réduction $\overline{\rho}_p$ de $\rho_p$ modulo $p$ est dans l'image (i.e. que $\overline{\rho}_p$ est modulaire), et de prouver que tous les relèvements $\rho'\colon \mathrm{Gal}_F \to \mathbb{G}(\overline{\mathbf{Z}}_p)$ de $\overline{\rho}_p$ ayant les mêmes propriétés locales[19] que $\rho_p$ sont dans l'image (et donc que $\rho_p$ aussi). Autrement dit, on cherche à prouver un *théorème de relèvement modulaire.*

• *Modularité potentielle modulo $p$.* — Les théorèmes de relèvement modulaire sont de la forme suivante : « $\overline{\rho}$ modulaire[20] + $\overline{\rho}(\mathrm{Gal}_F)$ assez gros + conditions de nature "théorie de Hodge $p$-adique" $\Rightarrow \rho$ modulaire ».

La condition de grosse image intervient dans l'utilisation du théorème de Čebotarev pour construire des ensembles de nombres premiers (dits de Taylor–Wiles) avec des propriétés spéciales.

La condition la plus sérieuse est « $\overline{\rho}$ modulaire » qui n'est pas vraiment plus facile à vérifier « $\rho$ modulaire » (de fait, $\rho = (\rho_p)_p$ est faiblement compatible, et si on pouvait montrer a priori que $\overline{\rho}_p$ est modulaire pour suffisamment de $p$, on en déduirait que $\rho_p$ est modulaire pour tout $p$). Dans le cas le plus simple $\overline{\rho}\colon \mathrm{Gal}_{\mathbf{Q}} \to \mathbf{GL}_2(\overline{\mathbf{F}}_p)$, cette modularité est une conjecture de Serre ([196], lettre du 01/05/1973, et [195]), qui a finalement été prouvée par Khare, Wintenberger [143, 144] avec l'aide de Kisin [153], mais la preuve utilise pleinement les théorèmes de relèvement modulaire.

Pour la modularité des courbes elliptiques qui est le point de départ de toute cette histoire, Wiles [222] pouvait s'appuyer sur la modularité de $\overline{\rho}\colon \mathrm{Gal}_{\mathbf{Q}} \to \mathbf{GL}_2(\mathbf{F}_3)$ découlant[21] des travaux de Langlands [161] et Tunnell [212], mais c'est à peu près le seul résultat général à notre disposition à part les $\mathrm{I}_F^K \chi$, où $F/K$ est cyclique et $\chi$ est un caractère, pour lesquels on peut utiliser les travaux d'Arthur, Clozel [5] mentionnés dans le n° 3.3.5.

Heureusement, Taylor [206, th. 1.6] a réalisé que l'on pouvait prouver des résultats de modularité potentielle pour $\overline{\rho}$ et en déduire la modularité potentielle de $\rho$, ce qui est suffisant pour prouver les propriétés analytiques des fonctions $L$ si on veut bien se contenter de la méromorphie au lieu de l'holomorphie[22].

---

19. I.e. en restriction à $\mathrm{Gal}_{F_v}$ : pour $v \mid p$, on impose une condition du genre "de Rham de type fixé".

20. Il arrive que l'on ait besoin de conditions plus strictes du genre « $\rho$ assez proche d'une représentation modulaire ».

21. À ceci près que les résultats en question fournissent une forme modulaire de poids 1 et que l'on a besoin d'une forme de poids 2, mais ceci s'arrange en utilisant une congruence modulo 3 entre formes modulaires.

22. L'idée est de partir d'une représentation auxiliaire $\overline{\rho}'$ modulo un autre nombre premier $q$ dont on sait prouver la modularité (ainsi que celle de toutes ses restrictions aux extensions de $F$, ce qui impose de prendre $\overline{\rho}' = \mathrm{I}_F^K \chi$ vu l'état actuel de nos connaissances), et de transférer cette modularité à $\overline{\rho}$ (ce qui, à première vue, semble n'avoir aucun sens mais a déjà été utilisé par Wiles [222] pour déduire de la modularité modulo 5 de résultats de modularité modulo 3). Voir le n° 2.2 pour des détails ; signalons que les $\overline{\rho}$ dont on peut prouver la modularité potentielle sont obtenus comme réduction de la composante en $p$ de systèmes compatibles dont les poids de Hodge–Tate ne sont pas arbitraires (i.e. ces poids sont $0, 1, \ldots, n-1$), ce qui impose certaines contorsions dans les applications (sur le modèle de « l'astuce de Harris » [126]).

• *Espaces de déformations de représentations galoisiennes.* — Une fois que l'on s'est mis dans une situation où $\overline{\rho}$ est modulaire, il reste à comparer l'ensemble $X_{\overline{\rho}}^{\mathrm{mod}}$ des relèvements modulaires de $\overline{\rho}$ et l'intersection de deux sous-variétés[23] de l'espace $X_{\overline{\rho}}^{\mathrm{loc}}$ des déformations locales[24] de $\overline{\rho}$ : la sous-variété $X_{\overline{\rho}}^{\mathrm{glob}}$ des déformations globales et celle $X^{\mathrm{dR}}$ des représentations de Rham à poids de Hodge–Tate égaux à ceux de $\rho_p$. La conjecture de Fontaine–Mazur implique que $\dim X_{\overline{\rho}}^{\mathrm{loc}} - (\dim X_{\overline{\rho}}^{\mathrm{glob}} + \dim X^{\mathrm{dR}}) \geq 0$ (sinon on aurait des familles continues de représentations géométriques, ce qui est incompatible avec le fait que les variétés algébriques définies sur $\mathbf{Q}$ sont en nombre dénombrable). Si $\ell_0 = 0$, on peut espérer que les variétés se coupent transversalement et que l'intersection ne soit pas trop compliquée à décrire, mais si $\ell_0 \geq 1$, on tombe dans le monde des intersections atypiques, et la situation devient nettement plus délicate.

Notons qu'en fait, les variétés ci-dessus sont les fibres génériques de schémas sur $\mathbf{Z}_p$ et qu'au lieu d'une collection de points, on obtient un schéma de dimension 1 dont la géométrie est nettement plus riche que celle d'un ensemble fini de points, ce qui donne des outils supplémentaires pour le cerner.

*1.3.4. Modularité potentielle de fonctions L, le cas $\ell_0 = 0$.* — Le seul cas où $\ell_0 = 0$ si on se restreint aux $\rho \colon \mathrm{Gal}_F \to \mathbf{GL}_n(\overline{\mathbf{Q}}_p)$, est le cas $n = 2$ et $F$ totalement réel qui est celui considéré par Taylor [206, 208].

La condition $\ell_0 = 0$ pour des $\rho \colon \mathrm{Gal}_F \to \mathbb{G}(\overline{\mathbf{Q}}_p)$ où $\mathbb{G}$ est un groupe algébrique réductif implique[25] l'existence de variétés de Shimura pour un dual bien choisi de $\mathbb{G}$, et on peut utiliser les représentations sortant de la cohomologie étale de ces variétés de Shimura pour essayer de prouver l'implication « modulaire + algébrique ⇒ motivique » nécessaire à la mise en route de la stratégie esquissée au n° 1.3.3.

La décennie qui a suivi l'introduction du concept de modularité potentielle a vu un certain nombre d'avancées significatives *dans le cas autodual impair*; en particulier, les travaux de Barnet-Lamb, Clozel, Gee, Geraghty, Harris, Shepherd-Barron, Taylor prouvant la modularité potentielle des puissances symétriques des systèmes compatibles réguliers, de dimension 2, pour les corps totalement réels ([56, 128, 209] pour les premiers résultats sur $\mathbf{Q}$, [8] pour le cas général sur $\mathbf{Q}$ et [6] pour les corps totalement réels). Ces avancées ont été rendues possibles par les progrès (cf. th. 2.1) sur l'implication « modulaire + algébrique ⇒ motivique » consécutifs à la preuve du lemme fondamental pour les groupes unitaires par Laumon, Ngô [162] et Waldspurger [214].

---

23. Cette présentation est un petit peu mensongère car l'application $X_{\overline{\rho}}^{\mathrm{glob}} \to X_{\overline{\rho}}^{\mathrm{loc}}$ n'est pas forcément injective et pourrait même, a priori, avoir des fibres de dimension $> 0$. Que ce ne soit pas le cas se démontre via des théorèmes de relèvement modulaires, cf. [1] pour des résultats de ce genre.

24. $X_{\overline{\rho}}^{\mathrm{loc}} = \prod_{v|p} X_{\overline{\rho}}^v$, où $X_{\overline{\rho}}^v$ est l'espace des $(\rho' \colon \mathrm{Gal}_{F_v} \to \mathbb{G}(\overline{\mathbf{Z}}_p))/\sim$ relevant $\overline{\rho}$. (Il faut considérer les $\rho'$ à conjugaison près par un élément de $\mathbb{G}(\overline{\mathbf{Z}}_p)$ s'envoyant sur 1 dans $\mathbb{G}(\overline{\mathbf{F}}_p)$, i.e. on ne fixe pas de base : on considère des *représentations non encadrées*. Dans les preuves, il est souvent commode de ne pas passer au quotient, ce qui augmente la dimension des espaces ci-dessus et il faut un petit peu changer la rhétorique.)

25. Elle implique aussi la régularité des poids de Hodge–Tate et une condition sur l'action des conjugaison complexes (il faut que $\rho$ soit « impaire »).

Tous les travaux ci-dessus concernent le cas autodual, impair, régulier, et le résultat définitif, dans ce cadre, est le suivant, fruit des efforts de Barnet-Lamb, Gee, Geraghty, Patrikis, Taylor ([7, th. A] et [181, th. 2.1]).

Théorème 1.24. — *Soit $F$ un corps CM, et soit $\rho$ un système compatible de dimension $n$, irréductible ou pur, régulier, autodual et impair*[26]. *Alors $\rho$ est potentiellement modulaire.*

1.3.5. *Modularité potentielle de fonctions $L$, le cas $\ell_0 \geq 1$.* — Calegari, Geraghty [43] ont développé une stratégie pour l'étude des intersections atypiques. Il peut y avoir plusieurs raisons pour lesquelles $\ell_0 \geq 1$.

⋄ La dimension de $X_{\overline{\rho}}^{\mathrm{glob}}$ dépend de l'action adjointe des conjugaisons complexes sur l'algèbre de Lie $\mathfrak{g}$ de $\mathbb{G}$ : la formule d'Euler–Poincaré pour la cohomologie galoisienne d'un corps de nombres fournit la formule (conjecturale)

$$\dim X_{\overline{\rho}}^{\mathrm{glob}} = (r_1 + r_2) \dim \mathfrak{g} - \sum_{v \text{ réelle}} \dim \mathfrak{g}^{\mathrm{Frob}_v = 1}$$

Chacun des $\mathrm{Frob}_v$ est une involution d'algèbre de Lie et la dimension ci-dessus est maximale si toutes ces involutions sont des involutions de Cartan. Par exemple, si $\mathfrak{g} = \mathfrak{sl}_2$ et $F = \mathbf{Q}$, une involution de $\mathfrak{g}$ est soit l'identité (correspondant à $\overline{\rho}$ paire), soit une involution de Cartan (correspondant à $\overline{\rho}$ impaire) ; dans le cas pair $\dim X_{\overline{\rho}}^{\mathrm{glob}} = 0$ et on ne sait, jusqu'ici, rien faire par les méthodes de relèvement modulaire.

⋄ La dimension de $X^{\mathrm{dR}}$ dépend de la régularité des poids de Hodge–Tate : moins les poids de Hodge–Tate sont réguliers, plus cette dimension est petite. Par exemple, pour les représentations de dimension 2 de $\mathrm{Gal}_{\mathbf{Q}_p}$, cette dimension est 1 si les poids de Hodge–Tate sont distincts (cas régulier) et 0 s'ils sont égaux (voire $-1$). Les représentations d'Artin de dimension 2 d'un corps totalement réel $F$ relèvent de cette rubrique (voir le premier point du nᵒ 1.3.6), mais le premier résultat concernant un objet géométrique de dimension $\geq 1$ est celui impliquant le th. 1.7, à savoir :

Théorème 1.25. — (Boxer, Calegari, Gee, Pilloni [28]). *Si $A$ est une variété abélienne de dimension 2 définie sur un corps totalement réel, alors $h^i(A)$ est potentiellement modulaire pour tout $i$.*

(Le résultat est immédiat si $i \neq 1, 2, 3$. Comme $h^3(A)$ est un twist de $h^1(A)$ et $h^2(A) = \wedge^2 h^1(A)$, il n'y a que la modularité potentielle de $h^1(A)$ à prouver, celle de $h^3(A)$ en découlant immédiatement et celle de $h^2(A)$ se déduisant de celle de $h^1(A)$ via les résultats de Kim [147] et Henniart [136] rappelés dans le dernier point du nᵒ 3.3.4.)

⋄ Enfin, on peut avoir $\ell_0 \geq 1$ parce qu'il n'y a pas de variété de Shimura attaché au dual de $\mathbb{G}$. Pendant longtemps cela a constitué un obstacle absolu car l'implication « modulaire ⇒ motivique » manquait à l'appel. Les travaux de Harris, Lan, Taylor,

---

[26]. Cela signifie que $\rho_p$ se factorise à travers $\mathbf{GSp}_n$ avec un multiplicateur totalement impair ou à travers $\mathbf{GO}_n$ avec un multiplicateur totalement pair.

THORNE [127] et ceux de SCHOLZE [191] ont débloqué la situation pour $\mathbf{GL}_n$ (qui est son propre dual) sur un corps CM. Cela a permis à [27] ALLEN, BOXER, CALEGARI, CARAIANI, GEE, HELM, LE HUNG, NEWTON, SCHOLZE, TAYLOR, THORNE de prouver le résultat suivant ([2, cor. 7.1.12] et [27, th. 6.2.1]) qui est le premier sans condition de parité :

THÉORÈME 1.26. — *Soit $F$ un corps CM imaginaire[28]. Soit $(\rho_p)_p$ un système compatible sur $F$, de dimension 2, à poids parallèles distincts[29]. Alors $\mathrm{Sym}^k \circ \rho$ est potentiellement modulaire, pour tout $k \geq 1$.*

Ce résultat est l'analogue, pour les corps CM, des résultats[30] de [8] (pour $F = \mathbf{Q}$) et [6] (pour $F$ totalement réel) ; on peut espérer qu'un renforcement des théorèmes de relèvement modulaire permettra de prouver un analogue du th. 1.24 dans un avenir pas trop lointain.

*Remarque 1.27.* — (i) Chacun des th. 1.25 et 1.26 implique que, si $E$ est une courbe elliptique définie sur un corps CM, alors $L(E, s)$ admet un prolongement méromorphe et une équation fonctionnelle (c'est clairement un cas particulier du th. 1.26 et c'est une conséquence du th. 1.25 appliqué à la restriction des scalaires de $K$ à $F$ de $E$, qui est une variété abélienne de dimension 2 définie sur $F$).

(ii) Si $C$ est une courbe de genre 2 définie sur un corps totalement réel, le th. 1.25 appliqué à la jacobienne de $C$ implique que $h^1(C)$ est potentiellement modulaire et donc que $C$ vérifie la conjecture de Hasse–Weil.

(iii) On peut aussi appliquer le th. 1.25 à la restriction des scalaires de $K$ à $F$ d'une courbe elliptique $E$ définie sur une extension quadratique $K$ de $F$, non CM. On en déduit que $E$ est potentiellement modulaire, ce qui semble être le premier résultat de ce type (en dehors de celui des représentations d'Artin) valable dans le cas non CM.

*1.3.6. Modularité de fonctions $L$*. — Comme le montre le cas des fonctions $L$ d'Artin, passer de la modularité potentielle à la modularité (i.e. de la méromorphie d'une fonction $L$ à son holomorphie) est loin d'être une promenade de santé. Il y a quand même eu quelques percées remarquables :

• Un des premiers grands succès du programme de Langlands a été la preuve de la conjecture d'Artin pour des représentations de dimension 2. Si $\rho \colon \mathrm{Gal}_F \to \mathrm{GL}_2(\mathbf{C})$ est une $\mathbf{C}$-représentation de dimension 2, l'image de $\mathrm{Gal}_F$ dans $\mathrm{PGL}_2(\mathbf{C})$ est abélienne, diédrale, ou bien isomorphe à $A_4$, $S_4$ ou $A_5$. La conjecture d'Artin dans les cas abélien et diédral a été prouvée par Artin lui-même, celle dans les cas $A_4$ et $S_4$ par Langlands et Tunnell [161, 212]. Le cas de $A_5$ a résisté, jusqu'ici, aux méthodes classiques.

---

27. Contrairement à ce que les dates de (pré)publication suggéreraient, il s'est passé pas mal de temps entre l'annonce des résultats de [2] (CIRM, décembre 2016) et [27].

28. Le cas $F$ totalement réel est couvert par les résultats de [6].

29. Cela équivaut à ce qu'il existe un caractère $\eta$ de $\mathbf{A}_F^*/F^*$ tel que les poids de Hodge–Tate de $\rho \otimes \eta$ soient $0, m$, avec $m \geq 1$, pour tout $\tau \in \mathrm{Hom}(F, \overline{\mathbf{Q}})$.

30. Comme pour ces résultats, la preuve dans le cas général utilise une variante de l'astuce de HARRIS [126] permettant de fabriquer un système compatible à poids de Hodge–Tate consécutifs à partir d'un système dont les poids de Hodge–Tate sont régulièrement espacés.

Par contre les méthodes $p$-adiques ont permis de traiter ce cas pour les représentations *impaires* (celles pour lesquelles det $\rho(c) = -1$, si $c$ est une conjugaison complexe, i.e. $\rho(c)$ est une symétrie et pas une homothétie) dans le cas où $F$ est totalement réel : résultat d'une série de travaux commençant avec l'article fondateur de Buzzard, Taylor [40] et se terminant par celui de Pilloni, Stroh [183]. Le cas des représentations paires résiste toujours : il n'y a *aucune représentation paire* de Gal$_\mathbf{Q}$, d'image $A_5$ dans $\mathrm{PGL}_2(\mathbf{C})$, dont on sache prouver la modularité.

• On a déjà mentionné la modularité [33] des courbes elliptiques définies sur $\mathbf{Q}$. Ce résultat a été étendu par Freitas, Le Hung, Siksek [113] aux courbes elliptiques définies sur un corps quadratique réel, et par Thorne [210] aux courbes elliptiques définies sur la $\mathbf{Z}_p$-extension cyclotomique de $\mathbf{Q}$ (encore un corps totalement réel).

• Les résultats de [2] ont permis de sortir du cadre totalement réel :

  ⋄ Allen, Khare, Thorne [3] ont prouvé que, si $K$ est un corps CM ne contenant pas $\mathbf{Q}(\boldsymbol{\mu}_5)$, une proportion positive de courbes elliptiques définies sur $K$ sont modulaires.

  ⋄ Caraiani, Newton [45] ont prouvé que les courbes elliptiques définies sur $\mathbf{Q}(\sqrt{-D})$, pour $D = 1, 2, 3, 5$, sont modulaires. Ils ont prouvé, plus généralement, le même énoncé pour les $D > 0$ pour lesquels la courbe modulaire $X_0(15)$ — c'est une courbe elliptique — n'a qu'un nombre fini de points rationnels sur $\mathbf{Q}(\sqrt{-D})$ (soit 50% des $\mathbf{Q}(\sqrt{-D})$ modulo des conjectures standard (i.e. Birch et Swinnerton–Dyer, et Goldfeld).

• Tous les résultats précédents concernent des systèmes de dimension 2. Le résultat suivant est le premier concernant des systèmes de dimension arbitraire.

Théorème 1.28. — (Newton, Thorne [171, 172]). *Si $f$ est une forme modulaire de poids $m \geq 2$, alors $\mathrm{Sym}^k \circ L(f, s)$ est modulaire pour tout $k \geq 1$.*

*Remarque 1.29.* — (i) On peut utiliser des méthodes purement automorphes pour attaquer ce genre de problèmes (Gelbart, Jacquet [120] pour $k = 2$, Kim, Shahidi [148, 149, 147] pour $k = 3, 4$). Les méthodes automorphes s'appliquent à n'importe quelle représentation automorphe de $\mathbf{GL}_2(\mathbf{A}_F)$ (pas seulement celles associées aux formes modulaires), sans restriction sur $F$, mais il n'y a pas vraiment eu de progrès dans cette direction depuis 20 ans.

(ii) L'approche automorphe semblant, provisoirement, dans une impasse, Newton, Thorne passent par les représentations galoisiennes, suivant la voie empruntée par Dieulefait [82] et par Clozel, Thorne [58, 59], qui avaient prouvé la modularité de $\mathrm{Sym}^k \circ L(f, s)$ pour $k \leq 8$.

(iii) Newton, Thorne ont, depuis, étendu le th. 1.28 aux corps totalement réels [173] (i.e., aux formes modulaires de Hilbert).

• Un autre résultat, purement automorphe, qu'on ne sait prouver qu'en passant par les représentations galoisiennes est le théorème de changement de base suivant.

Théorème 1.30. — (Dieulefait [81, 82], Clozel, Newton, Thorne [57]).

Si $f$ est une forme modulaire de poids $\geq 2$, alors $\mathrm{R}_{\mathbf{Q}}^F \circ \mathrm{Sym}^k \circ L(f,s)$ est modulaire pour tout $k \geq 1$ et tout corps totalement réel $F$.

• Enfin, mentionnons que Boxer, Calegari, Gee, Pilloni [28] ont, en sus de leur résultat de modularité potentielle, prouvé l'existence d'une infinité de courbes de genre 2, modulaires, définies sur $\mathbf{Q}$, mais sans vrai contrôle sur l'ensemble de courbes ainsi obtenues. Ils ont, depuis, trouvé une approche alternative qui permet de prouver la modularité dans beaucoup de cas [29] ; par exemple, on a le résultat suivant :

Théorème 1.31. — Soit $f \in \mathbf{Q}[x]$, unitaire, de degré 5, sans racine double, et soit $C$ la courbe d'équation $y^2 = f(x)$. On suppose que $C$ a bonne réduction ordinaire en 2 et 3 et que l'image de $\mathrm{Gal}_{\mathbf{Q}}$ sur les points de 3-torsion de la jacobienne de $C$ est aussi grosse que possible (i.e. $\mathbf{GSp}_4(\mathbf{F}_3)$). Alors $C$ est modulaire.

## 1.4. Applications

*1.4.1. Équirépartition.* — L'existence d'un prolongement analytique et la non annulation sur la droite $\mathrm{Re}(s) = 1 + \frac{w}{2}$ ont des conséquences sur la répartion des nombres premiers et des éléments de Frobenius. Par exemple :

• La non annulation de $\zeta(s)$ sur la droite $\mathrm{Re}(s) = 1$ est équivalente au théorème des nombres premiers $|\{\ell \leq x\}| \sim \frac{x}{\log x}$.

• La non annulation des fonctions $L$ de Dirichlet sur la droite $\mathrm{Re}(s) = 1$ implique le théorème de la progression arithmétique sous la forme : si $a$ et $N$ sont premiers entre eux, $|\{\ell \leq x,\ \ell \equiv a \bmod N\}| \sim \frac{1}{\varphi(N)} \frac{x}{\log x}$.

• La non annulation des fonctions $L$ d'Artin sur la droite $\mathrm{Re}(s) = 1$ implique le théorème de Čebotarev : si $K/F$ est une extension galoisienne finie de corps de nombres de groupe de Galois $G$, alors $|\{v,\ |v| \leq x \text{ et } \mathrm{Frob}_v \in C\}| \sim \frac{|C|}{|G|} \frac{x}{\log x}$ si $C$ est une classe de conjugaison de $G$. (Autrement dit les frobenius s'équirépartissent dans les groupes de Galois.)

• Soit $E$ une courbe elliptique définie sur un corps de nombres $F$. Si $v$ est une place de $F$, posons $a_v = |v| + 1 - |E(\mathcal{O}_v/\mathfrak{m}_v)|$. D'après Hasse [130], on a $|a_v| \leq 2|v|^{1/2}$ (hypothèse de Riemann pour les courbes elliptiques), et donc $\frac{a_v}{2|v|^{1/2}} = \cos\theta_v$ avec $\theta_v \in [0,\pi]$.

On a alors la conjecture suivante [31] dans laquelle $\delta_\theta$ désigne la masse de Dirac en $\theta$ :

---

31. Cette conjecture apparaît sous la forme d'un exercice [196] dans une lettre de Tate à Serre du 05/08/1963, avec des indications pour une solution dans une lettre du 28/08/1963. Le lecteur est invité à consulter aussi les échanges des 1-3/03/2008 pour des commentaires sur la genèse de cette conjecture.

Conjecture 1.32. — (Sato–Tate)

$$\lim_{x\to\infty}\left(\frac{\log x}{x}\sum_{|v|\leq x}\delta_{\theta_v}\right) = \begin{cases} \frac{2}{\pi}\sin^2\theta d\theta & \text{si } \mathrm{End}_{\overline{\mathbf{Q}}}E = \mathbf{Z}, \\ \frac{1}{\pi}d\theta & \text{si } \mathrm{End}_{\overline{\mathbf{Q}}}E \neq \mathbf{Z} \text{ mais } \mathrm{End}_F E = \mathbf{Z}, \\ \frac{1}{2\pi}(\delta_\pi + d\theta) & \text{si } \mathrm{End}_F E \neq \mathbf{Z}. \end{cases}$$

La quantité $a_v$ ci-dessus est la trace de $\mathrm{Frob}_v^{-1}$ agissant sur $H^1_{\text{ét}}(E_{\overline{\mathbf{Q}}}, \mathbf{Q}_\ell)$ pour $\ell$ distinct de la caractéristique de $\mathcal{O}_v/\mathfrak{m}_v$. L'énoncé ci-dessus équivaut à un énoncé d'équirépartition[32] des $\frac{1}{|v|^{1/2}}\mathrm{Frob}_v$ dans un certain groupe de Lie compact $H$ naturellement associé à $E$ (le groupe de Sato–Tate, variante du groupe de Mumford–Tate). Ce point de vue, mis en avant par Serre [193], se généralise à tous les systèmes compatibles. Serre a aussi expliqué comment déduire cette équirépartition de la non-annulation sur $\mathrm{Re}(s) = 1$ des fonctions $L$ attachées aux représentations de dimension finie de $H$.

Dans le cas des courbes elliptiques, on a $H = \mathbf{SU}_2$, $\mathbf{U}_1$ ($\cong \mathbf{SO}_2$) ou le normalisateur de $\mathbf{U}_1$ dans $\mathbf{SU}_2$, suivant que l'on est dans le premier, le second ou le troisième cas. Si $\mathrm{End}_{\overline{\mathbf{Q}}}E \neq \mathbf{Z}$ (i.e. si $E$ a de la multiplication complexe), le résultat est facile à prouver. Dans le cas $\mathrm{End}_{\overline{\mathbf{Q}}}E = \mathbf{Z}$, les premiers résultats n'ont été obtenus que récemment ; à ce jour, ce sont les suivants :
- $F = \mathbf{Q}$ (combinaison de [56, 128, 209] si $j(E) \notin \mathbf{Z}$, et [8] pour le cas général).
- $F$ est un corps totalement réel [6].
- $F$ est un corps CM [2].

*Remarque 1.33.* — Les résultats ci-dessus ne permettent pas de borner de manière effective le terme d'erreur dans la convergence vers la mesure limite à cause des choix, que l'on ne contrôle absolument pas, d'extensions auxiliaires dans la modularité potentielle. Le th. 1.28, qui évite ces choix, rend ce terme d'erreur effectif.

*Remarque 1.34.* — Pour les variétés abéliennes de dimension 2, il y a 52 sous-groupes $H$ de $\mathbf{USp}_4 := \mathbf{Sp}_4 \cap \mathbf{SU}_4$ possibles dont 34 pour les variétés définies sur $\mathbf{Q}$ (cf. [100]). La conjecture de Sato–Tate pour les variétés définies sur $\mathbf{Q}$ a en grande partie été démontrée [142, 204] sauf dans le cas $\mathrm{End}_{\overline{\mathbf{Q}}}A = \mathbf{Z}$ où il faudrait étendre le th. 1.25 aux systèmes $(r \circ V_p A)_p$, où $r$ décrit les représentations algébriques de $\mathbf{GSp}_4$ (le th. 1.25 correspond au cas où $r$ est la représentation standard de dimension 4, obtenue en définissant $\mathbf{GSp}_4$ comme sous-groupe de $\mathbf{GL}_4$). Cela semble totalement hors de portée à l'heure actuelle car les $r \circ V_p A$ ont des poids de Hodge–Tate de multiplicité arbitrairement grande.

*1.4.2. Pureté.* — Si $X$ est une variété propre et lisse définie sur $F$, il résulte de l'hypothèse de Riemann pour les variétés sur les corps finis (th. 1.9) que, pour tout $i$, $[L(h^i(X), s)]$ est pure (de poids $i$).

---

32. Comme $\mathrm{Frob}_v$ n'est défini qu'à conjugaison près, il s'agit d'une équirépartition sur l'espace des classes de conjugaison de $H$ pour la mesure image de la mesure de Haar sur $H$.

Ramanujan a conjecturé que $L(\Delta, s)$ est pure de poids 11 (où $\Delta$ est la forme modulaire de poids 12 habituelle : $\Delta(z) = q \prod_{n\geq 1}(1-q^n)^{24}$ avec $q = e^{2i\pi z}$). Ceci a été prouvé par Deligne comme conséquence de sa preuve de l'hypothèse de Riemann pour les variétés sur les corps finis. Si $\pi$ est la représentation automorphe de $\mathbf{GL}_2(\mathbf{A_Q})$ engendrée par $\Delta$, la pureté de $L(\Delta, s)$ équivaut à celle de $L(\pi, s)$. Cela conduit à l'énoncé suivant.

Conjecture 1.35. — (Conjecture de Ramanujan généralisée). *Soit $\pi$ une représentation automorphe cuspidale de $\mathbf{GL}_n(\mathbf{A}_F)$. Si le caractère central de $\pi$ est de poids $w$, alors $L(\pi, s)$ est pure, de poids $w$.*

Jusqu'à récemment, les rares cas connus de cette conjecture reposaient sur le fait que $L(\pi, s) = L(\rho, s)$ où $\rho$ est un système compatible sous-quotient de la cohomologie d'une variété de Shimura, auquel on peut appliquer les généralisations de l'hypothèse de Riemann pour les variétés sur les corps finis de [75].

Notons que cette conjecture est une conséquence de la fonctorialité de Langlands [160] : elle découle, via un argument analogue à celui que Deligne a utilisé pour prouver l'hypothèse de Riemann pour les variétés sur les corps finis, de la modularité de $L(\pi, r, s)$ pour suffisamment de représentations $r$ de $\mathbf{GL}_n(\mathbf{C})$. De fait, la modularité potentielle de suffisamment de $L(\pi, r, s)$ suffit, ce qui permet de prouver le résultat suivant ([2, cor. 7.1.15] et [27, th. 7.1.1]).

Théorème 1.36. — *Soit $F$ un corps CM non totalement réel. Soit $\pi$ une représentation automorphe cuspidale de $\mathbf{GL}_2(\mathbf{A}_F)$, algébrique et régulière. On suppose que $\pi$ est à poids parallèles*[(33)]. *Alors, $\pi$ vérifie la conjecture de Ramanujan généralisée.*

La preuve repose sur le th. 1.26, et utilise le système compatible associé à $\pi$ ; ce n'est donc pas une preuve purement automorphe (mais comme le système compatible associé à $\pi$ n'est pas construit comme sous-quotient de la cohomologie d'une variété algébrique, on ne peut utiliser les résultats de [75] pour conclure comme dans le cas de $\mathbf{Q}$ ou d'un corps totalement réel).

*1.4.3. L'équation de Fermat.* — Soient $F$ un corps de nombres et $p \geq 3$ un nombre premier. Si $a^p + b^p = c^p$, avec $a, b, c \in F$ et $abc \neq 0$, on peut utiliser les propriétés remarquables (i.e. la faible ramification de la représentation de $\mathrm{Gal}_F$ sur les points de $p$-torsion) de la courbe elliptique d'équation $y^2 = x(x - a^p)(x + b^p)$, introduite par Hellegouarch [134] et popularisée par Frey [115, 116], pour en déduire des renseignements sur $(a, b, c)$. Pour pouvoir exploiter ces propriétés, on a besoin de savoir *a priori* que cette courbe est modulaire ; la faible ramification des points de $p$-torsion implique alors une congruence modulo $p$ avec une "forme modulaire" de niveau petit (ne faisant pas intervenir $p$), et donc vivant dans un espace que l'on maîtrise bien.

Dans le cas $F = \mathbf{Q}$, la « conjecture $\varepsilon$ » de Serre [195], démontrée par Ribet [185],

---

33. Cela équivaut à l'existence d'un caractère de Hecke algébrique $\chi$ tel que $\pi \otimes \chi$ soit de poids $(0, m-1)$ pour tout $\tau \in \mathrm{Hom}(F, \overline{\mathbf{Q}})$, ce qui est automatique si $F$ est une extension quadratique de $\mathbf{Q}$.

implique une congruence modulo $p$ avec une forme modulaire parabolique de poids 2 et de niveau 1 ou 2. Comme il n'existe pas de telle forme modulaire, la modularité des courbes elliptiques définies sur $\mathbf{Q}$ prouve le théorème de Fermat [222].

FREITAS et SIKSEK [110, 111] ont utilisé la modularité des courbes elliptiques sur les corps quadratiques réels pour prouver les résultats suivants :

THÉORÈME 1.37. — *Soit $K = \mathbf{Q}(\sqrt{d})$, avec $d$ entier $> 0$, sans facteur carré.*

(i) *Si $d = 3, 6, 7, 10, 11, 13, 14, 15, 19, 21, 22, 23$, l'équation $a^p + b^p = c^p$ n'a que des solutions triviales*[34] *dans $K$.*

(ii) *Si $d$ est congru à 3, 6, 10 ou 11 modulo 16, il existe une constante $B_K$ effectivement calculable telle que, pour tout $p \geq B_K$, l'équation $a^p + b^p = c^p$ n'a que des solutions triviales dans $K$.*

FREITAS, KRAUS et SIKSEK [112] ont utilisé la modularité des courbes elliptiques sur la $\mathbf{Z}_p$-extension cyclotomique de $\mathbf{Q}$ pour prouver des résultats analogues pour l'équation de Fermat sur cette extension. Les résultats de CARAIANI, NEWTON laissent espérer des résultats analogues sur les corps quadratiques imaginaires [35].

*1.4.4. Une application exotique.* — Soit $\psi_p \colon \mathbf{F}_p \to \mathbf{C}^*$ un caractère non trivial. Si $q$ est une puissance de $p$ et si $a \in \mathbf{F}_q^*$, il existe [215] $\alpha_{q,a}, \beta_{q,a} \in \mathbf{C}$ vérifiant

$$|\alpha_{q,a}| = |\beta_{q,a}| = \sqrt{q}, \quad \alpha_{q,a}\beta_{q,a} = q \quad \text{et} \quad \alpha_{q,a} + \beta_{q,a} = -\sum_{x \in \mathbf{F}_q^*} \psi\Big(\mathrm{Tr}_{\mathbf{F}_q/\mathbf{F}_p}(x + \tfrac{a}{x})\Big)$$

(somme de Kloosterman). Si $k$ est un entier $\geq 1$, on pose alors

$$m_k(q) = \sum_{a \in \mathbf{F}_q^*} \sum_{i=0}^{k} \alpha_{q,a}^i \beta_{q,a}^{k-i}$$

et on fabrique la fonction génératrice

$$Z_{k,p}(T) = \exp\Big(\sum_{n=1}^{\infty} m_k(p^n) \tfrac{T^n}{n}\Big)$$

Alors $Z_{k,p} \in 1 + T\mathbf{C}[T]$ et se factorise sous la forme $E_{k,p}R_{k,p}$ avec $E_{k,p}, R_{k,p} \in 1 + T\mathbf{C}[T]$, toutes les racines de $E_{k,p}$ sont de valeur absolue $p^{-(k+1)/2}$ et celles de $R_{k,p}$ sont de valeur absolue $> p^{-(k+1)/2}$. Enfin, on définit une fonction $L_k(s)$ par $L_k(s) = \prod_{p>k} E_{k,p}(p^{-s})^{-1}$.

THÉORÈME 1.38. — (FRESÁN, SABBAH, YU [114]) *$L_k$ admet un prolongement analytique et une équation fonctionnelle reliant $s$ et $k + 2 - s$.*

La preuve consiste à vérifier que $L_k$ est la fonction $L$ d'un système compatible autodual, régulier, ce qui permet d'appliquer le th. 1.24 ; la vérification de la régularité utilise des idées venant de la théorie des motifs exponentiels.

---

34. Vérifiant $abc = 0$.

35. En incluant dans les solutions triviales, si $K = \mathbf{Q}(\sqrt{-3})$, les solutions découlant de l'identité $1 + \zeta + \zeta^2 = 0$ si $\zeta \neq 1$ est une racine cubique de l'unité (l'existence de ces solutions est probablement l'obstacle à une preuve "simple" du théorème de Fermat).

## 2. REPRÉSENTATIONS GALOISIENNES ET MODULARITÉ

Ce chapitre un peu décousu est consacré à l'existence ou la construction de représentations du groupe de Galois d'un corps de nombres [36] obéissant à des contraintes variées via diverses techniques.

### 2.1. Construction de représentations galoisiennes

*2.1.1. Représentations automorphes polarisables.* — La cohomologie des variétés de Shimura fournit une source importante de systèmes compatibles, dont on peut prouver la modularité. Par exemple, les travaux d'Eichler [87] et Shimura [199] associent à une forme modulaire primitive $f$, de poids 2, un système compatible découpé dans le $H^1$ d'une courbe modulaire, dont la fonction $[L]$ est la classe de celle de $L(f, s)$ ; un tel système compatible est donc modulaire. Ces résultats ont été étendus aux formes de poids $\geq 2$ par Deligne [73], et aux formes de poids 1 par Deligne, Serre [76] mais dans ce cas le système est obtenu par un passage à la limite à partir de systèmes associés à des formes de poids $\geq 2$ ; il ne se réalise pas dans la cohomologie étale d'une courbe modulaire [37]. Toutes les représentations obtenues de la sorte sont *impaires* ; pour obtenir des représentations paires, il faudrait pouvoir associer des systèmes compatibles aux formes de Maass de valeur propre $\frac{1}{4}$ pour le laplacien, mais on ne sait pas quel objet géométrique pourrait fournir ces systèmes compatibles.

Le th. 2.1 ci-dessous est une vaste généralisation des résultats précédents. Soient $F$ un corps CM, $F^+$ son sous-corps totalement réel maximal (donc $[F : F^+] = 1$ ou 2). Une représentation automorphe $\pi$ de $\mathbf{GL}_n(\mathbf{A}_F)$ est dite *polarisable* s'il existe un caractère $\chi \colon \mathbf{A}_{F^+}^*/(F^+)^* \to \mathbf{C}^*$ tel que $\chi_v(-1)$ ne dépend pas de $v \mid \infty$ et $\pi^c \cong \check{\pi} \otimes \chi \circ \mathrm{N}_{F/F^+}$ (où $c$ est la conjugaison complexe, que l'on fait agir sur $\mathbf{A}_F$ pour définir $\pi^c$ ; si $F = F^+$, alors $\pi^c = \pi$).

On a le résultat suivant [7, th. 2.1.1] qui regroupe des résultats de Clozel [54], Kottwitz [156], Clozel–Harris–Labesse [157, 55], Shin [201], Bellaïche–Chenevier [13], Caraiani [44],...

Théorème 2.1. — *Soit $\pi$ une représentation automorphe cuspidale de $\mathbf{GL}_n(\mathbf{A}_F)$, polarisable, algébrique et régulière. Alors, on peut associer à $\pi$ un système fortement compatible $\rho_\pi = (\rho_{\pi,p})$ qui est fortement modulaire et tel que $\Pi(\rho_\pi) = \pi$. De plus, $\rho_{\pi,p}$ est autoduale et impaire.*

---

36. Notons que les problèmes correspondant pour les courbes sur un corps fini ont été résolus par Drinfeld [83] et les frères Lafforgue [158, 159]. Comme, grâce à Grothendieck [124], on sait que les fonctions $L$ des systèmes locaux sur ces courbes ont les propriétés attendues, on en déduit la modularité de ces systèmes locaux via les théorèmes inverses du § 3.2 (ou plutôt leurs avatars sur les corps de fonctions). Les techniques de modularité potentielle permettent [23] de prouver l'automorphie potentielle de ces systèmes locaux pour le groupe prédit par la correspondance de Langlands.

37. Par contre, il se réalise dans la cohomologie complétée de la tour des courbes modulaires ou encore dans la cohomologie de la tour complétée des courbes modulaires, une courbe perfectoïde.

Beaucoup des représentations dont le théorème affirme l'existence apparaissent dans la cohomologie de variétés de Shimura (BLASIUS, ROGAWSKI [22] pour $n = 2$ et SHIN [201] pour $n$ général ; les résultats de SHIN utilisent le « lemme fondamental » pour les groupes unitaires prouvé par LAUMON, NGÔ [162] et WALDSPURGER [214], et complètent des travaux antérieurs de CLOZEL [54] et KOTTWITZ [156]). Pour les systèmes de représentations correspondants, la compatibilité est automatique puisque ces systèmes sont motiviques. Il y a malheureusement des cas où on ne sait pas comment trouver un motif dont les réalisations forment le système voulu ; la construction du sytème correspondant se fait par passage à la limite en utilisant des congruences avec les représentations apparaissant dans la cohomologie de variétés de Shimura (CHENEVIER, HARRIS [51], s'inspirant de techniques introduites par WILES [221] et TAYLOR [205] dans le cas $n = 2$) ; il faut alors se battre pour prouver la pureté des systèmes obtenus et le fait que les représentations du système sont de Rham, et pour vérifier la compatibilité avec la correspondance de Langlands locale (cf. [51] et [44]).

*Exemple 2.2.* — Soit $K$ un corps quadratique réel.

(i) Si $\pi$ est la représentation associée à une forme modulaire de Hilbert de poids $(2, 2)$ pour $K$ dont le $q$-développement est à coefficients dans $\mathbf{Q}$, et si $\pi$ est non ramifiée en toute place finie, alors $\rho_\pi$ devrait être le système associé à une courbe elliptique sur $K$ avec bonne réduction partout (il n'existe pas de telle courbe elliptique sur $\mathbf{Q}$, mais rien n'empêche qu'il en existe sur $K$) mais, dans ce cas, la construction de $\rho_\pi$ se fait par passage à la limite et ne fournit pas de motif que l'on pourrait relier à une courbe elliptique.

(ii) SHIMURA a raffiné la conjecture de Taniyama–Weil (modularité de la fonction $L$) sous la forme « si $E$ est une courbe elliptique, alors il existe $N$ tel que $E$ soit un quotient de la courbe modulaire $X_0(N)$ » : l'existence d'une telle paramétrisation modulaire implique la conjecture de Taniyama–Weil pour $E$ grâce aux travaux de EICHLER [87] et SHIMURA [199] ; réciproquement, la conjecture de Taniyama–Weil pour $E$, combinée à la conjecture de Tate (selon laquelle deux variétés abéliennes ayant même fonction $L$ sont isogènes), implique l'existence d'une paramétrisation modulaire, et comme FALTINGS [96] a prouvé la conjecture de Tate, les conjectures de Taniyama–Weil et de Shimura sont « équivalentes ».

L'analogue pour une courbe elliptique $E$ définie sur $K$ serait que $E$ est un quotient d'une courbe de Shimura ; c'est le cas pour beaucoup de $E$ (par exemple si $E$ a réduction multiplicative en au moins une place) mais ce n'est pas possible si $E$ a partout bonne réduction [38], et il est difficile de formuler un analogue de la conjecture de Shimura couvrant tous les cas.

*2.1.2. Le cas non polarisable.* — On conjecture que le th. 2.1 est valable sans les hypothèses restrictives sur $F$ ou sur $\pi$. Plus précisément, on a la conjecture suivante :

---

38. On aurait besoin d'une algèbre de quaternions $D$ sur $K$, déployée en une seule place à l'infini (pour avoir une courbe), et en toutes les places finies, ce qui contredit la nullité de $\sum_v \mathrm{inv}_v(D)$.

CONJECTURE 2.3. — (CLOZEL [53]). *Soit $F$ un corps de nombres et soit $\pi$ une représentation automorphe cuspidale algébrique de $\mathbf{GL}_n(\mathbf{A}_F)$. Alors il existe un système fortement compatible $\rho_\pi$ qui est fortement modulaire et tel que $\Pi(\rho_\pi) = \pi$.*

*Remarque 2.4.* — La philosophie de Langlands suggère que cette conjecture devrait avoir un pendant pour tout groupe réductif et pas seulement pour $\mathbf{GL}_n$; une telle conjecture a été formulée par BUZZARD, GEE [39] (c'est plus subtil que ce que l'on aurait pu croire).

Le point de départ de la construction des systèmes compatibles associés à une représentation automorphe cuspidale autoduale $\pi$ est le fait que l'on peut obtenir $\pi$ par fonctorialité à partir d'une représentation automorphe d'un groupe relié à des variétés de Shimura, et on construit le système cherché en utilisant la cohomologie étale de ces variétés de Shimura.

Sans cette condition d'autodualité, il est impossible de transférer $\pi$ sur un groupe donnant naissance à des variétés de Shimura. Ce qui reste vrai est que, si $\pi$ est algébrique régulière, alors $\pi$ se réalise dans la cohomologie d'un système local $\mathscr{L}_\pi$ sur la tour des espaces localement symétriques associés au groupe $\mathbb{G} = \mathbf{GL}_n$ sur $F$, i.e. la tour des

$$Y(K) := \mathbb{G}(F)\backslash\mathbb{G}(\mathbf{A}_F)/K$$

où $K = K_\infty K^{]\infty[}$, où $K_\infty$ est un sous-groupe compact maximal de $\mathbb{G}(\mathbf{R} \otimes F)$ fixé et $K^{]\infty[}$ décrit les sous-groupes ouverts du sous-groupe compact maximal $\prod_v \mathbb{G}(\mathscr{O}_v)$ de $\mathbb{G}(\mathbf{A}_F^{]\infty[})$. Par contre, comme les $Y(K)$ n'ont pas de structure de variété algébrique, il n'y a pas d'action galoisienne évidente sur $\overline{\mathbf{Q}}_\ell \otimes H^*(Y(K), \mathscr{L}_\pi)$ et la construction du système associé à $\pi$ est restée purement hypothétique jusqu'aux travaux de HARRIS, LAN, TAYLOR, THORNE [127] et ceux de SCHOLZE [191].

THÉORÈME 2.5. — *Si $F$ est un corps CM et si $\pi$ une représentation automorphe cuspidale algébrique régulière de $\mathbf{GL}_n(\mathbf{A}_F)$, alors il existe un système fortement compatible $\rho_\pi$ qui est modulaire, tel que $\Pi(\rho_\pi) = \pi$, et qui est compatible à la correspondance de Langlands locale en tous les $v$ non ramifiés.*

*Remarque 2.6.* — L'énoncé ci-dessus est dû à HARRIS, LAN, TAYLOR, THORNE. La preuve de SCHOLZE passe par celle d'un énoncé analogue pour des coefficients de torsion, ce qui a un intérêt propre car, dès que l'on quitte le monde des variétés de Shimura, la cohomologie des $Y(K)$ (à coefficients dans $\mathbf{Z}$) contient beaucoup de classes de torsion.

Si $K_v^0 = \mathbf{GL}_n(\mathscr{O}_v)$, soit $\mathbb{T}_v = \mathbf{Z}[K_v^0 \backslash \mathbf{GL}_n(F_v)/K_v^0]$ l'algèbre des fonctions bi-$K_v^0$-invariantes, à support compact (c'est une algèbre commutative), et soit $\mathbb{T}^{]S[} = \otimes_{v \notin S} \mathbb{T}_v$, où $S$ est un ensemble fini de places de $F$ contenant les places au-dessus de $p$. Alors $\mathbb{T}^{]S[}$ agit sur $H^i(Y(K), \overline{\mathbf{F}}_p)$ (cohomologie de Betti) pour tout $K = K_\infty K^{]\infty[}$ où $K^{]\infty[} = K_S \times \prod_{v \notin S} K_v^0$ et $K_S$ est un sous-groupe ouvert de $\prod_{v \in S} \mathbf{GL}_n(\mathscr{O}_v)$. Ce qui joue le rôle d'une représentation automorphe est un caractère $\lambda \colon \mathbb{T}^{]S[} \to \overline{\mathbf{F}}_p$ apparaissant dans la limite inductive des $H^i(Y(K), \overline{\mathbf{F}}_p)$, pour $K$ comme ci-dessus, i.e. il existe $v \in H^i(Y(K), \overline{\mathbf{F}}_p)$ non nul, tel que $T \cdot v = \lambda(T)v$ pour tout $T \in \mathbb{T}^{]S[}$.

SCHOLZE associe à un tel caractère une représentation $\rho_\lambda \colon \mathrm{Gal}_F \to \mathbf{GL}_n(\overline{\mathbf{F}}_p)$, la compatibilité avec la correspondance de Langlands locale pour $v \notin S$ étant remplacée par la relation [39] $a_i = |v|^{i(i-1)/2}\lambda(T_{v,i})$ entre les coefficients du polynôme caractéristique $X^n - a_1 X^{n-1} + \cdots + (-1)^n a_n$ de $\mathrm{Frob}_v^{-1}$ et les $\lambda(T_{v,i})$ où $T_{v,i} \in \mathbb{T}_v$ est la double classe $K_v^0 \Delta_{v,i} K_v^0$ de la matrice diagonale $\Delta_{v,i}$ ayant $i$ coefficients égaux à $\varpi_F$ et $n-i$ égaux à 1 (où $\varpi_v \in F_v$ est un générateur de l'idéal maximal de $\mathscr{O}_v$).

*2.1.3. Formes automorphes p-adiques.* — Si $\mathbb{G}$ est un groupe réductif sur [40] $\mathbf{Q}$, EMERTON [88] a introduit un objet, *la cohomologie complétée* $\widehat{H}_c^i(\mathbb{G})$, qui interpole $p$-adiquement les formes automorphes apparaissant dans la cohomologie de la tour des espaces localement symétriques $Y(K)$ associés à $\mathbb{G}$ (où $K = K_\infty K^{]\infty[}$ décrit les sous-groupes ouverts d'un sous-groupe compact maximal $\mathbb{K}$ de $\mathbb{G}(\mathbf{A_Q})$ comme ci-dessus). Si $K^{]p[}$ est un sous-groupe ouvert de $\mathbb{K} \cap \mathbb{G}(\mathbf{A_Q}^{]\infty,p[})$, on définit

$$\widehat{H}_c^i(K^{]p[}) := \varprojlim_k (\varinjlim_{K_p} H_c^i(Y(K^{]p[}K_p), \mathbf{Z}/p^k))$$

et on définit $\widehat{H}_c^i(\mathbb{G})$ comme la limite inductive des $\widehat{H}_c^i(K^{]p[})$. Une définition plus compacte [69] est $H_c^i(\mathbb{G}(\mathbf{Q}), \mathscr{C}^{(p)}(\mathbb{G}(\mathbf{A_Q}), \mathbf{Z}_p))$ où $\mathscr{C}^{(p)}(\mathbb{G}(\mathbf{A_Q}), \mathbf{Z}_p)$ est l'espace des fonctions continues, localement constantes comme fonctions de $x^{]\infty,p[}$, et l'action de $\mathbb{G}(\mathbf{Q})$ est donnée par $(\gamma \cdot \phi)(x) = \phi(\gamma^{-1}x)$.

Alors $\widehat{H}_c^i(\mathbb{G})$ est muni d'une action de $\mathbb{G}(\mathbf{A})$ induite par l'action $(g \star \phi)(x) = \phi(xg)$ sur les fonctions. CALEGARI, EMERTON [42, conj. 1.5] ont fait un certain nombre de conjectures concernant les $\widehat{H}_c^i(\mathbb{G})$. En particulier, si $\ell_0$ est la quantité apparaissant dans la discussion dans les n<sup>os</sup> 1.3.4 et 1.3.5 ($\ell_0 = 0$ si les $Y(K)$ peuvent être munis d'une structure de variété algébrique), et si $q_0 = \frac{1}{2}(d - \ell_0)$, où $d$ est la dimension des variétés différentielles $Y(K)$, alors on devrait avoir $\widehat{H}_c^i(\mathbb{G}) = 0$ si $i > q_0$, et $\widehat{H}_c^i(\mathbb{G})$ « petit » si $i < q_0$ (en résumé, le seul groupe vraiment intéressant devrait être $\widehat{H}_c^{q_0}(\mathbb{G})$).

Dans le cas $\mathbb{G} = \mathbf{GL}_2$, on a $\ell_0 = 0$ et $q_0 = 1$, et EMERTON [89] a analysé complètement la structure du $\mathbb{G}(\mathbf{A_Q})$-module $\widehat{H}_c^1(\mathbb{G})$. Un sous-produit de ses résultats [89, th. 1.2.1] est le suivant. Soit $L$ une extension finie de $\mathbf{Q}_p$, et soit $\rho \colon \mathrm{Gal}_\mathbf{Q} \to \mathbf{GL}_2(L)$ une représentation impaire, non ramifiée en dehors de $S$ fini, et dont la réduction $\overline{\rho}$ est « générique » ; pour tout $\ell \neq p$, soit $\pi_\ell(\rho)$ la $L$-représentation de $\mathbf{GL}_2(\mathbf{Q}_\ell)$ associée à la restriction de $\rho$ à $\mathrm{Gal}_{\mathbf{Q}_\ell}$ par la correspondance de Langlands locale. Alors

$$\mathbf{m}(\rho) := \mathrm{Hom}_{\mathbb{G}(\mathbf{A}^{]\infty,p[})}(\otimes'_{\ell \neq p} \pi_\ell(\rho), L \otimes_{\mathbf{Z}_p} \widehat{H}_c^1(\mathbb{G})) \neq 0$$

et donc $\rho$ est *p-adiquement modulaire* puisque $\otimes'_{\ell \neq p} \pi_\ell(\rho)$ est un facteur d'une représentation automorphe $p$-adique. EMERTON va plus loin en identifiant $\mathbf{m}(\rho)$ en termes de la correspondance de Langlands locale $p$-adique pour [41] $\mathbf{GL}_2(\mathbf{Q}_p)$. On peut se demander si l'énoncé ci-dessus ne pourrait pas être vrai en général : si $\rho\colon \mathrm{Gal}_F \to \mathbf{GL}_n(L)$ est impaire (au sens que $|\mathrm{Tr}(\rho(\sigma))| \leq 1$ si $\sigma$ est une conjugaison complexe), non ramifiée en dehors de $S$ fini ; a-t-on

$$\mathrm{Hom}_{\mathbb{G}(\mathbf{A}^{]\infty,p[})}(\otimes'_{v\nmid p}\pi_v(\rho), L \otimes_{\mathbf{Z}_p} \widehat{H}_c^{q_0}(\mathbb{G})) \neq 0\,?$$

---

39. Le facteur $|v|^{i(i-1)/2}$ dépend d'une normalisation ; il n'est pas garanti que cette normalisation soit compatible avec celle de la note 52.

40. Ceci n'est pas vraiment une restriction : si $\mathbb{G}$ est défini sur $F$, sa restriction des scalaires $\mathbb{G}'$ à $\mathbf{Q}$ est définie sur $\mathbf{Q}$, et on a $\mathbb{G}'(\mathbf{A_Q}) = \mathbb{G}(\mathbf{A}_F)$.

(On prend pour $\mathbb{G}$ la restriction des scalaires de $F$ à $\mathbf{Q}$ de $\mathbf{GL}_n$.)

Dans le cas où les $Y(K)$ sont des variétés algébriques (variétés de Shimura), elles sont définies sur des corps de nombres bien déterminés, et l'isomorphisme entre cohomologies de Betti et étale (à coefficients finis) munit les $\widehat{H}_c^i(\mathbb{G})$ d'une action galoisienne qui est beaucoup plus riche que ce qui provient d'un niveau fini [42]. C'est en particulier le cas si $\mathbb{G} = \mathbf{Sp}_{2n}$ où l'espace symétrique associé $\mathbb{G}(\mathbf{R})/K_\infty$ s'identifie au demi-espace de Siegel $\mathscr{H}_n$, sous-espace de $\mathbf{M}_n(\mathbf{C})$ des matrice symétriques dont la partie imaginaire est définie positive. Les quotients $Y(K)$ ont des compactifications minimales $Y^*(K)$ (de Satake–Baily–Borel).

Un résultat qui joue un grand rôle dans la preuve de SCHOLZE des résultats de la rem. 2.6 est que, si on fixe un sous-groupe ouvert compact $K^{]p[}$ de $\mathbb{G}(\widehat{\mathbf{Z}}^{]p[})$, la limite projective complétée $\widehat{Y}^*(K^{]p[})$ des $Y^*(K^{]p[}K_p)$, où $K_p$ décrit les sous-groupes ouverts de $\mathbb{G}(\mathbf{Z}_p)$, est une variété perfectoïde [191, th. III.3.18]. Un autre résultat fondamental est le lien suivant entre la cohomologie complétée et la cohomologie de la tour complétée.

THÉORÈME 2.7. — [191, th. IV.2.1]. *Il existe un presque-isomorphisme naturel*[43] $\mathscr{O}_{\mathbf{C}_p}\widehat{\otimes}_{\mathbf{Z}_p}\widehat{H}_c^i(K^{]p[}) \sim H^i(\widehat{Y}^*(K^{]p[}), \mathscr{I})$, *où $\mathscr{I}$ est le faisceau d'idéaux définissant le bord* (fermé complémentaire de $\widehat{Y}(K^{]p[})$).

Ce résultat fournit une description du $\mathbb{G}(\mathbf{Q}_p)$-module $\widehat{H}_c^i(K^{]p[})$ (après extension des scalaires à $\mathscr{O}_{\mathbf{C}_p}$ ou $\mathbf{C}_p$), en termes de cohomologie cohérente.

Si $n = 1$, cas des courbes modulaires, PAN [176, 177] en a déduit une description des vecteurs localement analytiques tués par l'action infinitésimale du sous-groupe unipotent supérieur, en termes de formes modulaires surconvergentes ; il a aussi montré que l'opérateur de Sen (cf. n° 4.2.5) sur cet espace est donné par l'action infinitésimale de $\begin{pmatrix} \mathbf{Q}_p^* & 0 \\ 0 & 1 \end{pmatrix}$, et en a déduit une nouvelle preuve, très prometteuse, de la conj. 1.22 dans le cas non régulier (poids de Hodge–Tate 0 et 0).

---

41. Et utilise les propriétés de cette dernière pour prouver la conj. 1.22, pour $n = 2$ et $F = \mathbf{Q}$, dans la plupart des cas.

42. Par exemple, les représentations associées aux variétés abéliennes de dimension 2 ne peuvent pas apparaître en niveau fini de la tour correspondant à $\mathbf{GSp}_4$, mais elles apparaissent dans la cohomologie complétée ; idem pour les représentations impaires de $\mathrm{Gal}_{\mathbf{Q}}$, d'image finie, qui apparaissent dans la cohomologie complétée des courbes modulaires, mais pas en niveau fini.

43. Ou un isomorphisme de presque-$\mathscr{O}_{\mathbf{C}_p}$-modules : la catégorie des presque-$\mathscr{O}_{\mathbf{C}_p}$-modules est la catégorie des $\mathscr{O}_{\mathbf{C}_p}$-modules quotientée par celle des modules presque nuls, i.e. tués par $p^r$, pour tout $r > 0$ ; dans cette catégorie, on a $\mathfrak{m}_{\mathbf{C}_p} \sim \mathscr{O}_{\mathbf{C}_p}$ par exemple.

Les résultats de PAN sur les vecteurs localement analytiques de la cohomologie complétée ont été généralisés en dimension supérieure par RODRÍGUEZ CAMARGO [186]. Cette généralisation est utilisée dans [29] pour donner une description géométrique de l'opérateur de Sen, ce qui permet d'exploiter sa semi-simplicité dans le cas de la représentation associée à une variété abélienne (c'est la semi-simplicité de l'opérateur de Sen qui permet de différencier les représentations géométriques des autres représentations à poids de Hodge–Tate généralisés $0, 0, 1, 1$, apparaissant dans la cohomologie complétée). L'approche utilisée dans [28] est différente : elle repose sur l'étude de la cohomologie cohérente du complémentaire du lieu supersingulier de $Y(K)$, pour un $K$ bien choisi (pas de montée en niveau infini), et sur la construction d'idempotents (à la HIDA [137]) découpant un sous-complexe avec de bonnes propriétés (théorie de Hida supérieure, introduite par PILLONI [182]).

## 2.2. Représentations modulo $p$

*2.2.1. Modularité potentielle modulo p.* — La preuve de la modularité potentielle de représentations modulo $p$ est une variante de l'échange 3-5 de WILES [222]. Rappelons ce dont il s'agit : on cherche à prouver qu'une courbe elliptique $E$ sur $\mathbf{Q}$ est modulaire (i.e. que le système $(\rho_{E,p})_p$ de ses modules de Tate l'est). On sait que $\overline{\rho}_{E,3}$ est modulaire grâce aux résultats de LANGLANDS et TUNNELL, mais la représentation $\overline{\rho}_{E,3}$ n'est pas forcément assez grosse pour pouvoir appliquer le théorème de relèvement modulaire 2.14. Si $\overline{\rho}_{E,3}$ n'est pas assez grosse, $\overline{\rho}_{E,5}$ l'est mais rien ne dit qu'elle soit modulaire. L'astuce est de construire une courbe elliptique $E'$ telle que $\overline{\rho}_{E',5} = \overline{\rho}_{E,5}$ et $\overline{\rho}_{E',3}$ est assez grosse : alors le th. 2.14 implique que $\rho_{E',3}$ est modulaire, et donc que le système $(\rho_{E',p})_p$ est modulaire, et donc que $\overline{\rho}_{E',5} = \overline{\rho}_{E,5}$ est modulaire. L'existence de $E'$ est due au fait que les $E'$ avec $\overline{\rho}_{E',5} = \overline{\rho}_{E,5}$ correspondent aux points rationnels de la courbe modulaire $Y$ de niveau $\overline{\rho}_{E,5}$, une forme tordue de la courbe modulaire $Y(5)$ qui se trouve être de genre 0, et donc $Y$ est isomorphe à $\mathbf{P}^1$ puisque $E$ fournit un point rationnel.

Tomber sur une courbe de genre 0 est un miracle qui se produit rarement et, en général, on ne peut pas exhiber de points rationnels sur le corps de base mais le th. 2.8 ci-dessous en fournit sur une extension $K$ de $F$ avec des propriétés locales imposées à condition qu'il n'y ait pas d'obstruction locale à leur existence (mais on ne contrôle pas le degré de l'extension $K/F$ par exemple).

THÉORÈME 2.8. — (MORET-BAILLY [168]). *Soient $F$ un corps de nombres, $S$ un ensemble fini de places de $F$, et $F_S$ l'extension maximale de $F$ dans laquelle toutes les places de $S$ sont totalement scindées*[44]. *Si $X$ est une variété quasi-projective sur $F$, et si $X(F_v) \neq \varnothing$ pour toute $v \in S$, alors $X(F_S)$ est Zariski-dense dans $X$.*

Soit maintenant $\overline{\rho}_p : \mathrm{Gal}_F \to \mathbf{GL}_n(\overline{\mathbf{F}}_p)$ dont on veut prouver la modularité potentielle, et soit $\overline{\rho}_q$, pour $q \neq p$ dont on sait prouver la modularité ainsi que celle de toutes ses restrictions aux extensions de $F$, ce qui, en l'état actuel de nos connaissances, impose de prendre $\overline{\rho}_q$ de la forme $\mathrm{I}_F^K \chi$ avec $K/F$ cyclique de degré $n$, et $\chi$ bien choisi.

---

44. Par exemple $\mathbf{Q}_{\{\infty\}}$ est le composé $\mathbf{Q}^{\text{t-r}}$ de tous les corps totalement réels.

Pour prouver la modularité de $\overline{\rho}_p$, on construit un système compatible $(\rho_\ell)_\ell$ tel que la réduction modulo $p$ de $\rho_p$ soit $\overline{\rho}_p$ et celle modulo $q$ soit $\overline{\rho}_q$ (c'est cette étape qui demande de remplacer $F$ par une extension non explicite $K$) ; ensuite on déduit la modularité de $\rho_q$ de celle de $\overline{\rho}_q$ via un théorème de relèvement modulaire, cette modularité implique celle de $\rho_p$ et donc aussi celle de $\overline{\rho}_p$ (restreinte à $\mathrm{Gal}_K$).

Pour cela, on fabrique une famille de motifs paramétrés par une variété algébrique $Y$ définie sur $F$ dont les réalisations étales modulo $p$ et $q$ sont $\overline{\rho}_p$ et $\overline{\rho}_q$ en tout point (i.e. si $y \in Y(K)$ où $K$ est une extension finie de $F$, et si $M_y$ est le motif correspondant, la réduction de $M_{y,p}$ modulo $p$ (resp. $q$) est $\mathrm{R}_F^K \overline{\rho}_p$ (resp. $\mathrm{R}_F^K \overline{\rho}_q$)). Dans le cas de Wiles $Y$ est une (forme tordue) de courbe modulaire ; en dimension $n$ on part de la famille d'hypersurfaces de Dwork d'équation $X_1^N + \cdots + X_N^N = NtX_1 \cdots X_N$ (si $t$ est fixé, cela définit une hypersurface de $\mathbf{P}^{N-1}$, de degré $N$), où $(N,n) = 1$ et $N$ est assez grand[45] par rapport à $n$. On utilise le groupe des automorphismes $\mathrm{Ker} : (\boldsymbol{\mu}_N)^N \to \boldsymbol{\mu}_N$ de cette famille pour fabriquer des systèmes locaux de dimension $n$ sur $Y_0 = \mathbf{P}^1 \setminus \{0, \boldsymbol{\mu}_N, \infty\}$, et alors $Y$ est un revêtement fini de $Y_0$ adéquat obtenu en imposant une structure de niveau $\overline{\rho}_p \times \overline{\rho}_q$ (cela impose des restrictions à l'image de $\overline{\rho}_p \times \overline{\rho}_q$).

La description précédente est très idéalisée ; la mouture ci-dessous du th. 2.8 donne une idée des acrobaties qu'il faut accomplir pour faire fonctionner la machine.

PROPOSITION 2.9. — [27, prop. 4.5.1]. *Soit $F$ un corps CM imaginaire, galoisien sur $\mathbf{Q}$, et soit $X$ une variété algébrique lisse, géométriquement connexe, définie sur $F$. On se donne :*
- *une extension finie $K$ de $F$,*
- *un ensemble fini $S_0$ de nombres premiers,*
- *pour tout $\ell \in S_0$ et toute place $v \mid \ell$ de $F$, une extension galoisienne finie $L_v/F_v$ avec $\sigma(L_v) = L_{\sigma(v)}$ pour $\sigma \in \mathrm{Gal}_{\mathbf{Q}_\ell}$, et un ouvert non vide $\Omega_v \subset X(L_v)$ stable par $\mathrm{Gal}(L_v/F_v)$.*

*Alors on peut trouver une extension galoisienne $F'$ de $F$, CM, avec $F' \cap K = F$, et $P \in X(F')$ tels que, pour toute place $w \mid v \mid \ell$ de $F'$, on ait $F'_w \cong L_v$ et $P \in \Omega_v \subset X(L_v) = X(F'_w)$.*

*De plus, si $G$ est un groupe fini, et si $f : \pi_1^{\mathrm{ét}}(X) \to G$ est un morphisme surjectif, on peut choisir $P$ de telle sorte que $f \circ P_* : \mathrm{Gal}_{F'} \to G$ soit surjectif.*

*2.2.2. Grosse image.* — La plupart des théorèmes de relèvement modulaire demande que la représentation que l'on cherche à relever soit d'image assez grosse, en particulier que la représentation soit irréductible (une exception est le th. 2.15 ci-dessous). Si on

---

45. L'introduction de $N$ permet [8, 27] de considérer des $\overline{\rho}_p$ à valeurs dans $\mathbf{GL}_n(\mathbf{F}_{p^k})$ ; si on ne s'intéresse qu'à des $\overline{\rho}_p$ à valeurs dans $\mathbf{GL}_n(\mathbf{F}_p)$, on peut prendre $N = n+1$ et considérer le système local fourni par le $H^{n-1}$ de la famille [128].

cherche à prouver la modularité d'un système compatible $(\rho_p)_p$, on a donc besoin de savoir que $\overline{\rho}_p$ a une grosse image pour au moins un $p$ (et préférablement plusieurs $p$ pour avoir un peu de flexibilité pour remplir d'autres conditions éventuelles). Dans cette direction, on a le résultat suivant :

THÉORÈME 2.10. — [7, prop. 5.3.2]. *Soit $\rho = (\rho_p)_p$ un système compatible de représentations de $\mathrm{Gal}_F$, irréductible, régulier. Alors pour $p$ dans un ensemble de densité de Dirichlet* 1, *la restriction de $\overline{\rho}_p$ à $\mathrm{Gal}_{F(\boldsymbol{\mu}_p)}$ est irréductible.*

Voir aussi [180, th. 1.2] pour une version sans l'hypothèse de régularité.

### 2.3. Relèvements de représentations globales

*2.3.1. Densité des déformations globales dans les déformations locales.* — Soit $F$ une extension de degré $n$ de $\mathbf{Q}$, et soit $\overline{\rho}\colon \mathrm{Gal}_F \to \mathbf{GL}_d(k)$ une représentation continue, où $k$ est une extension finie de $\mathbf{F}_p$, avec $p \neq 2$. Si $L$ est une extension finie de $\mathbf{Q}_p$ de corps résiduel $k$, on cherche à relever $\overline{\rho}$ en $\rho\colon \mathrm{Gal}_F \to \mathbf{GL}_d(\mathcal{O}_L)$ en imposant certaines conditions locales. Pour éviter les obstructions liées au déterminant, on fixe un relèvement $\chi\colon \mathrm{Gal}_F \to \mathcal{O}_L^*$ de $\det \overline{\rho}$.

Soit $S$ un ensemble fini de places de $F$ contenant celles divisant $p$ et $\infty$, en dehors desquelles $\overline{\rho}$ est non ramifiée. L'espace $\mathscr{X}_{\overline{\rho}}^{\chi}$ des déformations de $\overline{\rho}$, de déterminant $\chi$, non ramifiées en dehors de $S$ est, conjecturalement [165, prop. 5], de dimension

$$\delta_{\overline{\rho}}^{\chi} := \tfrac{n}{2}d^2 - \tfrac{1}{2}\sum_{v\,\text{réelle}} t_v^2 - r_2$$

où $r_2$ est le nombre de places complexes de $F$, et $t_v = \dim_k \overline{\rho}^+ - \dim_k \overline{\rho}^-$, où $\overline{\rho}^+$ et $\overline{\rho}^-$ sont les espaces propres de la conjugaison complexe de $\mathrm{Gal}(\overline{F}_v/F_v)$ pour les valeurs propres $1$ et $-1$. Cette dimension est donc maximale quand les $|t_v|$ sont minimaux et on dit que $\overline{\rho}$ est *impaire* si $|t_v| \leq 1$ pour toute place $v$ réelle.

Maintenant, il résulte du th. 4.21 que $\prod_{v|p} \mathscr{X}_{\overline{\rho}_v}^{\chi_v}$ est de dimension

$$\sum_{v|p}([F_v : \mathbf{Q}_p](d^2 - 1)) = n(d^2 - 1)$$

et on s'attend à ce que les fibres de l'application $\mathscr{X}_{\overline{\rho}}^{\chi} \to \prod_{v|p} \mathscr{X}_{\overline{\rho}_v}^{\chi_v}$ soient de dimension 0. Si $Z$ est une « sous-variété » de codimension $\delta_Z$ de $\prod_{v|p} \mathscr{X}_{\overline{\rho}_v}^{\chi_v}$, on peut donc espérer qu'il soit possible de relever $\overline{\rho}$ en $\rho\colon \mathrm{Gal}_F \to \mathbf{GL}_d(\mathcal{O}_L)$ telle que $(\rho_v)_{v|p} \in Z$ si $\delta_Z \leq \delta_{\overline{\rho}}^{\chi}$ et penser qu'il soit, en général, impossible de le faire si $\delta_Z > \delta_{\overline{\rho}}^{\chi}$.

Par exemple, si $Z$ est l'espace des représentations qui sont triangulines (cf. n° 4.3.3) pour tout $v \mid p$, alors $\delta_Z = \sum_{v|p}[F_v : \mathbf{Q}_p]\frac{d(d-1)}{2} = n\frac{d(d-1)}{2}$. Si $\overline{\rho}$ est impaire, $\delta_Z \leq \delta_{\overline{\rho}}^{\chi}$ et on peut donc espérer que $\overline{\rho}$ ait un relèvement triangulin.

Si $Z$ est l'espace des représentations potentiellement semi-stables de type fixé, alors $\delta_Z = n\left(\frac{d(d+1)}{2} - 1\right)$ (si $Z$ n'est pas vide) et, si $d \geq 2$, on a $\delta_Z > \delta_{\overline{\rho}}^{\chi}$ sauf si $d = 2$, $F$ est totalement réel et $\overline{\rho}$ est impaire. On s'attend donc à ce que $\overline{\rho}$ n'ait, en général, pas de relèvement potentiellement semi-stable sauf dans le cas particulier $d = 2$, $F$ est totalement réel et $\overline{\rho}$ est impaire. Dans le cas $F = \mathbf{Q}$ et $d = 2$, l'existence de relèvements potentiellement semi-stables de type imposé est contenu dans la « conjecture $\varepsilon$ » de Serre [195].

Théorème 2.11. — (Fakhruddin, Khare, Patrikis [95, th. 1.2]).

*Si $p$ est assez grand*[46]*, si $\overline{\rho}$ est impaire et si $\mathrm{R}_F^{F(\boldsymbol{\mu}_p)}\overline{\rho}$ est absolument irréductible, quitte à agrandir $L$, on peut trouver $S' \supset S$ fini, tel que $\overline{\rho}$ admette un relèvement triangulin $\rho\colon \mathrm{Gal}_F \to \mathbf{GL}_d(\mathscr{O}_L)$, non ramifié en dehors de $S'$, et on peut imposer de plus que $\rho(\mathrm{Gal}_F)$ contienne un sous-groupe ouvert de $\mathbf{SL}_d(\mathscr{O}_L)$.*

Le lecteur trouvera dans [95] un énoncé plus général (th. 3.8) qui aboutit à la même conclusion sous des hypothèses plus faibles (par exemple, en permettant à $\overline{\rho}$ d'être réductible) mais plus pénibles à énoncer.

Théorème 2.12. — (Fakhruddin, Khare, Patrikis [93, th. B]).

*Si $p$ est assez grand*[47] *si $\mathrm{R}_F^{F(\boldsymbol{\mu}_p)}\overline{\rho}$ est absolument irréductible et si, pour tout $v \in S$ on se donne un relèvement $\rho_v^0$ de $\overline{\rho}_v$, avec $\det \rho_v^0 = \chi_v$, alors, quitte à agrandir $L$, il existe, pour tout $n \geq 1$, $S^n \supset S$ fini et un relèvement $\rho^n\colon \mathrm{Gal}_F \to \mathbf{GL}_d(\mathscr{O}_L)$ de déterminant $\chi$, non ramifié en dehors de $S^n$, tel que $\rho_v^n$ soit $\mathbf{SL}_d(\mathscr{O}_L)$-conjuguée modulo $p^n$ à $\rho_v^0$, pour tout $v \in S$.*

*Remarque 2.13.* — (i) En faisant tendre $n$ vers $\infty$ dans l'énoncé du théorème, on obtient un relèvement dont les restrictions aux $F_v$, pour $v \in S$, sont les $\rho_v^0$ que l'on s'est fixés ; ce relèvement est, en général, ramifié en un nombre infini de places. On peut reformuler ce résultat en disant que l'espace des déformations de $\overline{\rho}$ (sans restriction sur la ramification) est dense dans le produit des espaces des déformations des $\overline{\rho}_v$, pour $v$ décrivant l'ensemble des places finies de $F$.

(ii) Fakhruddin, Khare et Patrikis démontrent des résultats beaucoup plus généraux et aussi plus précis que les énoncés ci-dessus : ils considèrent ([93, th. 6.21], [94, th. 5.6], [95, th. 3.6]) le problème de relever $\overline{\rho}\colon \mathrm{Gal}_F \to \mathbb{G}(k)$ en $\rho\colon \mathrm{Gal}_F \to \mathbb{G}(\mathscr{O}_L)$ où $\mathbb{G}$ est un groupe réductif quelconque (pas seulement $\mathbb{G} = \mathbf{GL}_d$) ; avec la bonne définition[48] de « impair », ils produisent des relèvements géométriques (au sens de la rem. 1.17), si $F$ est totalement réel (cf. [93, th. A], [94, th. E]).

*2.3.2. Relèvements modulaires.* — Dans ce n°, nous énonçons certains des théorèmes de relèvement modulaire les plus emblématiques que l'on peut trouver dans la littérature.

• *Les premiers théorèmes de relèvement.* — À tout seigneur, tout honneur, commençons par le théorème de relèvement modulaire de Wiles [222, th. 0.2] :

---

46. D'après les auteurs, il suffit que $[F(\boldsymbol{\mu}_p):F] > d\,f(d)^2$, où $f$ est une fonction explicite ($f(d) = (d+1)!$ si $d \geq 71$ et $f(d) = 60^d d!$, si $20 \leq d \leq 70$, etc.).

47. D'après les auteurs, la même borne que pour le th. 2.11 marche.

48. Correspondant à la condition $\ell_0 = 0$ des n°s 1.3.3, 1.3.4, 1.3.5 ; en particulier, avec cette « bonne » définition, il n'y a de représentations impaires pour $\mathbf{GL}_d$ que si $d = 2$ et si $F$ est totalement réel.

THÉORÈME 2.14. — *Soit $\rho\colon \mathrm{Gal}_\mathbf{Q} \to \mathbf{GL}_2(\overline{\mathbf{Q}}_p)$ une représentation continue, impaire, non ramifiée en dehors d'un ensemble fini de places, dont la restriction à $\mathrm{Gal}_{\mathbf{Q}_p}$ est ordinaire à poids $0$ et $k \geq 1$ ou cristalline à poids $0$ et $1$. On suppose que $\overline{\rho}$ est modulaire, que $\mathrm{R}_\mathbf{Q}^{\mathbf{Q}(\boldsymbol{\mu}_p)}\overline{\rho}$ est absolument irréductible, et que, si $q \equiv -1 \mod p$ est tel que $\overline{\rho}$ est ramifié en $q$, alors la restriction de $\overline{\rho}$ à $\mathrm{Gal}_{\mathbf{Q}_q}$ est soit absolument irréductible en restriction au sous-groupe d'inertie, soit une somme de deux caractères. Alors $\rho$ est modulaire.*

Le résultat précédent est le point clé de la preuve de la conjecture de Taniyama–Weil pour les courbes elliptiques semi-stables sur $\mathbf{Q}$. Le suivant a joué un grand rôle dans la preuve de la conjecture de Serre par KHARE, WINTENBERGER (qui ont aussi prouvé leur propre théorème de relèvement modulaire [143, th. 4.1]).

THÉORÈME 2.15. — (SKINNER–WILES [202]). *Soit $\rho\colon \mathrm{Gal}_\mathbf{Q} \to \mathbf{GL}_2(\overline{\mathbf{Q}}_p)$ une représentation continue impaire, non ramifiée en dehors d'un ensemble fini de places et ordinaire en $p$. On suppose que $\overline{\rho}$ n'est pas irréductible mais que la restriction de $\overline{\rho}^{\mathrm{ss}}$ à $\mathrm{Gal}_{\mathbf{Q}_p}$ n'est pas de la forme $\chi \oplus \chi$, alors $\rho$ est modulaire.*

• *Le cas autodual impair.* — Une représentation potentiellement cristalline $\rho$ de $\mathrm{Gal}_K$ est *potentiellement diagonalisable* s'il existe une extension $K'$ de $K$ telle que la composante irréductible de l'espace des déformations cristallines de $\mathrm{R}_K^{K'}\overline{\rho}$ de poids de Hodge–Tate fixés contenant $\rho$ contient aussi une somme de caractères.

Il est possible que toute représentation potentiellement cristalline soit potentiellement diagonalisable (en tout cas, on n'a pas de contrexemple). Le résultat suivant est un ingrédient clé dans la preuve du th. 1.24.

THÉORÈME 2.16. — [7, th. C] *Soit $F$ un corps totalement réel. Soit $p \geq 2n+1$ et soit $\rho\colon \mathrm{Gal}_F \to \mathbf{GL}_n(\overline{\mathbf{Q}}_p)$ une représentation continue, non ramifiée presque partout, autoduale, potentiellement diagonalisable à poids de Hodge–Tate réguliers, telle que $\mathrm{R}_F^{F(\boldsymbol{\mu}_p)}\overline{\rho}$ soit irréductible. Alors $\rho$ est potentiellement modulaire.*

• *Le cas non autodual.* — Le résultat suivant (et ses variantes mentionnées dans la rem. 2.18) est un des ingrédients principaux dans la preuve du th. 1.26.

THÉORÈME 2.17. — ([2, th. 6.1.1]). *Soit $F$ un corps CM et soit $p$ un nombre premier non ramifié dans $F$. Soit $\rho\colon \mathrm{Gal}_F \to \mathbf{GL}_n(\overline{\mathbf{Q}}_p)$ une représentation non ramifiée en presque toute place de $F$ et cristalline en les places divisant $p$, à poids de Hodge–Tate réguliers. On suppose que :*

• *$p$ est assez grand par rapport à $n$ et aux poids de Hodge–Tate de $\rho$,*

• *$\overline{\rho}$ est irréductible, assez générique et d'image assez grosse,*

• *il existe une représentation automorphe cuspidale algébrique $\pi$ de caractère infinitésimal correspondant aux poids de Hodge–Tate de $\rho$, telle que $\overline{\rho}_\pi \cong \overline{\rho}$.*

*Alors $\rho$ est modulaire et compatible à la correspondance de Langlands locale pour toute $v \mid p$ et toute place $v$ en laquelle $\pi$ et $\rho$ sont non ramifiées.*

*Remarque 2.18.* — (i) Il y a une variante [2, th. 6.1.2], où on remplace la condition de cristallinité en les places divisant $p$ par celle d'ordinarité (sans supposer que $p$ est non ramifié dans $F$) : on demande juste à $\pi$ d'être de caractère infinitésimal régulier (sans imposer de lien avec les poids de Hodge–Tate de $\rho$), et on garde les autres conditions. Alors $\rho$ est modulaire et compatible à la correspondance de Langlands locale pour toute place $v$ en laquelle $\pi$ et $\rho$ sont non ramifiées.

(ii) Dans [27, th. 3.2.1], les auteurs énoncent une autre variante du th. 2.17, dans laquelle on autorise $p$ à être arbitrairement ramifié mais on suppose que les poids de Hodge–Tate de $\rho$ sont $0, 1, \ldots, n-1$ pour tout $\tau \in \mathrm{Hom}(F, \overline{\mathbf{Q}})$ et on impose des liens un peu plus forts entre $\rho_\pi$ et $\rho$ que dans le th. 2.17.

(iii) Un ingrédient fondamental dans la preuve du th. 2.17 est le théorème d'annulation de CARAIANI, SCHOLZE [46] (voir aussi [155, 125]).

- *Un cas non régulier.* — Pour finir, voici un résultat [49] un peu différent puisqu'on construit une forme automorphe pour $\mathbf{GSp}_4$ (que l'on peut transférer à $\mathbf{GL}_4$, si besoin).

THÉORÈME 2.19. — [28, th. 1.1.7]. *Soit $F$ un corps totalement réel dans lequel $p > 2$ se scinde complètement. Soit $A$ une surface abélienne sur $F$ avec bonne réduction ordinaire en toute $v \mid p$. Soit $\rho_{A,p}$ la représentation associée à $H^1_{\mathrm{ét}}(A_{\overline{F}}, \mathbf{Z}_p)$. On suppose que :*

- $\overline{\rho}_{A,p}$ *est d'image assez grosse et assez générique,*
- *il existe une représentation automorphe $\pi$ de $\mathbf{GSp}_4/F$ de poids 2, caractère central $\mid \mid^2_{\mathbf{A}}$, ordinaire en les $v \mid p$, telle que $\overline{\rho}_{\pi,p} = \overline{\rho}_{A,p}$ et $\rho_\pi$ vérifie des conditions locales de pureté.*

*Alors $\rho_{A,p}$ est la représentation associée à une forme modulaire de Siegel de poids 2.*

## 3. BOITE À OUTILS AUTOMORPHE

Les résultats de ce chapitre interviennent à plusieurs niveaux dans la preuve de l'existence de prolongements analytiques de fonctions $[L]$ ; en particulier, les théorèmes de transfert ou de changement de base du § 3.3 permettent de circuler d'un groupe que l'on comprend à un groupe qui nous intéresse (comme ils ne sont pas vraiment aussi généraux que ce que l'on aimerait, cela oblige à certaines acrobaties).

### 3.1. Représentations automorphes de $\mathbf{GL}_n(\mathbf{A}_F)$

*3.1.1. Formes automorphes.* — Soit $F$ un corps de nombres et soit $\mathbb{G} = \mathbf{GL}_n$. Notons $\mathscr{C}^\infty(\mathbb{G})$ l'espace des $\phi \colon \mathbb{G}(\mathbf{A}_F) \to \mathbf{C}$ telles que $\phi(x_\infty, x^{]\infty[})$ soit localement constante comme fonction de $x^{]\infty[}$ et de classe $\mathscr{C}^\infty$ comme fonction de $x_\infty \in \mathbb{G}_\infty = \prod_{v \mid \infty} \mathbb{G}(F_v)$.

---

49. Voir aussi [220, th. 7.12].

On munit $\mathscr{C}^\infty(\mathbb{G})$ d'actions de $\mathbb{G}(F)$, $\mathbb{G}(\mathbf{A}_F)$ et de l'algèbre de Lie $\mathfrak{g}$ de $\mathbb{G}_\infty$ données par les formules :

$$(\gamma * \phi)(x) = \phi(\gamma^{-1} x), \text{ si } \gamma \in \mathbb{G}(F), \quad (g \star \phi)(x) = \phi(xg), \text{ si } g \in \mathbb{G}(\mathbf{A}_F)$$

$$(X \star \phi)(x) = \left(\tfrac{d}{dt}\phi\bigl(xe^{tX}\bigr)\right)_{|t=0}, \text{ si } X \in \mathfrak{g}$$

On choisit un sous-groupe compact maximal $K_\infty$ de $\mathbb{G}_\infty$ : on a $\mathbb{G}_\infty \cong \mathbb{G}(\mathbf{R})^{r_1} \times \mathbb{G}(\mathbf{C})^{r_2}$ et $K_\infty$ est conjugué à $\mathbf{O}_n^{r_1} \times \mathbf{U}_n^{r_2}$ où $\mathbf{O}_n$ et $\mathbf{U}_n$ sont les groupes orthogonal et unitaire. On choisit aussi un sous-groupe compact maximal $K^{]\infty[}$ de $\mathbb{G}(\mathbf{A}^{]\infty[})$ (on peut prendre pour $K^{]\infty[}$ le produit des $\mathbb{G}(\mathscr{O}_v)$, pour $v$ finie), et on pose $K = K_\infty \times K^{]\infty[}$.

L'espace $\mathrm{A}(\mathbb{G})$ des *formes automorphes* est l'espace des $\phi \in \mathscr{C}^\infty(\mathbb{G})$ vérifiant les conditions suivantes :
- ($K$-finitude) L'espace engendré par les $g \star \phi$, pour $g \in K$, est de dimension finie.
- (Finitude de l'action du centre) : L'espace engendré par les $z \star \phi$ et les $Z \star \phi$, pour [50] $z \in Z(\mathbb{G}(\mathbf{A}_F))$ et $Z \in Z(U(\mathfrak{g}))$, est de dimension finie.
- $\phi$ est « à croissance lente » à l'infini.
- (Automorphie) $\phi$ est fixe par l'action de $\mathbb{G}(F)$.

On dit que $\phi \in \mathrm{A}(\mathbb{G})$ est *cuspidale* si $\phi$ est « nulle à l'infini », ce qui se traduit par $\int_{\mathbb{U}(F)\backslash\mathbb{U}(\mathbf{A}_F)} \phi(u)\,du = 0$ pour tout sous-groupe unipotent $\mathbb{U}$ de $\mathbb{G}$. On note $\mathrm{A}^0(\mathbb{G})$ le sous-espace de $\mathrm{A}(\mathbb{G})$ des formes cuspidales.

*3.1.2. Représentations automorphes.* — L'espace $\mathrm{A}(\mathbb{G})$ n'est pas stable par l'action de $\mathbb{G}_\infty$ à cause du choix de $K_\infty$ (par contre, comme deux sous-groupes ouverts compacts de $\mathbb{G}(\mathbf{A}^{]\infty[})$ sont commensurables, $\mathrm{A}(\mathbb{G})$ est stable par $\mathbb{G}(\mathbf{A}^{]\infty[})$). L'action de $\mathbb{G}_\infty$ est remplacée par une structure de $(\mathfrak{g}, K_\infty)$-module (les actions de $\mathfrak{g}$ et $K_\infty$ stabilisent $\mathrm{A}(\mathbb{G})$ et on a $g \star (X \star \phi) = (gXg^{-1}) \star (g \star \phi)$ si $g \in K_\infty$ et $X \in \mathfrak{g}$ et $(X \star \phi)(x) = \left(\tfrac{d}{dt}\phi\bigl(xe^{tX}\bigr)\right)_{|t=0}$, si $X \in \mathrm{Lie}(K_\infty) \subset \mathfrak{g}$).

La condition sur la finitude de l'action du centre fait que $\mathrm{A}(\mathbb{G})$ se décompose comme une somme directe de sous-espaces sur lesquels le centre agit par un caractère. Un tel caractère est donné par :
- sa restriction $\psi$ à $Z(\mathbb{G}(\mathbf{A}_F))$ (le *caractère central*, qui est un caractère lisse de $\mathbf{A}_F^*$, automorphe, i.e. qui se factorise à travers $\mathbf{A}_F^*/F^*$),
- pour tout $\tau \in \mathrm{Hom}(F, \overline{\mathbf{Q}})$, la donnée d'un caractère du centre de l'algèbre enveloppante (le *caractère infinitésimal*) de l'algèbre de Lie $\mathbf{M}_n(\mathbf{C})$ (avec $[X,Y] = XY - YX$) ; un tel caractère est lui-même encodé par la donnée de $n$ nombres complexes à permutation près (i.e par un polynôme unitaire de degré $n$) ; on note $H_\tau$ le multi-ensemble [51] de cardinal $n$ correspondant [52].

---

50. $Z(\mathbb{G}(\mathbf{A}_F))$ désigne le centre de $\mathbb{G}(\mathbf{A}_F)$ et $Z(U(\mathfrak{g}))$ celui de l'algèbre enveloppante universelle $U(\mathfrak{g})$ de $\mathfrak{g}$.

51. Un multi-ensemble est un ensemble où les éléments ont une multiplicité ; si $m_x$ est la multiplicité de $x$, on note le multi-ensemble sous la forme $\sum_x m_x\{x\}$ ; son cardinal est $\sum_x m_x$.

52. Pour fixer les idées, si $\lambda_1 \geq \lambda_2 \geq \cdots \geq \lambda_n$ sont des entiers, et si $V$ est la représentation algébrique de $\mathbb{G}$ obtenue en induisant le caractère $g \mapsto a_{1,1}^{\lambda_1} \cdots a_{n,n}^{\lambda_n}$ du Borel inférieur, les $n$ nombres

Si $H = (H_\tau)_\tau$ est la collection de ces multi-ensembles de nombres complexes, on note $\mathrm{A}_{\psi,H}(\mathbb{G})$ et $\mathrm{A}^0_{\psi,H}(\mathbb{G})$ les sous-espaces correspondants.

Une *représentation automorphe* de $\mathbb{G}$ est un sous-$\mathbb{G}(\mathbf{A}_F^{]\infty[}) \times (\mathfrak{g}, K_\infty)$-module irréductible de $\mathrm{A}^0(\mathbb{G})$. Une représentation automorphe $\pi$ est *supercuspidale* si $\pi \subset \mathrm{A}^0(\mathbb{G})$. Si $\pi$ est une représentation automorphe de $\mathbb{G}$, il existe $(\psi_\pi, H_\pi = (H_{\pi,\tau})_\tau)$, unique, tel que $\pi \subset \mathrm{A}_{\psi_\pi, H_\pi}(\mathbb{G})$; on dit que $\pi$ est de type $(\psi_\pi, H_\pi)$.

Une représentation automorphe $\pi$ de $\mathbb{G}$ admet une factorisation sous la forme d'un produit tensoriel restreint $\pi \cong \otimes'_v \pi_v$ où $\pi_v$ est :

• une représentation admissible irréductible de $\mathbb{G}(F_v)$, si $v$ est une place finie, de caractère central $\psi_{\pi,v}$ (restriction de $\psi_\pi$ à $F_v^*$), non ramifiée[53] pour tout $v$ sauf un nombre fini,

• un $(\mathfrak{g}_v, K_v)$-module admissible de caractère infinitésimal encodé par $H_{\pi,\tau}$ si $v$ est une place réelle induite par $\tau \colon F \to \overline{\mathbf{Q}} \hookrightarrow \mathbf{C}$ ou par $(H_{\pi,\tau}, H_{\pi,\overline{\tau}})$ si $v$ est une place complexe induite par $\tau$ et $\overline{\tau}$.

*Exemple 3.1.* — Le cas $n=1$. Ce cas est un peu caricatural car $\mathbb{G}$ est égal à son centre, et les représentations automorphes sont donc exactement les *caractères automorphes* (aussi connus sous le nom de caractères de Hecke), i.e. les caractères continus $\eta \colon \mathbf{A}_F^* \to \mathbf{C}^*$ triviaux sur $F^*$ (par exemple, le caractère $|\ |_\mathbf{A}$ défini par $|(x_v)_v|_\mathbf{A} = \prod_v |x_v|_v$ ; la trivialité de $|\ |_\mathbf{A}$ sur $F^*$ est équivalente à la formule du produit).

Dans ce cas la factorisation prend la forme simple $\eta((x_v)_v) = \prod_v \eta_v(x_v)$, où $\eta_v \colon F_v^* \to \mathbf{C}^*$ est un caractère continu.

Si $v \mid \infty$ est réelle induite par $\tau \colon F \to \overline{\mathbf{Q}} \hookrightarrow \mathbf{C}$, la restriction de $\eta_v$ à $\mathbf{R}_+^* \subset F_v^*$ est de la forme $x \mapsto x^{\lambda_\tau}$, et $H_{\eta,\tau}$ est la valeur propre (avec multiplicité 1) de $x \frac{d}{dx}$, à savoir $\lambda_\tau$. Si $v$ est complexe, induite par $\tau$ et $\overline{\tau}$, alors $F_v = \mathbf{C}$ (identification choisie de telle sorte que $F \hookrightarrow F_v = \mathbf{C}$ coïncide avec $\iota_\infty \circ \tau$) et $\eta_v$ est de la forme $z \mapsto z^{\lambda_{\eta,\tau}} \overline{z}^{\lambda_{\eta,\overline{\tau}}}$ (et $\lambda_\tau - \lambda_{\overline{\tau}} \in \mathbf{Z}$) ; alors $H_{\eta,\tau}$ et $H_{\eta,\overline{\tau}}$ sont les valeurs propres respectives (avec multiplicité 1) de $z \frac{\partial}{\partial z}$ et $\overline{z} \frac{\partial}{\partial \overline{z}}$, à savoir $\lambda_{\eta,\tau}$ et $\lambda_{\eta,\overline{\tau}}$.

*Remarque 3.2.* — Si $\pi \subset \mathrm{A}(\mathbb{G})$ est une représentation automorphe, et si $\eta$ est un caractère automorphe, alors $\{(\eta \circ \det)\phi,\ \phi \in \pi\}$ est une représentation automorphe ; on la note $\pi \otimes \eta$. On a $\psi_{\pi \otimes \eta} = \psi_\pi \eta^n$, et $H_{\pi \otimes \eta, \tau} = H_{\pi, \eta} + \lambda_{\eta,\tau}$.

*3.1.3. Fonctions L modulaires.* — Si $\pi_v$ est non ramifiée, et si $K_v = \mathbb{G}(F_v)$, alors $\pi_v^{K_v}$ est de dimension 1, et l'algèbre (de Hecke) $\mathbf{C}[K_v \backslash \mathbb{G}(F_v)/K_v]$ des fonctions à support compact, bi-invariantes par $K_v$, agit sur $\pi_v^{K_v}$ par un caractère. Or on dispose d'un isomorphisme naturel (de Satake [187]) $\mathbf{C}[K_v \backslash \mathbb{G}(F_v)/K_v] \cong \mathbf{C}[X_1^\pm, \ldots, X_n^\pm]^{S_n}$. Un

---

complexes correspondants sont $\lambda_1 + n - 1, \lambda_2 + n - 2, \ldots, \lambda_n$ (une normalisation plus classique donnerait $\lambda_1 + \frac{n-1}{2}, \lambda_2 + \frac{n-3}{2}, \ldots, \lambda_n + \frac{1-n}{2}$). Par exemple, si $V$ est la représentation $\det^a \otimes \mathrm{Sym}^k$ de $\mathbf{GL}_2(\mathbf{R})$, les deux nombres complexes correspondants sont $a + k + 1, a$.

53. Cela signifie que $H^0(\mathbb{G}(\mathcal{O}_v), \pi) \neq 0$.

caractère de cette algèbre de Hecke est donc équivalent à la donnée de $n$ nombres complexes non nuls $\alpha_{v,1},\ldots,\alpha_{v,n}$ à permutation près. On pose alors

$$L(\pi_v, s) = \prod_{i=1}^{n}(1 - \alpha_{v,i}|v|^{-s})^{-1}$$

et $\pi_v$ est complètement décrite par le facteur d'Euler $L(\pi_v, s)$.

Si $\pi_v$ est ramifiée, il y a une recette plus compliquée fournissant $L(\pi_v, s)$, un conducteur $f(\pi_v) \in \mathbf{N}$ et une constante locale $\varepsilon(\pi_v) \in \mathbf{C}^*$. On pose alors $L(\pi, s) = \prod_v L(\pi_v, s)$, le produit portant sur les places finies.

Il y a aussi une recette faisant intervenir les $H_{\pi,\tau}$ permettant de définir un facteur $\Gamma(\pi, s)$ produit de facteurs du type $\Gamma_{\mathbf{R}}(s-a)$. On a le résultat suivant (pour $n=2$, ce résultat a été prouvé par JACQUET et LANGLANDS [139], en adélisant à la TATE les résultats classiques de HECKE [133] pour les formes modulaires).

THÉORÈME 3.3. — (GODEMENT, JACQUET [121]). *Si $\pi$ est une représentation automorphe, alors $L(\pi, s)$ converge pour $\mathrm{Re}(s) > 1$, admet un prolongement méromorphe à $\mathbf{C}$, holomorphe en dehors d'un ensemble fini de pôles*[54] *et une équation fonctionnelle de la forme*

$$\Gamma(\pi, s)\, L(\pi, s) = \varepsilon(\pi)\, N(\pi)^{-s}\, \Gamma(\check{\pi}, s)\, L(\check{\pi}, 1-s)$$

*où $\check{\pi}$ est la contragrédiente de $\pi$, $\varepsilon(\pi) = \prod_v \varepsilon(\pi_v)$ et $N(\pi) = \prod_v |v|^{f(\pi_v)}$.*

CONJECTURE 3.4. — *Si $\rho$ est un système compatible irréductible de poids de Hodge–Tate $(H_\tau)_\tau$, alors il existe $\pi$ automorphe cuspidale, telle que $L(\rho, s) \sim L(\pi, s)$; de plus*[55] *$H_{\pi,\tau} = -H_\tau$ pour tout $\tau \in \mathrm{Hom}(F, \overline{\mathbf{Q}})$.*

Compte-tenu du th. 3.3, cette conjecture, cas particulier de la correspondance de Langlands globale, implique la conjecture 1.18 sous une forme particulièrement satisfaisante. Notons que la représentation $\pi$ dont la conjecture affirme l'existence est unique d'après le résultat suivant.

THÉORÈME 3.5. — (JACQUET, SHALIKA [140]). *Soient $\pi$ et $\pi'$ des représentations cuspidales automorphes de $\mathbf{GL}_n(\mathbf{A}_F)$. Si $L(\pi, s) \sim L(\pi', s)$, alors $\pi = \pi'$ (théorème de multiplicité 1).*

### 3.2. Théorèmes inverses

On dispose de théorèmes inverses inspirés par un résultat de WEIL [219], disant que, si une fonction $L$ et suffisamment de ses tordues ont de bonnes propriétés, alors cette fonction $L$ est la fonction $L$ d'une représentation automorphe. La torsion dont il est question fait l'objet du théorème suivant.

---

54. L'ensemble des pôles est presque toujours vide ; c'est en particulier le cas si $\pi$ est cuspidale.

55. $-H_\tau = \sum m_i\{-i\}$ si $H_\tau = \sum_i m_i\{i\}$. Ce signe ennuyeux est dû à la normalisation des poids de Hodge–Tate (le caractère cyclotomique a poids 1 dans nos conventions) ou bien au fait que l'on a privilégié les frobenius géométriques $\mathrm{Frob}_v^{-1}$.

THÉORÈME 3.6. — (JACQUET, PIATETSKI-SHAPIRO, SHALIKA [141]). *Soient $\pi$ et $\pi'$ des représentations automorphes cuspidales de $\mathbf{GL}_r(\mathbf{A})$ et $\mathbf{GL}_s(\mathbf{A})$ respectivement. Alors $[L(\pi,s) \otimes L(\pi',s)]$ admet un prolongement analytique. De plus :*

• (SHAHIDI [197]). *$[L(\pi,s) \otimes L(\pi',s)]$ admet une équation fonctionnelle (et donc a un représentant privilégié $L(\pi \times \pi', s)$).*

• (MOEGLIN, WALDSPURGER [167, Appendice]). *Le prolongement de $L(\pi \times \pi', s)$ est holomorphe sauf si $r = s$ et il existe $t$ tel que $\pi = \check{\pi}' \otimes |\ |_{\mathbf{A}}^t$ où il y a un pôle simple en $s = 1 - t$.*

Les théorèmes inverses proprement dits sont résumés dans l'énoncé suivant.

THÉORÈME 3.7. — (COGDELL, PIATETSKI-SHAPIRO [63, 64]). *Soit $\Pi = \otimes'_v \Pi_v$ une représentation admissible irréductible de $\mathbb{G}(\mathbf{A}_F)$, dont le caractère central est automorphe, et soit $S$ un ensemble fini de places finies de $K$.*

(i) *Si pour tout $m \leq n - 2$, et pour toute représentation automorphe $\pi$ de $\mathbf{GL}_m(\mathbf{A}_F)$, non ramifiée pour $v \in S$, la fonction $L(\Pi \times \pi, s)$ vérifie les conclusions du th. 3.6. Alors, il existe $\Pi'$ automorphe telle que $\Pi'_v \cong \Pi_v$ pour tout $v \notin S$.*

(ii) *Si $\mathscr{O}[\frac{1}{S}]$ est principal, si pour tout $m \leq n - 1$, et pour toute représentation automorphe $\pi$ de $\mathbf{GL}_m(\mathbf{A}_F)$, non ramifiée pour $v \notin S$, la fonction $L(\Pi \times \pi, s)$ vérifie les conclusions du th. 3.6. Alors, il existe $\Pi'$ automorphe telle que $\Pi'_v \cong \Pi_v$ pour tout $v \in S$ et tout $v \notin S$ tel que $\Pi_v$ et $\Pi'_v$ sont toutes les deux non ramifiées.*

*Remarque 3.8.* — (i) COGDELL et PIATETSKI-SHAPIRO conjecturent [63] qu'il suffit de tordre par des caractères pour que la conclusion du (i) soit valide. C'est connu si $n = 2$ et $n = 3$. Combinée avec le th. 3.6, une preuve de cette conjecture aurait des conséquences remarquables sur la structure tannakienne des représentations automorphes.

(ii) Ce théorème est utilisé pour étudier la fonctorialité, principalement l'implication « automorphe ⇒ modulaire » (cf. n° 3.3.4 ci-dessous).

### 3.3. Fonctorialité

Ce qui suit est fortement inspiré d'un article de survol de COGDELL [60].

*3.3.1. Dual de Langlands.* — Soit $\mathbb{G}$ un groupe algébrique réductif [56], défini sur $F$. On note $^L\mathbb{G}$ le dual de Langlands de $\mathbb{G}$; c'est un produit semi-direct de $^L\mathbb{G}_0$ (composante connexe de l'identité de $^L\mathbb{G}$, qui est un $\mathbf{C}$-groupe algébrique) par $\mathrm{Gal}_F$.

• Si $\mathbb{G} = \mathbf{GL}_n$, alors $^L\mathbb{G} = \mathbf{GL}_n(\mathbf{C}) \times \mathrm{Gal}_F$.
• Si $\mathbb{G} = \mathbf{SO}_{2n+1}, \mathbf{Sp}_{2n}, \mathbf{SO}_{2n}$, alors $^L\mathbb{G}_0 = \mathbf{Sp}_{2n}(\mathbf{C}), \mathbf{SO}_{2n+1}(\mathbf{C}), \mathbf{SO}_{2n}(\mathbf{C})$.
• Si $\mathbb{G}$ et $\mathbb{G}'$ sont des formes intérieures [57] l'un de l'autre, alors $^L\mathbb{G} = {^L\mathbb{G}'}$.

---

56. Sans sous-groupe unipotent distingué non trivial.

57. Cela signifie qu'il existe un isomorphisme de groupes algébriques $\phi : \mathbb{G} \xrightarrow{\sim} \mathbb{G}'$, défini sur $\overline{F}$, tel que, pour tout $\sigma \in \mathrm{Gal}_F$, l'automorphisme $\phi^{-1} \circ \phi^\sigma$ de $\mathbb{G}$ soit la conjugaison par un élément de $\mathbb{G}(\overline{F})$.

*3.3.2. Fonctions L automorphes*. — On définit les formes automorphes pour $\mathbb{G}(\mathbf{A}_F)$ et la notion de cuspidalité de la même manière exactement que pour $\mathbb{G} = \mathbf{GL}_n$. On note $\mathscr{A}(\mathbb{G})$ (ou $\mathscr{A}(\mathbb{G}(\mathbf{A}_F))$ si on veut expliciter le corps $F$) l'ensemble des représentations automorphes de $\mathbb{G}(\mathbf{A}_F)$, et $\mathscr{A}^0(\mathbb{G})$ (ou $\mathscr{A}^0(\mathbb{G}(\mathbf{A}_F))$) le sous-ensemble des représentations cuspidales.

Comme ci-dessus une représentation automorphe $\pi$ de $\mathbb{G}(\mathbf{A}_F)$ se factorise sous la forme $\otimes'_v \pi_v$, où $v$ est non ramifiée pour tout $v$ sauf un nombre fini. La différence est que l'isomorphisme de Satake ne produit pas une famille de nombres complexes $\alpha_{v,i}$ mais une classe de conjugaison [58] $\alpha(\pi_v)$ dans $^L\mathbb{G}$, semi-simple, se projetant sur la classe de Frob$_v$ via la flèche $^L\mathbb{G} \to \mathrm{Gal}_F$.

Pour fabriquer une fonction $L$, on a besoin d'un paramètre supplémentaire, à savoir une représentation continue $r\colon {}^L\mathbb{G} \to \mathbf{GL}_n(\mathbf{C})$. Alors, pour presque tout $v$, $r(\alpha(\pi_v))$ fournit une classe de conjugaison de $\mathbf{GL}_n(\mathbf{C})$, et donc un facteur d'Euler $L(\pi_v, r, s) = E_v(|v|^{-s})^{-1}$ où $E_v(T)$ est le déterminant de $1 - Tr(\alpha(\pi_v))$. D'où une fonction $[L(\pi, r, s)]$, classe de $\prod_v L(\pi_v, r, s)$. Une fonction $[L]$ de cette forme est une *fonction $[L]$ automorphe.*

*3.3.3. Principe de fonctorialité*. — Le principe de fonctorialité de Langlands affirme que, si $\mathbb{G}$ et $\mathbb{H}$ sont des $F$-groupes réductifs avec $\mathbb{G}$ quasi-déployé [59] et si $\rho\colon {}^L\mathbb{H} \to {}^L\mathbb{G}$ est un morphisme de $L$-groupes (on demande que la restriction à $^L\mathbb{H}_0$ soit un morphisme de groupes de Lie complexes et que la flèche $\mathrm{Gal}_F \to \mathrm{Gal}_F$, induite par $\rho$ par passage au quotient, soit l'identité), alors il existe une application naturelle $\pi \mapsto \Pi(\pi, \rho)$ de $\mathscr{A}(\mathbb{H})$ dans $\mathscr{A}(\mathbb{G})$ compatible aux fonctions $L$ au sens que, si $r\colon {}^L\mathbb{G} \to \mathbf{GL}_n(\mathbf{C})$ est un morphisme continu de groupes, alors $[L(\pi, r \circ \rho, s)] = [L(\Pi(\pi, \rho), r, s)]$.

*Remarque 3.9*. — (i) Si $\mathbb{H} = \{1\}$ et $\mathbb{G} = \mathbf{GL}_n$, les $\rho\colon {}^L\mathbb{H} \to {}^L\mathbb{G}$ correspondent aux représentations de dimension $n$ de $\mathrm{Gal}_F$. Comme les représentations automorphes de $\{1\}$ se réduisent à la représentation triviale, le principe de fonctorialité infère qu'il existe $\Pi(\rho) \in \mathscr{A}(\mathbb{G})$ telle que $L(\Pi(\rho), s) \sim L(\rho, s)$. Si on utilise les équations fonctionnelles de $L(\Pi(\rho), s)$ et $L(\rho, s)$, on en déduit que $L(\Pi(\rho), s) = L(\rho, s)$ et donc, en particulier, la conjecture d'Artin pour $\rho$.

(ii) Si $\mathbb{H} = \mathbf{GL}_n$, $\mathbb{G} = \mathbf{GL}_1$ et $\rho\colon {}^L\mathbb{H} \to {}^L\mathbb{G}$ est le déterminant sur $^L\mathbb{H}_0$, l'application $\mathscr{A}(\mathbb{H}) \to \mathscr{A}(\mathbb{G})$ prédite par le principe de fonctorialité est celle qui envoie une représentation automorphe sur son caractère central.

(iii) Si $\mathbb{G}$ n'est pas quasi-déployé, pour que $\Pi(\pi, \rho)$ existe, il y a des restrictions sur les $\pi_v$ pour les $v$ tels que $\mathbb{G}$ n'est pas quasi-déployé sur $F_v$. Typiquement, si $\mathbb{H}$ est le groupe algébrique associé au groupe multiplicatif d'une algèbre de quaternions sur $F$ et si $\mathbb{G} = \mathbf{GL}_2$, alors $\pi_v$ ne peut pas être une série principale si $\mathbb{H}(F_v)$ n'est pas isomorphe à $\mathbf{GL}_2(F_v)$.

---

58. Pour faire le lien avec le cas $\mathbb{G} = \mathbf{GL}_n$, notons que $n$ nombres complexes non nuls définissent une classe de conjugaison semi-simple de $\mathbf{GL}_n(\mathbf{C})$.

59. Un groupe algébrique sur un corps $k$ est déployé s'il contient un tore maximal défini sur $k$ qui est un produit de groupes multiplicatifs $\mathbf{G}_m$; il est quasi-déployé s'il contient un sous-groupe de Borel défini sur $k$. Si $G$ est un $k$-groupe algébrique réductif, il existe une unique forme intérieure de $G$ qui est quasi-déployée.

*3.3.4. Exemples.* — On est très loin d'avoir établi le principe de fonctorialité en toute généralité mais il y a quand même des avancées significatives très utiles. Toutes les fonctorialités ci-dessous sont des transferts vers $\mathbf{GL}_n$ ; elles reposent sur le th. 3.7 qui demande d'avoir préétabli les propriétés des fonctions $L$ impliquées par l'existence du transfert. L'étude de ces fonctions $L$ se fait via celle des coefficients de Fourier de séries d'Eisenstein (méthode de Langlands–Shahidi, voir l'article de survol de Shahidi [198]) ou bien via des variations sur la méthode de Rankin–Selberg (voir l'article de survol de Bump [38]).

• *Groupes classiques.* — On prend $\mathbb{H} = \mathbf{SO}_{2n+1}, \mathbf{Sp}_{2n}, \mathbf{SO}_{2n}$, de telle sorte que ${}^L\mathbb{H}_0 = \mathbf{Sp}_{2n}(\mathbf{C}), \mathbf{SO}_{2n+1}(\mathbf{C}), \mathbf{SO}_{2n}(\mathbf{C})$, et $\mathbb{G} = \mathbf{GL}_{2n}, \mathbf{GL}_{2n+1}, \mathbf{GL}_{2n}$, de telle sorte que ${}^L\mathbb{H}_0$ se plonge naturellement dans ${}^L\mathbb{G}_0$ et ${}^L\mathbb{H}$ dans ${}^L\mathbb{G}$. L'existence du transfert $\mathscr{A}(\mathbb{H}) \to \mathscr{A}(\mathbb{G})$ correspondant à ce plongement naturel a été établie par Cogdell, Kim, Piatetski-Shapiro, Shahidi ([61] pour $\mathbb{H} = \mathbf{SO}_{2n+1}$ et [62] pour les autres cas). Dans tous ces cas, les représentations automorphes de $\mathbb{G}$ que l'on obtient vérifient une condition d'autodualité.

• *Produits tensoriels.* — On prend $\mathbb{H} = \mathbf{GL}_n \times \mathbf{GL}_m$ de telle sorte que ${}^L\mathbb{H}_0 = \mathbf{GL}_n(\mathbf{C}) \times \mathbf{GL}_m(\mathbf{C})$ et $\mathbb{G} = \mathbf{GL}_{nm}$ de telle sorte que[60] $\otimes \colon \mathbf{GL}_n(\mathbf{C}) \times \mathbf{GL}_m(\mathbf{C}) \to \mathbf{GL}_{nm}(\mathbf{C})$ induise un plongement de ${}^L\mathbb{H}_0$ dans ${}^L\mathbb{G}_0$ et de ${}^L\mathbb{H}$ dans ${}^L\mathbb{G}$. L'existence du produit tensoriel automorphe correspondant à ce plongement a été établie pour $n = m = 2$ par Ramakrishnan [184] et pour $n = 2$ et $m = 3$ par Kim, Shahidi [148].

• *Puissances symétriques.* — On prend $\mathbb{H} = \mathbf{GL}_2$ de telle sorte que ${}^L\mathbb{H}_0 = \mathbf{GL}_2(\mathbf{C})$. Si $k \geq 1$, l'action de $\mathbf{GL}_2(\mathbf{C})$ sur la puissance symétrique $k$-ième $\mathrm{Sym}^k$ de la représentation standard de $\mathbf{GL}_2(\mathbf{C})$ sur $\mathbf{C}^2$ induit un plongement $\mathrm{Sym}^k \colon \mathbf{GL}_2(\mathbf{C}) \to \mathbf{GL}_{k+1}(\mathbf{C})$. Si $\mathbb{G} = \mathbf{GL}_{k+1}$, l'existence de $\mathrm{Sym}^k \colon \mathscr{A}(\mathbb{H}) \to \mathscr{A}(\mathbb{G})$ correspondant à ce plongement prédite par le principe de fonctorialité a été établie par Gelbart, Jacquet [120] pour $k = 2$, par Kim, Shahidi [148, 149] pour $k = 3$, et par Kim pour $k = 4$ [147].

Cela prouve que, si $\pi$ est une représentation automorphe de $\mathbf{GL}_2$ (par exemple une représentation attachée à une forme modulaire ou une forme modulaire de Hilbert), les $L(\pi, \mathrm{Sym}^k, s)$ ont un prolongement holomorphe et une équation fonctionnelle pour $k = 1, 2, 3, 4$. Mais l'automorphie de $\mathrm{Sym}^k \pi$ combinée avec le th. 3.6 et la décomposition de $\mathrm{Sym}^r \otimes \mathrm{Sym}^s$ permet de montrer l'existence de prolongements méromorphes et d'équations fonctionnelles des $L(\pi, \mathrm{Sym}^k, s)$ pour $k = 5, 6, 7, 8$ et même $k = 9$ en injectant d'autres ingrédients automorphes [149, prop. 4.7].

Dans [147], Kim traite aussi le carré extérieur $\wedge^2 \colon \mathscr{A}(\mathbf{GL}_4) \to \mathscr{A}(\mathbf{GL}_6)$ avec quelques restrictions levées par Henniart [136].

---

60. On fait agir $\mathbf{GL}_n(\mathbf{C}) \times \mathbf{GL}_m(\mathbf{C})$ de manière naturelle sur $\mathbf{C}^n$ (resp. $\mathbf{C}^m$), le second (resp. premier) facteur agissant trivialement ; l'action sur $\mathbf{C}^n \otimes \mathbf{C}^m$ qui s'en déduit fournit l'application $\otimes$ ci-dessus.

*3.3.5. Changement de corps*

• Changement de base. — On part d'une extension $K/F$ de corps de nombres et d'un groupe réductif $\mathbb{H}$ déployé sur $F$. Soit $\mathbb{G}$ la restriction des scalaires [61] de $K$ à $F$ de $\mathbb{H}$ : si $\Lambda$ est une $F$-algèbre, alors $\mathbb{G}(\Lambda) = \mathbb{H}(K \otimes_F \Lambda)$ ; en particulier $\mathbb{G}(\mathbf{A}_F) = \mathbb{H}(\mathbf{A}_K)$ et les représentations automorphes de $\mathbb{G}(\mathbf{A}_F)$ sont les mêmes que celles de $\mathbb{H}(\mathbf{A}_K)$. Une fois n'est pas coutume, le dual de Langlands de $\mathbb{G}$ est un vrai produit semi-direct et pas un produit : ${}^L\mathbb{G}_0 = \prod_{\mathrm{Gal}_K \backslash \mathrm{Gal}_F} {}^L\mathbb{H}_0$ et $\mathrm{Gal}_F$ agit par permutation des facteurs (par translation sur $\mathrm{Gal}_K \backslash \mathrm{Gal}_F$). On a ${}^L\mathbb{H} = {}^L\mathbb{H}_0 \times \mathrm{Gal}_F$, et on dispose d'une injection diagonale ${}^L\mathbb{H}_0 \hookrightarrow {}^L\mathbb{G}_0$, et cette injection s'étend en une injection ${}^L\mathbb{H} \hookrightarrow {}^L\mathbb{G}$. L'application $\mathscr{A}(\mathbb{H}(\mathbf{A}_F)) \to \mathscr{A}(\mathbb{H}(\mathbf{A}_K))$ correspondant à cette injection, dont l'existence est prédite par le principe de fonctorialité est *le changement de base de $F$ à $K$*, noté $\mathrm{BC}_F^K$.

L'existence du changement de base a été prouvée par ARTHUR et CLOZEL [5] dans le cas où $\mathbb{H} = \mathbf{GL}_n$ et l'extension $K/F$ est galoisienne de groupe de Galois résoluble par dévissage à partir du cas où l'extension est cyclique (le cas $n = 2$ avait été traité auparavant par Langlands [161]). De plus, si $K/F$ est cyclique, et si $\Pi \in \mathscr{A}^0(\mathbb{G})$ vérifie $\Pi \circ \sigma = \Pi$ pour tout $\sigma \in \mathrm{Gal}(F/K)$, alors il existe $\pi \in \mathscr{A}^0(\mathbb{H})$ telle que $\Pi = \mathrm{BC}_F^K(\pi)$ ; par contre, il n'existe pas de recette complètement explicite décrivant $\pi$ en termes de $\Pi$.

• Induction automorphe. — On part d'une extension $K/F$ de corps de nombres de degré $d$. Soit $\mathbb{H}$ la restriction des scalaires de $K$ à $F$ de $\mathbf{GL}_n$ et soit $\mathbb{G} = \mathbf{GL}_{nd}$ (vus tous deux comme des $F$-groupes). On dispose d'une injection naturelle de $L$-groupes ${}^L\mathbb{H} \hookrightarrow {}^L\mathbb{G}$ envoyant ${}^L\mathbb{H}_0 = \prod_{\mathrm{Gal}_K \backslash \mathrm{Gal}_F} \mathrm{GL}_n(\mathbf{C})$ dans ${}^L\mathbb{G} = \mathrm{GL}_{nd}(\mathbf{C})$ sur le sous-groupe des matrices diagonales par blocs $n \times n$ (cela demande de numéroter les éléments de $\mathrm{Gal}_K \backslash \mathrm{Gal}_F$), et $\mathrm{Gal}_F$ dans $S_d \times \mathrm{Gal}_F \subset {}^L\mathbb{G}$, où $S_d \subset \mathrm{GL}_{nd}(\mathbf{C})$ est le groupe des matrices de permutations par blocs $n \times n$. L'application $\mathscr{A}(\mathbb{H}) \to \mathscr{A}(\mathbb{G})$ correspondant à cette injection, dont l'existence est prédite par le principe de fonctorialité, est *l'induction automorphe de $K$ à $F$*.

L'existence de l'induction automorphe dans le cas où l'extension $K/F$ est galoisienne de groupe de Galois résoluble a été prouvée par ARTHUR et CLOZEL [5] (le cas $n = 2$ avait été traité auparavant par JACQUET et LANGLANDS [139]).

# 4. BOITE À OUTILS DE THÉORIE DE HODGE $p$-ADIQUE

La théorie de Hodge $p$-adique fournit certains des outils les plus puissants utilisés dans la preuve des théorèmes de relèvement modulaire. En retour la quête de ces derniers a puissamment poussé au développement de certains aspects de la théorie (en particulier, les questions d'intégralité et l'étude des familles de représentations $p$-adiques). Par exemple, les premiers résultats de modularité n'utilisaient que la théorie de Fontaine–Laffaille (limitée au cas non ramifié et aux petits poids) ; la preuve de des derniers cas de la conjecture de Taniyama–Weil [33] a eu besoin de la classification [62] de BREUIL [32] des schémas en groupe finis sur $\mathscr{O}_K$ pour $p > 3$ (mais sans restriction sur l'extension finie $K$ de $\mathbf{Q}_p$). La classification de BREUIL a été revisitée par KISIN [151] ce qui a conduit à la définition des modules de Breuil–Kisin avec, comme applications : la description des réseaux des représentations cristallines à poids 0 et 1 (cette restriction sur les poids a été levée récemment [21, 117, 85] grâce à l'introduction du point de vue prismatique), l'existence des variétés de Kisin [152] paramétrant les représentations potentiellement semi-stables de type fixé et représentation résiduelle fixée, et des théorèmes de modularité [152, 153, 154]. La géométrie des variétés de Kisin est étroitement liée à la conjecture de Breuil–Mézard, un des points de départ de la correspondance de Langlands locale $p$-adique pour $\mathbf{GL}_2(\mathbf{Q}_p)$. Le champ des $(\varphi, \Gamma)$-modules d'EMERTON–GEE [91] est un objet fascinant qui permet de voir les représentations résiduelles comme les points d'une variété algébrique et les variétés de Kisin correspondant à un type fixé comme les tubes de points de la fibre spéciale d'une sous-variété analytique de l'espace de toutes les représentations de $\mathrm{Gal}_K$ de dimension fixée ; ce champ offre une perspective complètement naturelle sur la conjecture de Breuil–Mézard ; il est aussi utilisé dans [27] pour analyser les représentations résiduelles récalcitrantes par « prolongement analytique ».

---

61. Cf. [26, I.4, I.5] pour la construction de $\mathbb{G}$ et de son dual de Langlands.

Dans une autre direction, les représentations trianguline ont été introduites [66] pour décrire les représentations de la série principale unitaire de $\mathbf{GL}_2(\mathbf{Q}_p)$, elles étaient inspirées par des travaux de KISIN [150] sur la conjecture de Fontaine–Mazur ; les familles analytiques de représentations triangulines ont des applications aux théorèmes de modularité [34, 35, 171].

Enfin, le versant géométrique de la théorie de Hodge $p$-adique a été le témoin d'une véritable révolution depuis les travaux de SCHOLZE [190]. En particulier, l'étude des tours complétées de variétés de Shimura par SCHOLZE (pour la construction de représentations galoisienne de torsion associées aux classes de cohomologie singulière de torsion de certains espaces symétriques) a fourni de nouveaux outils pour étudier les représentations obtenues comme limites de représentations apparaissant dans la cohomologie des variétés de Shimura (ces limites n'apparaissent pas en niveau fini, mais elles apparaissent dans la cohomologie de la tour complétée) ; voir le n° 2.1.3 pour des applications de ces idées.

### 4.1. Théorie de Hodge des variétés $p$-adiques

On note $\mathbf{C}_p$ le complété de $\overline{\mathbf{Q}}_p$ et $\mathbf{Q}_p^{\mathrm{nr}} \subset \overline{\mathbf{Q}}_p$ l'extension maximale non ramifiée de $\mathbf{Q}_p$ (c'est la réunion des $\mathbf{Q}_p(\boldsymbol{\mu}_N)$, pour $(N, p) = 1$).

Si $K \subset \overline{\mathbf{Q}}_p$ est une extension finie de $\mathbf{Q}_p$, on note $K_\infty$ l'extension cyclotomique $K(\boldsymbol{\mu}_{p^\infty})$, et $H_K \subset \mathrm{Gal}_K$ le groupe $\mathrm{Gal}(\overline{\mathbf{Q}}_p/K_\infty)$ et $\Gamma_K$ le quotient $\mathrm{Gal}_K/\Gamma_K$. Alors le caractère cyclotomique $\chi \colon \mathrm{Gal}_K \to \mathbf{Z}_p^*$ identifie $\Gamma_K$ à un sous-groupe ouvert de $\mathbf{Z}_p^*$.

---

62. Que cette classification ait été disponible à ce moment précis est une coïncidence heureuse.

*4.1.1. Construction des anneaux de périodes p-adiques* [103]

- *L'anneau* $\widetilde{\mathbf{A}}^+$. — Soit

$$\widetilde{\mathbf{E}}^+ = \mathscr{O}_{\mathbf{C}_p^\flat} := \{(x_n)_{n\in\mathbf{N}},\ x_n \in \mathscr{O}_{\mathbf{C}_p}/p,\ x_{n+1}^p = x_n,\ \forall n \in \mathbf{N}\} = \varprojlim\nolimits_{x\mapsto x^p} \mathscr{O}_{\mathbf{C}_p}/p$$

Alors $\widetilde{\mathbf{E}}^+$ est un anneau parfait de caractéristique $p$. Si $x = (x_n)_{n\in\mathbf{N}} \in \widetilde{\mathbf{E}}^+$, et si $\hat{x}_n$ est un relèvement de $x_n$ dans $\mathscr{O}_{\mathbf{C}_p}$, alors $\hat{x}_n^{p^n}$ converge dans $\mathscr{O}_{\mathbf{C}_p}$, et la limite $x^\sharp$ ne dépend pas du choix des $\hat{x}_n$. Si on pose $v^\flat(x) = v_p(x^\sharp)$, alors $v^\flat$ est une valuation sur $\widetilde{\mathbf{E}}^+$ pour laquelle il est complet. Il s'ensuit que, si $\alpha \in \widetilde{\mathbf{E}}^+$ vérifie $v^\flat(\alpha) > 0$, alors $\mathbf{C}_p^\flat = \widetilde{\mathbf{E}} := \widetilde{\mathbf{E}}^+[\frac{1}{\alpha}]$ est le corps des fractions de $\widetilde{\mathbf{E}}^+$ et que $v^\flat$ s'étend en une valuation de $\widetilde{\mathbf{E}}$ pour laquelle il est complet; de plus $\widetilde{\mathbf{E}}$ est algébriquement clos.

On fait agir $\mathrm{Gal}_{\mathbf{Q}_p}$ sur $\widetilde{\mathbf{E}}^+$ composante par composante (via son action sur $\mathscr{O}_{\mathbf{C}_p}$); cette action s'étend naturellement à $\widetilde{\mathbf{E}}$, et cette action est continue.

Soit [63] $\varepsilon_n = e^{2i\pi/p^n}$, de telle sorte que $\varepsilon = (\varepsilon_n)_{n\in\mathbf{N}}$ est un élément de $\widetilde{\mathbf{E}}^+$ sur lequel $\sigma \in \mathrm{Gal}_{\mathbf{Q}_p}$ agit par $\sigma(\varepsilon) = \varepsilon^{\chi(\sigma)}$ (car $\sigma(\varepsilon_n) = \varepsilon_n^{\chi(\sigma)}$, par définition de $\chi$).

Soit [64] $\widetilde{\mathbf{A}}^+ = W(\widetilde{\mathbf{E}}^+)$ l'anneau des vecteurs de Witt à coefficients dans $\widetilde{\mathbf{E}}^+$. Il existe un unique système multiplicatif de représentants de $\widetilde{\mathbf{E}}^+$ dans $\widetilde{\mathbf{A}}^+$ (les *représentants de Teichmüller*). Si $x \in \widetilde{\mathbf{E}}^+$, notons $[x]$ son représentant de Teichmüller. Tout élément de $\mathbf{A}_{\mathrm{inf}}$ peut s'écrire, de manière unique, sous la forme $\sum_{k\in\mathbf{N}} p^k[x_k]$, où les $x_k$ sont des éléments arbitraires de $\widetilde{\mathbf{E}}^+$. Par fonctorialité des vecteurs de Witt, $\widetilde{\mathbf{A}}^+$ est muni d'un frobenius $\varphi$ donné par $\varphi(\sum_{k\in\mathbf{N}} p^k[x_k]) = \sum_{k\in\mathbf{N}} p^k[x_k^p]$, et d'une action de $\mathrm{Gal}_{\mathbf{Q}_p}$ commutant à $\varphi$.

On définit $\theta\colon \widetilde{\mathbf{A}}^+ \to \mathscr{O}_{\mathbf{C}_p}$ par $\theta(\sum_{k\in\mathbf{N}} p^k[x_k]) = \sum_{k\in\mathbf{N}} p^k x_k^\sharp$. Alors $\theta$ est un morphisme surjectif d'anneaux dont le noyau est engendré par $\xi = \frac{[\varepsilon]-1}{[\varepsilon^{1/p}]-1}$.

- *Les anneaux* $\mathbf{B}_{\mathrm{cris}}$, $\mathbf{B}_{\mathrm{st}}$, $\mathbf{B}_{\mathrm{dR}}$. — Soit $\mathbf{B}_{\mathrm{dR}}^+ = \varprojlim(\mathbf{A}_{\mathrm{inf}}[\frac{1}{p}]/\xi^k)$. C'est un anneau de valuation discrète complet, de corps résiduel [65] $\mathbf{C}_p$, contenant le complété $\mathbf{A}_{\mathrm{cris}}$ de $\widetilde{\mathbf{A}}^+[\frac{\xi^k}{k!},\ k \in \mathbf{N}]$ pour la topologie $p$-adique. L'action de $\mathrm{Gal}_{\mathbf{Q}_p}$ s'étend à tous ces anneaux et, si on pose

$$t = \log[\varepsilon] = -\sum_{k\geq 1} \tfrac{(1-[\varepsilon])^k}{k},$$

alors $t \in \mathbf{A}_{\mathrm{cris}}$ est une uniformisante de $\mathbf{B}_{\mathrm{dR}}^+$, et

$$\sigma(t) = \log[\varepsilon^{\chi(\sigma)}] = \log\left([\varepsilon]^{\chi(\sigma)}\right) = \chi(\sigma)t,$$

ce qui fait de $t$ un analogue $p$-adique de $2i\pi$.

Le frobenius $\varphi$ s'étend par continuité à $\mathbf{A}_{\mathrm{cris}}$, et $\varphi(t) = pt$. L'action de $\varphi$ s'étend donc au sous-anneau $\mathbf{B}_{\mathrm{cris}} = \mathbf{A}_{\mathrm{cris}}[\frac{1}{t}] = \mathbf{B}_{\mathrm{cris}}^+[\frac{1}{t}]$ (avec $\mathbf{B}_{\mathrm{cris}}^+ = \mathbf{A}_{\mathrm{cris}}[\frac{1}{p}]$) de $\mathbf{B}_{\mathrm{dR}} = \mathbf{B}_{\mathrm{dR}}^+[\frac{1}{t}]$. On munit $\mathbf{B}_{\mathrm{dR}}$ de la filtration décroissante par les $\mathbf{B}_{\mathrm{dR}}^i = t^i \mathbf{B}_{\mathrm{dR}}^+$, pour $i \in \mathbf{Z}$ ; cette filtration est stable par $\mathrm{Gal}_{\mathbf{Q}_p}$.

---

63. Rappelons que l'on a fixé des plongements de $\overline{\mathbf{Q}}$ dans $\mathbf{C}$ et $\overline{\mathbf{Q}}_p$.

64. L'anneau $\widetilde{\mathbf{A}}^+$ est souvent noté $\mathbf{A}_{\mathrm{inf}}$.

65. Et donc, si on croit à l'axiome du choix, isomorphe à $\mathbf{B}_{\mathrm{dR}}^+ \cong \mathbf{C}_p[[t]]$, mais il y a au moins deux raisons pour lesquelles c'est une très mauvaise idée de croire que $\mathbf{B}_{\mathrm{dR}}^+ = \mathbf{C}_p[[t]]$ : une telle identification ne peut être ni continue, ni compatible à l'action de Galois.

Soit $u \in \mathbf{B}_{\mathrm{dR}}$ l'analogue $p$-adique de $\log p$ défini par
$$u = \log \tfrac{[p^\flat]}{p} = \sum_{k\geq 1} \tfrac{(-1)^{k-1}}{k}\bigl(\tfrac{[p^\flat]}{p}-1\bigr)^k,$$
où $p^\flat = (p, p^{1/p}, p^{1/p^2}, \dots) \in \mathscr{O}_{\mathbf{C}_p^\flat}$. Alors $u$ est transcendant sur $\mathbf{B}_{\mathrm{cris}}$ et le sous-anneau $\mathbf{B}_{\mathrm{st}} := \mathbf{B}_{\mathrm{cris}}[u]$ de $\mathbf{B}_{\mathrm{dR}}$ est stable par $\mathrm{Gal}_{\mathbf{Q}_p}$ : il existe $c\colon \mathrm{Gal}_{\mathbf{Q}_p} \to \mathbf{Z}_p$ tel que $\sigma(u) = u + c(\sigma)t$, si $\sigma \in \mathrm{Gal}_{\mathbf{Q}_p}$. On munit $\mathbf{B}_{\mathrm{st}}$ d'un frobenius en posant $\varphi(u) = pu$ et d'un « opérateur de monodromie » $N = \tfrac{-d}{du}$. On a la relation $N\varphi = p\varphi N$.

Si $[K : \mathbf{Q}_p] < \infty$ et $K_0 = K \cap \mathbf{Q}_p^{\mathrm{nr}}$, les points fixes sous l'action de $\mathrm{Gal}_K$ sont :
$$\mathbf{B}_{\mathrm{dR}}^{\mathrm{Gal}_K} = K \quad \text{et} \quad \mathbf{B}_{\mathrm{st}}^{\mathrm{Gal}_K} = \mathbf{B}_{\mathrm{cris}}^{\mathrm{Gal}_K} = K_0.$$

Ces anneaux sont reliés par *les suites exactes fondamentales*
$$0 \to \mathbf{B}_{\mathrm{cris}} \to \mathbf{B}_{\mathrm{st}} \xrightarrow{N} \mathbf{B}_{\mathrm{st}} \to 0,$$
$$0 \to \mathbf{B}_{\mathrm{cris}}^{\varphi=1} \to \mathbf{B}_{\mathrm{cris}} \xrightarrow{\varphi-1} \mathbf{B}_{\mathrm{cris}} \to 0,$$
$$0 \to \mathbf{Q}_p \to \mathbf{B}_{\mathrm{cris}}^{\varphi=1} \longrightarrow \mathbf{B}_{\mathrm{dR}}/\mathbf{B}_{\mathrm{dR}}^0 \to 0.$$

Il en ressort que *l'on peut retrouver $\mathbf{Q}_p$ à l'intérieur de $\mathbf{B}_{\mathrm{st}}$* (ou $\mathbf{B}_{\mathrm{cris}}$) en utilisant les structures additionnelles ($\varphi$, $N$ et la filtration).

• L'anneau $\mathbf{B}_{\mathrm{HT}}$. — On note $\mathbf{B}_{\mathrm{HT}}$ le gradué de $\mathbf{B}_{\mathrm{dR}}$ ; on a donc $\mathbf{B}_{\mathrm{HT}} = \mathbf{C}_p[t, t^{-1}]$, et $\mathrm{Gal}_{\mathbf{Q}_p}$ agit sur $t$ par le caractère cyclotomique (i.e. $\sigma(t) = \chi(\sigma)t$ si $\sigma \in \mathrm{Gal}_{\mathbf{Q}_p}$) et la graduation $\mathrm{Gr}^i \mathbf{B}_{\mathrm{HT}} = t^i \mathbf{C}_p$ est stable par $\mathrm{Gal}_{\mathbf{Q}_p}$. Les points fixes de $\mathbf{B}_{\mathrm{HT}}$ sous l'action de $\mathrm{Gal}_K$ sont $\mathbf{B}_{\mathrm{HT}}^{\mathrm{Gal}_K} = K$.

*4.1.2. La conjecture $C_{\mathrm{st}}$.* — Si $X$ est une variété algébrique définie sur $K$, on dispose d'une cohomologie $H_{\mathrm{HK}}^\bullet$, la *cohomologie de Hyodo–Kato*, construite par HYODO, KATO [138] dans le cas semi-stable et par BEILINSON [10] dans le cas général (la construction de BEILINSON marche aussi pour les variétés singulières).

Les $H_{\mathrm{HK}}^i(X)$ sont des $\mathbf{Q}_p^{\mathrm{nr}}$-espaces vectoriels de dimension finie, munis d'un frobenius $\varphi$ bijectif, $\mathbf{Q}_p^{\mathrm{nr}}$-semi-linéaire, d'une action semi-linéaire localement constante[66] de $\mathrm{Gal}_K$ commutant à $\varphi$ et d'un opérateur $N$ « de monodromie » $\mathbf{Q}_p^{\mathrm{nr}}$-linéaire, nilpotent, vérifiant la relation $N\varphi = p\,\varphi N$ et commutant à $\mathrm{Gal}_K$. (Si $X$ est propre et a bonne réduction, alors $N = 0$ et l'action de $\mathrm{Gal}_K$ est non ramifiée.) De plus, on a un isomorphisme naturel « de Hyodo–Kato » de $K$-espaces vectoriels :
$$\iota_{\mathrm{HK}}\colon (\overline{\mathbf{Q}}_p \otimes_{\mathbf{Q}_p^{\mathrm{nr}}} H_{\mathrm{HK}}^i(X))^{\mathrm{Gal}_K} \cong H_{\mathrm{dR}}^i(X).$$

---

66. Cette action provient du fait qu'il faut étendre les scalaires pour trouver des modèles semi-stables, localement pour la topologie étale (l'existence d'un modèle semi-stable global, après extension des scalaires, n'est pas connue en général).

L'énoncé suivant a été conjecturé par FONTAINE (dans un cadre moins général), et a suscité énormément de travaux (par exemple [9, 10, 19, 97, 108, 174, 190, 211], liste non exhaustive, loin s'en faut) ; sous sa forme la plus générale, il est dû à BEILINSON [9, 10].

THÉORÈME 4.1. — *Soit $K$ une extension finie de $\mathbf{Q}_p$, et soit $X$ une variété projective lisse définie sur $K$. On a des isomorphismes naturels :*

$$\iota_{\mathrm{st}} \colon \mathbf{B}_{\mathrm{st}} \otimes_{\mathbf{Q}_p} H^i_{\mathrm{\acute{e}t}}(X_{\overline{K}}, \mathbf{Q}_p) \cong \mathbf{B}_{\mathrm{st}} \otimes_{\mathbf{Q}_p^{\mathrm{nr}}} H^i_{\mathrm{HK}}(X),$$
$$\iota_{\mathrm{dR}} \colon \mathbf{B}_{\mathrm{dR}} \otimes_{\mathbf{Q}_p} H^i_{\mathrm{\acute{e}t}}(X_{\overline{K}}, \mathbf{Q}_p) \cong \mathbf{B}_{\mathrm{dR}} \otimes_K H^i_{\mathrm{dR}}(X),$$

*commutant aux actions de $\mathrm{Gal}_K$, $\varphi$, $N$ et respectant les filtrations.*

*De plus, $\iota_{\mathrm{dR}}$ s'obtient à partir de $\iota_{\mathrm{st}}$ par extension des scalaires de $\mathbf{B}_{\mathrm{st}}$ à $\mathbf{B}_{\mathrm{dR}}$.*

### 4.2. La hiérarchie des représentations de $\mathrm{Gal}_K$

FONTAINE a développé, à partir de la fin des années 1970, un programme visant à classifier et décrire les $\mathbf{Q}_p$-*représentations de* $\mathrm{Gal}_K$ (i.e. les $\mathbf{Q}_p$-espaces vectoriels de dimension finie, munis d'une action $\mathbf{Q}_p$-linéaire continue de $\mathrm{Gal}_K$), en termes plus concrets (cf. [102, 104]).

*4.2.1. Représentations B-admissibles.* — La stratégie de FONTAINE part de l'observation suivante : si on dispose d'une $\mathbf{Q}_p$-algèbre topologique $B$, munie d'une action $\mathbf{Q}_p$-linéaire continue de $\mathrm{Gal}_K$ et de structures additionnelles stables sous l'action de $\mathrm{Gal}_K$, on peut associer à toute $\mathbf{Q}_p$-représentation $V$ de $\mathrm{Gal}_K$ un invariant $D_B(V)$ en prenant les points fixes $(B \otimes_{\mathbf{Q}_p} V)^{\mathrm{Gal}_K}$ de $B \otimes_{\mathbf{Q}_p} V$ sous l'action de $\mathrm{Gal}_K$. Alors $D_B(V)$ est un $B^{\mathrm{Gal}_K}$-module muni des structures additionnelles sur $B$ et qui est souvent plus facile à décrire que la représentation $V$ dont on est parti. Un tel anneau $B$ permet en outre de découper la sous-catégorie des *représentations B-admissibles* :

*Définition 4.2.* — Une $\mathbf{Q}_p$-représentation $V$ de $\mathrm{Gal}_K$ est *B-admissible* si $B \otimes_{\mathbf{Q}_p} V$ est triviale, i.e. isomorphe à $B^{\dim V}$ en tant que représentation de $\mathrm{Gal}_K$.

*4.2.2. La hiérarchie des représentations p-adiques.* — On peut appliquer le programme ci-dessus avec les anneaux $\mathbf{B}_{\mathrm{cris}}, \mathbf{B}_{\mathrm{st}}, \mathbf{B}_{\mathrm{dR}}, \mathbf{B}_{\mathrm{HT}}, \overline{\mathbf{Q}}_p \cdot \mathbf{B}_{\mathrm{cris}}$ (sous-anneau de $\mathbf{B}_{\mathrm{dR}}$ engendré par $\overline{\mathbf{Q}}_p$ et $\mathbf{B}_{\mathrm{cris}}$), $\overline{\mathbf{Q}}_p \cdot \mathbf{B}_{\mathrm{st}}$... Ceci permet de définir les notions suivantes pour une $\mathbf{Q}_p$-représentation $V$ de $\mathrm{Gal}_K$ :

• $V$ est dite *cristalline* si elle est $\mathbf{B}_{\mathrm{cris}}$-admissible et *potentiellement cristalline* si elle est $\overline{\mathbf{Q}}_p \cdot \mathbf{B}_{\mathrm{cris}}$-admissible.

• $V$ est dite *semi-stable* si elle est $\mathbf{B}_{\mathrm{st}}$-admissible et *potentiellement semi-stable* si elle est $\overline{\mathbf{Q}}_p \cdot \mathbf{B}_{\mathrm{st}}$-admissible.

• $V$ est dite *de Rham* si elle est $\mathbf{B}_{\mathrm{dR}}$-admissible.

• $V$ est dite *Hodge–Tate* si elle est $\mathbf{B}_{\mathrm{HT}}$-admissible.

Les relations entre les différents anneaux fournissent les implications suivantes :

$$
\begin{array}{ccc}
\text{cristalline} & \Longrightarrow & \text{potentiellement cristalline} \\
\Downarrow & & \Downarrow \\
\text{semi-stable} & \Rightarrow & \text{potentiellement semi-stable} \\
& & \Downarrow \\
& \text{de Rham} & \Longrightarrow \quad \text{Hodge–Tate}
\end{array}
$$

*Remarque 4.3.* — Toutes les implications ci-dessus sont strictes sauf « pst ⇒ dR » qui est, en fait, une équivalence. Nous renvoyons à [65, 68] pour une discussion des travaux consacrés à cette équivalence (connue sous le nom de « conjecture de monodromie $p$-adique de Fontaine ») et au th. 4.5 ci-dessous.

*4.2.3. Représentations potentiellement semi-stables et $(\varphi, N)$-modules filtrés.* — Soit $K$ une extension finie de $\mathbf{Q}_p$. Un $(\mathrm{Gal}_K, \varphi, N)$-*module filtré sur* $K$ est la donnée de :

- *un $(\mathrm{Gal}_K, \varphi, N)$-module $D$ sur $\mathbf{Q}_p^{\mathrm{nr}}$*, i.e. un $\mathbf{Q}_p^{\mathrm{nr}}$-espace vectoriel $D$ de dimension finie, muni :
  - ⋄ d'une action semi-linéaire de $\mathrm{Gal}_K$, localement constante,
  - ⋄ d'une action semi-linéaire bijective d'un frobenius $\varphi$ commutant à celle de $\mathrm{Gal}_K$,
  - ⋄ d'un opérateur $N$, commutant à $\mathrm{Gal}_K$, et vérifiant $N\varphi = p\varphi N$,
- une structure de $K$-*module filtré* sur $D_K = (\overline{\mathbf{Q}}_p \otimes_{\mathbf{Q}_p^{\mathrm{nr}}} D)^{\mathrm{Gal}_K}$, i.e. une filtration décroissante sur $D_K$ par des sous-$K$-espaces vectoriels $D_K^i$, pour $i \in \mathbf{Z}$, avec $D_K^i = D_K$ si $i$ est suffisamment petit, et $D_K^i = 0$ si $i$ est suffisamment grand.

Si $V$ est une représentation potentiellement semi-stable de $\mathrm{Gal}_K$, de dimension $d$, on pose

$$\mathbf{D}_{\mathrm{pst}}(V) := \varinjlim_{[L:K]<\infty}(\mathbf{B}_{\mathrm{st}} \otimes_{\mathbf{Q}_p} V)^{\mathrm{Gal}_L} \quad \text{et} \quad \mathbf{D}_{\mathrm{dR}}(V) := (\mathbf{B}_{\mathrm{dR}} \otimes_{\mathbf{Q}_p} V)^{\mathrm{Gal}_K}.$$

Alors $\mathbf{D}_{\mathrm{pst}}(V)$ est, naturellement, un $(\mathrm{Gal}_K, \varphi, N)$-module filtré sur $K$, de rang $d$ : $\mathbf{D}_{\mathrm{pst}}(V)$ est un $\mathbf{Q}_p^{\mathrm{nr}}$-module de rang $d$ puisque $\varinjlim_{[L:K]<\infty} \mathbf{B}_{\mathrm{st}}^{\mathrm{Gal}_L} = \mathbf{Q}_p^{\mathrm{nr}}$, est muni des actions de $\varphi$ et $N$ existant sur $\mathbf{B}_{\mathrm{st}}$ et l'action de $\mathrm{Gal}_K$ est lisse puisque, par construction, tout élément est fixe par $\mathrm{Gal}_L$ pour $L$ assez grand ; $\mathbf{D}_{\mathrm{dR}}(V)$ est un $K$-module de rang $d$ puisque $\mathbf{B}_{\mathrm{dR}}^{\mathrm{Gal}_K} = K$, et est muni de la filtration de $\mathbf{B}_{\mathrm{dR}}$ ; l'inclusion $(\overline{\mathbf{Q}}_p \otimes_{\mathbf{Q}_p^{\mathrm{nr}}} \mathbf{D}_{\mathrm{pst}}(V))^{\mathrm{Gal}_K} \hookrightarrow \mathbf{D}_{\mathrm{dR}}(V)$, induite par l'inclusion $\overline{\mathbf{Q}}_p \otimes_{\mathbf{Q}_p^{\mathrm{nr}}} \mathbf{B}_{\mathrm{st}} \hookrightarrow \mathbf{B}_{\mathrm{dR}}$, est une bijection pour des raisons de dimension.

Enfin, le gradué de $\mathbf{D}_{\mathrm{dR}}(V)$ est $\mathbf{D}_{\mathrm{HT}}(V) := (\mathbf{B}_{\mathrm{HT}} \otimes_{\mathbf{Q}_p} V)^{\mathrm{Gal}_K}$.

*Remarque 4.4.* — Si $V$ est cristalline, on a $\mathbf{D}_{\mathrm{pst}}(V) = \mathbf{Q}_p^{\mathrm{nr}} \otimes_{K_0} \mathbf{D}_{\mathrm{cris}}(V)$, où $\mathbf{D}_{\mathrm{cris}}(V) = (\mathbf{B}_{\mathrm{cris}} \otimes_{\mathbf{Q}_p} V)^{\mathrm{Gal}_K}$ est un $\varphi$-module filtré (i.e. $N = 0$ et l'action de $\mathrm{Gal}_K$ est triviale).

Si $D$ est un $(\mathrm{Gal}_K, \varphi, N)$-module filtré sur $K$, *le rang* $\mathrm{rg}(D)$ *de* $D$ est la dimension de $D_K$ sur $K$. Si $D$ est de rang 1, on définit le *degré* $\deg(D)$ *de* $D$ par la formule

$$\deg(D) = t_N(D) - t_H(D),$$

où $t_N(D)$ et $t_H(D)$ sont définis en choisissant une base $e$ de $D$ sur $\mathbf{Q}_p^{\mathrm{nr}}$ :

- il existe $\lambda \in (\mathbf{Q}_p^{\mathrm{nr}})^*$ tel que $\varphi(e) = \lambda e$, et on pose $t_N(D) = v_p(\lambda)$ ;
- il existe $i \in \mathbf{Z}$, unique, tel que $e \in (\overline{\mathbf{Q}}_p \otimes_K D_K^i) \setminus (\overline{\mathbf{Q}}_p \otimes_K D_K^{i+1})$, et on pose $t_H(D) = i$.

Si $D$ est de rang $r \geq 2$, alors $\det D = \wedge^r D$ est de rang 1, et on définit *le degré de $D$* par $\deg(D) = \deg(\det D)$ et *la pente de $D$* par $\mu(D) = \frac{\deg(D)}{\mathrm{rg}(D)}$. On dit que $D$ est *faiblement admissible*, si $\mu(D) = 0$, et si $\mu(D') \leq 0$, pour tout sous-objet $D'$ de $D$.

THÉORÈME 4.5. — $V \mapsto \mathbf{D}_{\mathrm{pst}}(V)$ *induit une équivalence de catégories de la catégorie des représentations potentiellement semi-stables de* $\mathrm{Gal}_K$ *sur celle des* $(\mathrm{Gal}_K, \varphi, N)$-*modules filtrés sur $K$ faiblement admissibles, le foncteur inverse étant*

$$D \mapsto \mathbf{V}_{\mathrm{st}}(D) = (\mathbf{B}_{\mathrm{st}} \otimes_{\mathbf{Q}_p^{\mathrm{nr}}} D)^{\varphi=1, N=0} \cap \mathrm{Fil}^0(\mathbf{B}_{\mathrm{dR}} \otimes_K D_K).$$

*Remarque 4.6.* — (i) Ce résultat [71] donne une description concrète des représentations potentiellement semi-stables mais, pour les applications aux théorèmes de relèvement modulaire, on a besoin de résultats plus fins, décrivant les réseaux de ces représentations (l'idéal serait d'avoir une description de ces réseaux modulo $p^n$) ; voir § 4.4 pour des résultats de ce type.

(ii) Si on part de $L$-représentations au lieu de $\mathbf{Q}_p$-représentations, on a le même résultat en rajoutant une action de $L$ commutant à celles de $\mathrm{Gal}_K$, $\varphi$ et $N$, et respectant la filtration sur $D_K$ (on fera attention que le $L \otimes_{\mathbf{Q}_p} K$-module $D_K^i$ n'a aucune raison d'être libre a priori).

(iii) On peut transformer un $(\mathrm{Gal}_K, \varphi, N)$-module en une représentation du groupe de Weil–Deligne $\mathrm{WD}_K$ de $K$ grâce à une recette de FONTAINE [105] : si $w \in \mathrm{W}_K$, on fait agir $w$ par $w \circ \varphi^{-\deg w}$ (ce qui rend l'action $\mathbf{Q}_p^{\mathrm{nr}}$-linéaire car les semi-linéarités de $w$ et $\varphi^{-\deg w}$ se compensent par définition de $\deg w$), et on a alors $N \circ w = p^{-\deg w} w \circ N$, ce qui est la relation requise pour une représentation de $\mathrm{WD}_K$. Cette construction peut s'inverser [37, prop. 4.1].

Notons que la recette de Fontaine marche tout aussi bien pour une représentation $\ell$-adique de $\mathrm{Gal}_K$, avec $\ell \neq p$, la différence est que, si $\ell \neq p$, toute représentation est potentiellement semi-stable.

(iv) On peut combiner la recette précédente avec la correspondance de Langlands locale pour associer une représentation lisse $\mathrm{LL}(V)$ de $\mathbf{GL}_n(K)$ à toute représentation $\ell$-adique $V$ (avec $\ell = p$ ou $\ell \neq p$) de $\mathrm{Gal}_K$, potentiellement semi-stable, de dimension $n$.

*4.2.4. Poids de Hodge–Tate.* — Soit $V$ une $\mathbf{Q}_p$-représentation de Rham de $\mathrm{Gal}_K$, de dimension $d$. L'ensemble des poids de Hodge–Tate de $V$ est le multi-ensemble $\mathrm{HT}(V) = \sum_i m_i\{i\}$, où $m_i = \dim_K \mathbf{D}_{\mathrm{HT}}^{-i}(V) = \dim_K(t^{-i}\mathbf{C}_p \otimes_{\mathbf{Q}_p} V)^{\mathrm{Gal}_K}$ (on a aussi $\mathbf{D}_{\mathrm{HT}}^{-i}(V) = \mathrm{gr}^{-i}\mathbf{D}_{\mathrm{dR}}(V)$). Le cardinal $\sum_i m_i$ de $\mathrm{HT}(V)$ est $d$ puisque $V$ est supposée de Rham.

Maintenant, si $V$ est une $L$-représentation de de Rham de $\mathrm{Gal}_K$, de dimension $n$, les $\mathbf{D}_{\mathrm{HT}}^{-i}(V)$ sont des $L \otimes K$-modules, pas nécessairement libres. Si $L$ contient tous les conjugués de $K$, chaque $\mathbf{D}_{\mathrm{HT}}^{-i}(V)$ est de la forme $\oplus_{\tau \in \mathrm{Hom}(K,L)} L^{m_{\tau,i}}$, où $a \in K$ agit sur $L^{m_{\tau,i}}$ par $\tau(a)$. Si $\tau \in \mathrm{Hom}(K,L)$, on définit l'ensemble $\mathrm{HT}_\tau(V)$ des *$\tau$-poids de Hodge–Tate de $V$* comme le multi-ensemble $\mathrm{HT}(V) = \sum_i m_{\tau,i}\{i\}$ (il est de cardinal $d$ pour tout $\tau$). On définit alors l'ensemble des poids de Hodge–Tate de $V$ comme la famille $(\mathrm{HT}_\tau(V))_\tau$ (où $\tau$ parcourt $\mathrm{Hom}(K,L)$).

*Remarque 4.7*. — Si $F$ est un corps de nombres, $\mathrm{Hom}(F,\overline{\mathbf{Q}}) = \sqcup_{v|p}\mathrm{Hom}(F_v,\overline{\mathbf{Q}}_p)$ (si $\tau : F \to \overline{\mathbf{Q}}$, alors $|\iota_p \circ \tau|_p$ définit une place $v$ de $F$ au-dessus de $p$, et $\iota_p \circ \tau$ s'étend par continuité à $F_v$). Si $\rho \colon \mathrm{Gal}_F \to \mathbf{GL}_n(\overline{\mathbf{Q}}_p)$ est continue, on peut trouver une extension finie $L$ de $\mathbf{Q}_p$ telle que $\rho$ se factorise par $\mathbf{GL}_n(L)$, et on peut supposer que $L$ contient tous les conjugués de $F$ (et donc tous ceux des $F_v$). Si $\tau \in \mathrm{Hom}(F,\overline{\mathbf{Q}})$, on définit l'ensemble $\mathrm{HT}_\tau(\rho)$ des $\tau$-*poids de Hodge–Tate* de $\rho$ comme celui de la restriction de $\rho$ à $\mathrm{Gal}_{F_v}$, où $v$ est l'unique place de $F$ avec $v \mid p$ telle que $\tau \in \mathrm{Hom}(F_v,\overline{\mathbf{Q}}_p)$.

*4.2.5. Opérateur de Sen*. — Si $V$ est une $\mathbf{Q}_p$-représentation de $\mathrm{Gal}_K$, de dimension $d$, alors $(\mathbf{C}_p \otimes_{\mathbf{Q}_p} V)^{H_K}$ est une représentation de Banach de $\Gamma_K$. On définit $\mathbf{D}_{\mathrm{Sen}}(V)$ comme le sous-espace des vecteurs localement analytiques pour l'action de $\Gamma_K$ (c'est un $K_\infty$-module de rang $d$), et $\Theta_{\mathrm{Sen}}$ comme $\lim_{\gamma \to 1} \frac{\gamma - 1}{\chi(\gamma) - 1}$ (générateur de l'algèbre de Lie de $\Gamma_K$) agissant sur $\mathbf{D}_{\mathrm{Sen}}(V)$ (c'est *l'opérateur de Sen* défini par Sen [192]; la description ci-dessus est issue de [16]). Les valeurs propres de $\Theta_{\mathrm{Sen}}$ sont les *poids de Hodge–Tate généralisés*[67] de $V$, et $V$ est Hodge–Tate si et seulement si $\Theta_{\mathrm{Sen}}$ est semi-simple et les poids de Hodge–Tate généralisés de $V$ sont des entiers[68].

## 4.3. La théorie des $(\varphi,\Gamma)$-modules

*4.3.1. Quelques anneaux gnomiques*. — Si $[K : \mathbf{Q}_p] < \infty$, on note $\mathbf{E}_K$ le corps des normes de $K_\infty$ vu comme sous-corps de $\widetilde{\mathbf{E}}$ : l'anneau $\mathbf{E}_K^+$ de ses entiers est l'ensemble des $(x_n)_n \in \widetilde{\mathbf{E}}^+$, avec $x_n \in \mathcal{O}_{K_n}$ modulo $p^{1/p}$, pour tout $n$ assez grand. En particulier, $\mathbf{E}_{\mathbf{Q}_p} = \mathbf{F}_p((\varepsilon - 1))$. La clôture séparable $\mathbf{E}$ de $\mathbf{E}_{\mathbf{Q}_p}$ est la réunion des $\mathbf{E}_K$ et on a $\mathrm{Gal}(\mathbf{E}/\mathbf{E}_K) = H_K$.

Soit $\widetilde{\mathbf{A}} = W(\widetilde{\mathbf{E}})$ l'anneau des vecteurs de Witt à coefficients dans $\widetilde{\mathbf{E}}$ (c'est aussi le complété de $\widetilde{\mathbf{A}}^+[\frac{1}{[p^\flat]}]$ pour la topologie $p$-adique). Soit $\pi = [\varepsilon] - 1 \in \widetilde{\mathbf{A}}$. On note $\mathbf{A}_{\mathbf{Q}_p}$ l'adhérence dans $\widetilde{\mathbf{A}}$ de $\mathbf{Z}_p[\pi, \pi^{-1}]$ ; c'est l'anneau des $\sum_{k \in \mathbf{Z}} a_k \pi^k$, avec $a_k \in \mathbf{Z}_p$ et $a_k \to 0$ quand $k \to -\infty$. On a

$$\varphi(\pi) = [\varepsilon]^p - 1 = (1+\pi)^p - 1 \quad \text{et} \quad \sigma(\pi) = [\varepsilon]^{\chi(\sigma)} - 1 = (1+\pi)^{\chi(\sigma)} - 1, \text{ si } \sigma \in \mathrm{Gal}_{\mathbf{Q}_p}.$$

Il s'ensuit que $\mathbf{A}_{\mathbf{Q}_p}$ est stable par $\varphi$ et par $\mathrm{Gal}_{\mathbf{Q}_p}$ qui agit à travers $\Gamma_{\mathbf{Q}_p}$.

Il existe un unique sous-anneau $\mathbf{A}$ de $\widetilde{\mathbf{A}}$, $p$-saturé ($x \in \widetilde{\mathbf{A}}$ et $px \in \mathbf{A} \Rightarrow x \in \mathbf{A}$) et complet pour la topologie $p$-adique, contenant $\mathbf{A}_{\mathbf{Q}_p}$ et tel que $\mathbf{A}/p\mathbf{A} = \mathbf{E} \subset \widetilde{\mathbf{E}} = \widetilde{\mathbf{A}}/p\widetilde{\mathbf{A}}$. Cet anneau est stable par $\mathrm{Gal}_{\mathbf{Q}_p}$ et par $\varphi$, et on a $\mathbf{A}^{H_{\mathbf{Q}_p}} = \mathbf{A}_{\mathbf{Q}_p}$ et $\mathbf{A}^{\varphi=1} = \mathbf{Z}_p$.

Si $[K : \mathbf{Q}_p] < \infty$, on pose $\mathbf{A}_K = \mathbf{A}^{H_K}$. Alors $\mathbf{A}_K$ est stable par $\varphi$ et par $\mathrm{Gal}_K$ agissant à travers $\Gamma_K$, et on a $\mathbf{A}_K/p\mathbf{A}_K = \mathbf{E}_K$.

---

67. Ou « poids de Hodge–Tate–Sen » ou simplement « poids de Hodge–Tate », voire juste « poids ».
68. Auquel cas, les poids de Hodge–Tate généralisés de $V$ ne sont rien d'autre que les poids de Hodge–Tate de $V$.

Enfin, on pose $\mathbf{B} = \mathbf{A}[\frac{1}{p}]$, $\mathbf{B}_K = \mathbf{A}_K[\frac{1}{p}]$. Alors $\mathbf{B}$ et $\mathbf{B}_K$ sont des corps munis d'actions de $\varphi$ et de $\mathrm{Gal}_K$ commutant entre elles, $\mathrm{Gal}_K$ agit à travers $\Gamma_K$ sur $\mathbf{B}_K$, et

$$\mathbf{B}^{\varphi=1} = \mathbf{Q}_p \quad \text{et} \quad \mathbf{B}^{H_K} = \mathbf{B}_K.$$

*Remarque 4.8.* — (i) En général, si $F = K_\infty \cap \mathbf{Q}_p^{\mathrm{nr}}$, il existe $\pi_K \in \mathbf{A}_K$ tel que $\mathbf{A}_K$ soit l'ensemble des $\sum_{k \in \mathbf{Z}} a_k \pi_K^k$, avec $a_k \in \mathscr{O}_F$ et $a_k \to 0$ quand $k \to -\infty$ (i.e. $\mathbf{A}_K$ a la même forme que $\mathbf{A}_{\mathbf{Q}_p}$), mais les formules donnant l'action de $\varphi$ et $\Gamma_K$ sur $\pi_K$ ne peuvent pas vraiment s'expliciter sauf si $K$ est non ramifiée (et alors $K = F$) où on peut prendre $\pi_K = \pi$ (l'expérience montre que ce n'est pas un problème).

(ii) Si $K$ est non ramifiée le sous-anneau $\mathbf{A}_K^+ = \mathscr{O}_K[[\pi]]$ de $\mathbf{A}_K$ est stable par $\varphi$ et $\Gamma_K$ et se surjecte sur $\mathbf{E}_K^+$; il n'existe pas de tel sous-anneau de $\mathbf{A}_K$ si $K$ n'est pas une extension abélienne de $\mathbf{Q}_p$.

*4.3.2. L'équivalence de catégories de Fontaine.* — Un $(\varphi, \Gamma)$-module $D$ sur $\Lambda = \mathbf{A}_K, \mathbf{B}_K$ est un $\Lambda$-module de type fini muni d'actions semi-linéaires de $\varphi$ et $\Gamma_K$ commutant entre elles. Un $(\varphi, \Gamma)$-module $D$ sur $\mathbf{A}_K$ (resp. $\mathbf{B}_K$) est *étale* si l'application naturelle[69] $\varphi^*D \to D$ est un isomorphisme (resp. si $D$ possède un $\mathbf{A}_K$-réseau, stable par $\varphi$ et $\Gamma_K$, qui est étale comme $(\varphi, \Gamma)$-module sur $\mathbf{A}_K$).

Si $V$ est une $\mathbf{Z}_p$ ou $\mathbf{Q}_p$-représentation de $\mathrm{Gal}_K$, on pose $D(V) := (\mathbf{A} \otimes_{\mathbf{Z}_p} V)^{H_K}$. C'est un $\mathbf{A}_K$ ou $\mathbf{B}_K$-module muni d'une action résiduelle de $\mathrm{Gal}_K/H_K = \Gamma_K$ et d'une action de $\varphi$ provenant du frobenius $\varphi$ sur $\mathbf{B}$; c'est donc un $(\varphi, \Gamma)$-module sur $\mathbf{A}_K$ ou $\mathbf{B}_K$.

THÉORÈME 4.9. — (FONTAINE [102]). *Si $V$ est une $\mathbf{Z}_p$ ou $\mathbf{Q}_p$-représentation de $\mathrm{Gal}_K$, le $(\varphi, \Gamma)$-module $D(V)$ est étale et $V \mapsto D(V)$ induit une équivalence de catégories de la catégorie des $\mathbf{Z}_p$ ou $\mathbf{Q}_p$-représentations de $\mathrm{Gal}_K$ sur celle des $(\varphi, \Gamma)$-modules étales sur $\mathbf{A}_K$ ou $\mathbf{B}_K$, le foncteur inverse étant $D \mapsto V(D) = (\mathbf{A} \otimes_{\mathbf{A}_K} D)^{\varphi=1}$.*

*4.3.3. Représentations triangulines.* — Soit $\mathbf{B}_{\mathbf{Q}_p}^\dagger \subset \mathbf{B}_{\mathbf{Q}_p}$ le sous-anneau des éléments surconvergents (i.e des $\sum_{k \in \mathbf{Z}} a_k \pi^k$ tels qu'il existe $r > 0$, tel que $\sum_{k \in \mathbf{Z}} a_k T^k$ définisse une fonction analytique bornée sur $0 < v_p(T) \leq r$).

Si $K$ est une extension finie de $\mathbf{Q}_p$, soit $\mathbf{B}_K^\dagger$ la clôture intégrale de $\mathbf{B}_{\mathbf{Q}_p}^\dagger$ dans $\mathbf{B}_K$. Alors $\mathbf{B}_K^\dagger$ est un corps et $[\mathbf{B}_K^\dagger : \mathbf{B}_{\mathbf{Q}_p}^\dagger] = [\mathbf{B}_K : \mathbf{B}_{\mathbf{Q}_p}]$. On dispose du résultat de descente suivant [52] :

THÉORÈME 4.10. — *Si $D$ est un $(\varphi, \Gamma)$-module sur $\mathbf{B}_K$, il existe un unique sous-$\mathbf{B}_K^\dagger$-module $D^\dagger$ de $D$, stable par $\varphi$ et $\Gamma_K$, tel que $D = \mathbf{B}_K \otimes_{\mathbf{B}_K^\dagger} D^\dagger$.*

Soit $\mathbf{B}_{\mathrm{rig},K}$ l'anneau de Robba : $\mathbf{B}_{\mathrm{rig},\mathbf{Q}_p}$ est l'ensemble des $\sum_{k \in \mathbf{Z}} a_k \pi^k$ tels qu'il existe $r > 0$ tel que $\sum_{k \in \mathbf{Z}} a_k T^k$ définisse une fonction analytique sur $0 < v_p(T) \leq r$, et $\mathbf{B}_{\mathrm{rig},K} = \mathbf{B}_K^\dagger \otimes_{\mathbf{B}_{\mathbf{Q}_p}^\dagger} \mathbf{B}_{\mathrm{rig},\mathbf{Q}_p}$.

Si $L$ est une extension finie de $\mathbf{Q}_p$, et si $V$ est une $L$-représentation de $\mathrm{Gal}_K$, le $(\varphi, \Gamma)$-module $D^\dagger(V)$ obtenu à partir de $D(V)$ par le th. 4.10 est un $L \otimes_{\mathbf{Q}_p} \mathbf{B}_K^\dagger$-module.

---

69. Si $M$ est un $\varphi$-module sur un anneau $\Lambda$, on pose $\varphi^*M := \Lambda \otimes_{\varphi(\Lambda)} \varphi(M)$.

Soit $\mathscr{R}_L := L \otimes_{\mathbf{Q}_p} \mathbf{B}_{\text{rig},K}$, et soit $D_{\text{rig}}(V) := \mathbf{B}_{\text{rig},K} \otimes_{\mathbf{B}_K^\dagger} D^\dagger(V)$ ; c'est un $\mathscr{R}_L$-module. On dit que $V$ est *trianguline* si $D_{\text{rig}}(V)$ est triangulable, i.e. est une extension successive de $(\varphi, \Gamma)$-modules de rang 1 sur $\mathscr{R}_L$.

*Remarque 4.11.* — (i) Si $K = \mathbf{Q}_p$, tout $(\varphi, \Gamma)$-module de rang 1 sur $\mathscr{R}_L$ est de la forme $\mathscr{R}(\delta) = \mathscr{R} \cdot e_\delta$, où $\delta \colon \mathbf{Q}_p^* \to L^*$ est un caractère continu, et $\varphi(e_\delta) = \delta(p) e_\delta$, $\sigma(e_\delta) = \delta(\chi(\sigma)) e_\delta$ si $\sigma \in \Gamma_{\mathbf{Q}_p}$.

(ii) Si $K \neq \mathbf{Q}_p$, tout $(\varphi, \Gamma)$-module de rang 1 sur $\mathscr{R}_L$ est de la forme $\mathscr{R}_L(\delta) = \mathscr{R}_L \cdot e_\delta$, où $\delta \colon K^* \to L^*$ est un caractère continu, mais les formules pour $\varphi(e_\delta)$ et $\sigma(e_\delta)$, pour $\sigma \in \Gamma_K$, ne sont pas vraiment explicites (cf. [169]).

(iii) Si $V$ est trianguline, et si $D_{\text{rig}}(V)$ est une extension $\mathscr{R}_L(\delta_1) - \cdots - \mathscr{R}_L(\delta_n)$ (avec $\delta_1$ dans le socle et $\delta_n$ dans le cosocle), le th. 4.12 ci-dessous fournit des contraintes assez fortes sur les $v_p(\delta_i(p))$ : par exemple, si $n = 2$, on doit avoir $v_p(\delta_1(p)) + v_p(\delta_2(p)) = 0$ et $v_p(\delta_1(p)) \geq v_p(\delta_2(p))$ (ces deux conditions sont presque suffisantes pour qu'une extension de $\mathscr{R}_L(\delta_2)$ par $\mathscr{R}_L(\delta_1)$ soit de la forme $D_{\text{rig}}(V)$).

THÉORÈME 4.12. — (KEDLAYA [145]). *Si $\Delta$ est un $(\varphi, \Gamma)$-module sur $\mathbf{B}_{\text{rig},K}$, il existe une unique filtration $0 = \Delta_0 \subset \Delta_1 \subset \cdots \subset \Delta_h = \Delta$ et, pour tout $1 \leq i \leq h$, un unique sous-$\mathbf{B}_K^\dagger$-module $D_i^\dagger$ de $\Delta_i/\Delta_{i-1}$, non nul, stable par $\varphi$ et $\Gamma_K$ et isocline[70] de pente $r_i \in \mathbf{Q}$, tels que :*
- $r_1 < r_2 < \cdots < r_h$,
- $\Delta_i/\Delta_{i-1} = \mathbf{B}_{\text{rig},K} \otimes_{\mathbf{B}_K^\dagger} D_i^\dagger$.

*Remarque 4.13.* — (i) Une représentation potentiellement semi-stable est trianguline si et seulement si elle devient semi-stable sur une extension abélienne (ce qui équivaut à ce que la représentation de $\mathbf{GL}_n(K)$ qui lui correspond via la recette du (iv) de la rem. 4.6 soit de la série principale, i.e. sous-quotient d'une induite d'un caractère du Borel).

(ii) Contrairement aux représentations du (i) dont les poids de Hodge–Tate sont des entiers, et donc sont fixes dans une famille analytique, les représentations trianguliness peuvent se mettre dans des familles de poids variable : si $K = \mathbf{Q}_p$, CHENEVIER [50, th. B] a prouvé que l'espace des représentations trianguliness est une variété analytique de dimension $\frac{n(n+1)}{2} + 1$ ; NAKAMURA [170] a généralisé le résultat à $K$ quelconque (et la dimension devient $[K : \mathbf{Q}_p] \frac{n(n+1)}{2} + 1$).

(iii) LIU [164] a prouvé qu'une famille analytique de $(\varphi, \Gamma)$-modules, point par point triangulable, admet une triangulation en famille (voir aussi [146]).

### 4.4. Théorie de Hodge $p$-adique entière

Les théorèmes de relèvement modulaire demandent une bonne compréhension des propriétés de théorie de Hodge $p$-adique (cristallinité, semi-stabilité, etc.) au niveau entier (i.e. pour des $\mathbf{Z}_p$ (ou $\mathscr{O}_L$)-représentations de $\text{Gal}_K$ au lieu de $\mathbf{Q}_p$ (ou $L$)-représentations), ou même, idéalement, pour les représentations de torsion (puisqu'on part d'une représentation $\overline{\rho}$ en caractéristique $p$ et que l'on essaie de comprendre ses relèvements mod $p^2$, $p^3$, etc.). La théorie des $(\varphi, \Gamma)$-modules fournit, en principe, une telle description pour toutes les représentations, mais cette universalité fait qu'il est difficile de reconnaître celles qui nous intéressent vraiment, à savoir les représentations potentiellement semi-stables.

---

70. Cela signifie que $\mathbf{B}_K \otimes_{\mathbf{B}_K^\dagger} D_i^\dagger$, muni de $p^{-r_i}\varphi$ (au lieu de $\varphi$), est étale.

*4.4.1. La théorie de Fontaine–Laffaille.* — Cette théorie fournit une description des réseaux des représentations cristallines et de leurs quotients de torsion, mais seulement dans le cas où $K$ est non ramifié, et les poids de Hodge–Tate sont dans un intervalle de la forme $[a, a+p-2]$ (il y a une obstruction évidente à aller plus loin : les caractères 1 et $\chi^{p-1}$ — respectivement de poids de Hodge–Tate 0 et $p-1$ — ont la même réduction modulo $p$ ; si on impose qu'il n'y ait pas de quotient de la forme $(\mathbf{Z}/p^n)(\chi^{a+p-1})$, on peut passer à un intervalle de la forme $[a, a+p-1]$).

On définit un *module de Fontaine–Laffaille* sur $\mathscr{O}_K$ (où $K$ est non ramifié), comme un $\mathscr{O}_K$-module $M$, de type fini, muni :
- d'une filtration décroissante par des sous-$\mathscr{O}_K$-modules $M^i$, $i \in \mathbf{Z}$, avec $\cup_i M^i = M$ et $\cap_i M^i = 0$,
- de frobenius $\varphi_i \colon M^i \to M$, semi-linéaires, avec $\varphi_i = p\varphi_{i+1}$ sur $M^{i+1}$, et [71] $M = \sum_i \varphi_i(M^i)$.

On dit que $M$ est à poids dans $[a, b]$, si $M^a = M$ et $M^{b+1} = 0$.

THÉORÈME 4.14. — (FONTAINE, LAFFAILLE [106]). *La catégorie des modules de Fontaine–Laffaille de torsion, à poids dans $[a, a+p-2]$, est équivalente à celle des sous-quotients de torsion des réseaux de représentations cristallines de $\mathrm{Gal}_K$, à poids dans $[a, a+p-2]$.*

*Remarque 4.15.* — (i) On peut aller jusqu'à $a+p-1$, si on interdit les sous-quotients $N$ avec $N = N^{a+p-1}$ du côté des modules de Fontaine–Laffaille, et ceux de la forme $(\mathbf{Z}/p^n)(\chi^{a+p-1})$ du côté des représentations de $\mathrm{Gal}_K$.

(ii) Si [72] $a = 0$, l'équivalence est $M \mapsto \mathbf{V}^*_{\mathrm{cris}}(M) := \mathrm{Hom}(M, \mathbf{A}_{\mathrm{cris}} \otimes (\mathbf{Q}_p/\mathbf{Z}_p))$, où l'on considère les morphismes de $\mathscr{O}_K$-modules, compatibles aux filtrations [73], et commutant avec [74] les $\varphi_i$.

En passant à la limite, on obtient une équivalence de catégories $M \mapsto \mathbf{V}^*_{\mathrm{cris}}(M)$ entre les modules de Fontaine–Laffaille, sans torsion, à poids dans $[0, p-2]$ (voire $[0, p-1]$ avec une condition supplémentaire) et les réseaux de représentations cristallines de $\mathrm{Gal}_K$ à poids de Hodge–Tate dans $[0, p-2]$ (voire $[0, p-1]$) ; si on inverse $p$, un module de Fontaine–Laffaille devient un $\varphi$-module filtré et $\mathbf{V}^*_{\mathrm{cris}}(M)[\frac{1}{p}]$ est la duale de la représentation $\mathbf{V}_{\mathrm{st}}(M[\frac{1}{p}])$ du th. 4.5.

---

71. Cette dernière condition est l'analogue en niveau entier, de la faible admissibilité.
72. Le cas général se déduit du cas $a = 0$ par torsion.
73. La filtration sur $\mathbf{A}_{\mathrm{cris}} \otimes (\mathbf{Q}_p/\mathbf{Z}_p)$ est induite par celle sur $\mathbf{A}_{\mathrm{cris}}$ pour $i \leq p-1$ et nulle pour $i \geq p$.
74. $\varphi_i = p^{-i}\varphi$ sur $\mathbf{A}^i_{\mathrm{cris}}$ — la raison pour imposer $(\mathbf{A}_{\mathrm{cris}} \otimes (\mathbf{Q}_p/\mathbf{Z}_p))^i = 0$ si $i \geq p$ est que $\varphi$ n'est pas divisible par $p^i$ sur $\mathbf{A}^i_{\mathrm{cris}}$, si $i \geq p$.

(iii) Wach [213] a donné une autre preuve du th. 4.14, via la théorie des $(\varphi, \Gamma)$-modules. Berger [14] a étendu la classification des réseaux des représentations cristallines qui en résulte aux représentations de poids arbitraires (pas limités à $[0, p-2]$) mais toujours pour $K$ non ramifié.

Pour énoncer le résultat de Berger, définissons un *module de Wach* sur $\Lambda = \mathbf{A}_K^+, \mathbf{A}_K^+[\frac{1}{p}]$ comme un $\Lambda$-module libre $M$, muni d'une action de $\Gamma_K$ triviale sur $M/\pi M$, et d'un $\varphi \colon M[\frac{1}{\pi}] \to M[\frac{1}{\varphi(\pi)}]$ commutant à l'action de $\Gamma_F$, tel que, si $b$ est assez grand, il existe $a$ tel que $\varphi(\pi)^a M/\varphi^*(\pi^b M)$ soit tué par une puissance de $\frac{\varphi(\pi)}{\pi}$ $(= \varphi(\xi))$.

**Théorème 4.16.** — (Berger [14]). *Soit $K$ une extension non ramifiée de $\mathbf{Q}_p$.*

(i) *Si $T$ est un réseau d'une représentation cristalline de $\mathrm{Gal}_K$, alors $D(T)$ contient un unique sous-module de Wach $N(T)$ tel que $D(T) = \mathbf{A}_K \otimes_{\mathbf{A}_K^+} N(T)$.*

(ii) *Le foncteur $T[\frac{1}{p}] \to N(T)[\frac{1}{p}]$ est une équivalence de catégories entre les représentations cristallines de $\mathrm{Gal}_K$ et les modules de Wach sur $\mathbf{B}_K^+$.*

*4.4.2. Modules de Breuil–Kisin.* — Le (ii) du th. 4.16 précise une conjecture de Fontaine selon laquelle les représentations cristallines d'une extension non ramifiée $K$ de $\mathbf{Q}_p$ sont *de hauteur finie*. Son handicap pour les applications aux théorèmes de relèvement modulaire est la restriction « $K$ non ramifiée ».

Dans le cas où $K$ est ramifiée, on dispose des modules de Breuil–Kisin, une variante des modules de Wach (la définition, due à Kisin [151], est une interpolation entre des objets définis par Breuil [31, 32] et les modules de Wach) utilisant l'extension de Kummer $K_\infty^{\mathrm{Kum}} = K(\varpi^{1/p^\infty})$ (où $\varpi$ est une uniformisante fixée de $K$) au lieu de l'extension cyclotomique. Une des raisons *a posteriori* pour lesquelles la théorie marche est que l'on ne perd pas d'information en restreignant une représentation cristalline de $\mathrm{Gal}_K$ au sous-groupe $H_K^{\mathrm{Kum}}$ fixant $K_\infty^{\mathrm{Kum}}$ (Kisin [151] déduit ceci de la classification via les modules de Breuil–Kisin mais c'est en fait élémentaire [11]).

Soit $E \in \mathscr{O}_{K_0}[X]$ le polynôme minimal de $\varpi$ sur $K_0$ (normalisé par $E(0) = p$). Notons $\mathbf{A}_K^{\mathrm{Kum},+}$ le sous-anneau $\mathscr{O}_{K_0}[[u]]$ de $\widetilde{\mathbf{A}}^+$, où $u = [\varpi^\flat]$. Alors $E(u)$ est un générateur[75] de $\mathrm{Ker}\,\theta$ (i.e. $\widetilde{\mathbf{A}}^+/E(u) = \mathscr{O}_{\mathbf{C}_p}$), et $\mathbf{A}_K^{\mathrm{Kum},+}/E(u) = \mathscr{O}_K$.

Un module de Breuil–Kisin $M$ est un $\mathbf{A}_K^{\mathrm{Kum},+}$-module de la forme $M_1/M_2$ où $M_1$ et $M_2$ sont libres et de type fini ($M_2$ peut être nul), muni d'un $\varphi$ semi-linéaire, tel que $M/\varphi^*M$ soit tué par une puissance de $E(u)$. On dit que $M$ est de hauteur $\leq h$ si $M/\varphi^*M$ est tué par $E(u)^h$. Le résultat suivant est dû à Kisin [151] ; il est modelé sur les résultats de Breuil [32].

---

75. C'est l'équivalent de $\xi = \varphi^{-1}(\frac{\varphi(\pi)}{\pi})$ ; il y a donc une torsion par $\varphi$ par rapport à ce qui sort des $(\varphi, \Gamma)$-modules ; notons aussi que l'élément $\lambda = \prod_{n \geq 0} \varphi^n(p^{-1}E(u))$ qui joue un grand rôle dans les formules de Kisin [151] est relié à $t$ : si $\alpha \in \widetilde{\mathbf{A}}^+ \setminus p\widetilde{\mathbf{A}}^+$ vérifie $\varphi(\alpha) = E(u)\alpha$, alors $t = \alpha\lambda$ (à multiplication près par un élément de $\mathbf{Z}_p^*$), et on a $\sigma(\alpha) = \chi(\sigma)\alpha$ pour tout $\sigma \in H_K^{\mathrm{Kum}}$.

Théorème 4.17. — *Si $p \geq 3$, il y a des équivalences de catégories entre :*
- *modules de Breuil–Kisin de torsion, de hauteur $\leq 1$,*
- *sous-quotients de torsion des représentations cristallines de $\mathrm{Gal}_K$ à poids $0$ et $-1$,*
- *schémas en groupe finis, plats sur $\mathscr{O}_K$.*

*Remarque 4.18.* — (i) En passant à la limite, on obtient des équivalences entre modules de Breuil–Kisin libres de hauteur $\leq 1$, réseaux de représentations cristallines à poids $0$ et $-1$, et groupes $p$-divisibles sur $\mathscr{O}_K$.

(ii) Si $p = 2$, on obtient une résultat analogue en se restreignant aux objets connexes [153].

Kisin [151] donne aussi une classification [76] des représentations semi-stables de $\mathrm{Gal}_K$ en termes de modules de Breuil–Kisin, libres sur $\mathbf{A}_K^{\mathrm{Kum},+}$, munis d'un opérateur $N$ sur $(M/u)\frac{1}{p}$ vérifant $N\varphi = p\varphi N$. Cette classification fournit une description des réseaux d'une représentation cristalline stable par $H_K$ mais pas une description des réseaux stables par $\mathrm{Gal}_K$. Des classifications vraiment en niveau entier ont été obtenues récemment par Bhatt, Scholze [21] (via une approche prismatique) et par Gao [117] (en termes de modules de Breuil–Kisin enrichis, cf. ci-dessous), en se basant sur son travail avec Poyeton [118], ainsi que par Du, Liu [85].

On définit *un $\mathrm{Gal}_K$-module de Breuil–Kisin* comme un module de Breuil–Kisin $M$, libre sur $\mathbf{A}_K^{\mathrm{Kum},+}$, muni d'une action semi-linéaire de $\mathrm{Gal}_K$ sur $\widetilde{\mathbf{A}}^+ \otimes M$ telle que $M$ soit fixe par $H_K^{\mathrm{Kum}} \subset \mathrm{Gal}_K$ et $M/uM$ soit fixe par $\mathrm{Gal}_K$.

Théorème 4.19. — (Gao [117]). *Les $\mathrm{Gal}_K$-modules de Breuil–Kisin sont en équivalence de catégorie avec les réseaux des représentations semi-stables de $\mathrm{Gal}_K$ à poids de Hodge–Tate $\leq 0$.*

*Remarque 4.20.* — (i) Il est un peu tôt pour prédire quel impact ces résultats auront sur les théorèmes de relèvement modulaire. Ces derniers demandent plutôt une description des sous-quotients de torsion des réseaux des représentations semi-stables et les bons objets à considérer [21, note 3] seraient les $\varphi$-jauges, un vieux rêve de Fontaine, Jannsen [107] revisité par Drinfeld [84] et Bhatt, Lurie [17, 18].

(ii) Les théorèmes de comparaison au niveau entier, raffinant le th. 4.1, ont fleuri ces dernières années avec les travaux de Bhatt, Česnavičius, Koshikawa, Lurie, Morrow, Scholze [19, 49, 20, 17] (entre autres), et les inventions de la $\mathbf{A}_{\mathrm{inf}}$-cohomologie et de la cohomologie prismatique.

### 4.5. Familles de représentations

*4.5.1. Anneaux de déformations.* — Soit $K$ une extension finie de $\mathbf{Q}_p$, et soient $k$ un corps fini de caractéristique $p$, $\overline{\rho}\colon \mathrm{Gal}_K \to \mathbf{GL}_d(k)$ une représentation continue, et $L$ une extension finie de $\mathbf{Q}_p$ de corps résiduel $k$ et d'anneau des entiers $\mathscr{O}_L$. On suppose que $\mathscr{O}_L$ est suffisamment grand ; en particulier qu'il contient $\boldsymbol{\mu}_{p^\infty}(K)$.

---

76. Cette classification revient à dire que toute représentation semi-stable est de hauteur finie vue de l'extension de Kummer (alors que, vue de l'extension cyclotomique, ce n'est le cas que pour certaines, si $K$ est ramifiée).

Notons $R_{\bar\rho}^\square$ l'anneau des $\mathscr{O}_L$-déformations encadrées de $\bar\rho$ et $\mathscr{X}_{\bar\rho}^\square := \operatorname{Spf} R_{\bar\rho}^\square$ le schéma formel associé.

Théorème 4.21. — (Böckle, Iyengar, Paškūnas [24, 25]).

(i) *L'anneau $R_{\bar\rho}^\square$ est une $\mathscr{O}_L$-algèbre locale réduite, d'intersection complète, plate sur $\mathscr{O}_L$, de dimension relative $d^2 + d^2[K : \mathbf{Q}_p]$. En particulier, $\bar\rho$ a un relèvement en caractéristique 0.*

(ii) *Le morphisme $R_{\det\bar\rho} \to R_{\bar\rho}^\square$ induit par $\rho \mapsto \det\rho$ est plat, et fournit une bijection entre composantes connexes des fibres génériques de $\mathscr{X}_{\bar\rho}^\square$ et $\mathscr{X}_{\det\bar\rho}$ qui, elles-mêmes, sont en bijection avec les caractères $\chi \colon \boldsymbol{\mu}_{p^\infty}(K) \to \mathscr{O}_L^*$.*

*Ces composantes connexes sont aussi en bijection avec les composantes irréductibles de $\operatorname{Spec} R_{\bar\rho}^\square$ et, si $X$ est une de ces composantes irréductibles, $\mathscr{O}(X)$ et $\mathscr{O}(X)/\varpi_L$ sont normales, d'intersection complète.*

(iii) *Les points cristallins sont Zariski-denses dans l'espace analytique associé à la fibre générique de $\mathscr{X}_{\bar\rho}^\square$.*

*Remarque 4.22.* — (i) L'existence d'un relèvement de $\bar\rho$ en caractéristique 0 est très facile à prouver en dimension 2, mais surprenamment difficile en dimension plus grande. La première preuve de ce résultat a été obtenue par Emerton et Gee [91] (qui ont de plus prouvé qu'il existait un relèvement cristallin).

(ii) Chenevier [50], Nakamura [170] ont démontré des résultats partiels concernant le (iii) ; cet énoncé admet une variante à déterminant fixé, qui est souvent utile.

La densité des cristallines permet de déduire certains résultats pour une représentation encadrée $\rho \colon \operatorname{Gal}_K \to \mathbf{GL}_d(L)$ quelconque par prolongement analytique à partir du cas des représentations cristallines (c'est ainsi que la correspondance de Langlands locale $p$-adique pour $\mathbf{GL}_2(\mathbf{Q}_p)$ a été établie [67]).

(iii) On trouve des sous-ensembles Zariski-denses plus petits dans [24, th. 1.6], sous la condition $p \nmid 2d$.

*4.5.2. Anneaux de Kisin.* — Soit $H = (H_\tau)_\tau$, où $\tau$ décrit $\operatorname{Hom}(K, L)$, une famille de multi-ensembles d'entiers, avec $H_\tau = \sum_i m_{\tau,i}\{i\}$ de cardinal $n$ pour tout $\tau$. On va s'intéresser aux $L$-représentations potentiellement semi-stables de $\operatorname{Gal}_K$ de poids de Hodge–Tate $H$.

Soit $\alpha \colon I_K \to \mathbf{GL}_n(L)$ une représentation localement constante. On dit que $V$ est *de type $(H, \alpha)$*, si $V$ est de Rham (et donc potentiellement semi-stable), de poids de Hodge–Tate $H$ et si $\mathbf{D}_{\mathrm{pst}}(V) \cong \alpha$ comme représentation de $I_K$.

Théorème 4.23. — (Kisin [152]). *Soit $\rho^{\mathrm{univ}} \colon \operatorname{Gal}_K \to \mathbf{GL}_n(R_{\bar\rho}^\square)$ la déformation universelle encadrée de $\bar\rho$.*

(i) *Il existe un unique quotient réduit $R_{\bar\rho}^{\square, H, \alpha}$ de $R_{\bar\rho}^\square$, sans $p$-torsion, tel que $\rho^{\mathrm{univ}} \otimes_{R_{\bar\rho}^\square} E$ est de type $(H, \alpha)$ si et seulement si $R_{\bar\rho}^\square \to E$ se factorise par $R_{\bar\rho}^{\square, H, \alpha}$.*

(ii) $\mathrm{Spec}(R_{\overline{\rho}}^{\square,H,\alpha}[\frac{1}{p}])$ *est équidimensionnel, de dimension* $n^2+\sum_\tau\left(\frac{n(n-1)}{2}-\sum_i\frac{m_{\tau,i}(m_{\tau,i}-1)}{2}\right)$.

*Remarque 4.24*. — (i) Ce résultat peut se reformuler en : les $\rho\colon \mathrm{Gal}_K \to \mathbf{GL}_n(\overline{\mathbf{Q}}_p)$ relevant $\overline{\rho}$, de type $(H,\alpha)$, forment une sous-variété algébrique de $\mathscr{X}_{\overline{\rho}}^\square$ (ce qui est plus précis que les résultat de [15] qui prouvent que ces $\rho$ forment une sous-variété analytique de l'espace analytique associé à la fibre générique de $\mathscr{X}_{\overline{\rho}}^\square$).

(ii) La dimension est maximum si les poids de Hodge–Tate sont réguliers ; dans ce cas cette dimension est $n^2 + [K:\mathbf{Q}_p]\frac{n(n-1)}{2}$.

(iii) Si on s'intéresse aux représentations non encadrées, il faut quotienter par l'action de $\mathbf{GL}_n$ par conjugaison et la dimension chute de $n^2-1$ (si $\mathrm{End}\,\overline{\rho}=\kappa$, on obtient une variété ; dans le cas contraire, on obtient un champ algébrique).

(iv) Le sous-espace des représentations potentiellement cristallines (i.e le sous-espace sur lequel $N=0$) est une réunion de composantes irréductibles et est lisse.

*4.5.3. La conjecture de Breuil–Mézard*. — Si les poids de Hodge–Tate sont petits et réguliers, et si $\alpha$ n'est pas trop ramifié, $\mathscr{X}_{\overline{\rho}}^{\square,H,\alpha}$ est juste une boule ouverte de dimension $n^2 + [K:\mathbf{Q}_p]\frac{n(n-1)}{2}$. Cette simplicité est trompeuse, et la géométrie des variétés de Kisin est, en général, très compliquée. Dans le cas général, on dispose de la conjecture de Breuil–Mézard donnant, sous sa forme géométrique due à EMERTON, GEE [90], des renseignements précieux sur les composantes irréductibles de la fibre spéciale (et donc aussi sur celles de la fibre générique, ce qui est très utile pour prouver des théorèmes de relèvement modulaire).

Soit $K$ une extension finie de $\mathbf{Q}_p$ et $n \geq 2$. Un *poids de Serre* (pour $K,n$) est une $\overline{\mathbf{F}}_p$-représentation irréductible $W$ de $\mathbf{GL}_n(\mathscr{O}_K)$ (ou, ce qui revient au même, de $\mathbf{GL}_n(\kappa_K)$). Si $K = \mathbf{Q}_p$ et $n=2$, ce sont les $\mathrm{Sym}^k \otimes \det^a$, pour $0 \leq k \leq p-1$ et $0 \leq a \leq p-2$, où $\mathrm{Sym}^k$ est la puissance symétrique $k$-ième de la représentation de dimension 2 de $\mathbf{GL}_2(\mathbf{Z}_p)$ agissant à travers son quotient $\mathbf{GL}_2(\mathbf{F}_p) \hookrightarrow \mathbf{GL}_2(\overline{\mathbf{F}}_p)$ ; il y a une description analogue en général, faisant intervenir les représentations algébriques de $\mathbf{GL}_n$.

Si $\lambda = (\lambda_1,\ldots,\lambda_n)$ avec $\lambda_i \in \mathbf{Z}$ et $\lambda_1 \geq \lambda_2 \geq \cdots \geq \lambda_n$, on associe à $\lambda$ une représentation algébrique $M_\lambda$ du groupe $\mathbf{GL}_n$ (induite d'un caractère algébrique du Borel). Si $H = (H_\tau)_\tau$, avec $H_\tau$ régulier, on peut écrire $H_\tau$ sous la forme $(\lambda_{\tau,1}+n-1, \lambda_{\tau,2}+n-2, \ldots, \lambda_{\tau,n})$ avec $\lambda_{\tau,1} \geq \lambda_{\tau,2} \geq \cdots \geq \lambda_{\tau,n}$. On note $M_{H,\tau}$ la représentation de $\mathbf{GL}_n(\mathscr{O}_K)$ obtenue via le plongement $\tau\colon \mathbf{GL}_n(\mathscr{O}_K) \hookrightarrow \mathbf{GL}_n(L)$, et en faisant agir $\mathbf{GL}_n(L)$ sur $M_{\lambda(\tau)}(L)$ où $\lambda(\tau) = (\lambda_{\tau,1},\ldots \lambda_{\tau,n})$. On pose alors $M_H = \otimes_\tau M_{H,\tau}$ ; c'est une représentation de dimension finie de $\mathbf{GL}_d(\mathscr{O}_K)$.

THÉORÈME 4.25. — (Correspondance de Langlands locale inertielle [12, pro. 6.5.3]).

*Si $\alpha\colon I_K \to \mathbf{GL}_n(L)$ est une représentation lisse, il existe une $L$-représentation lisse irréductible $\mathrm{LL}(\alpha)$ de $\mathbf{GL}_n(\mathscr{O}_K)$ avec la propriété suivante : si $\rho$ est une représentation de dimension $n$ de $\mathrm{WD}_K$, et si $\mathrm{LL}(\rho)$ est la représentation de $\mathbf{GL}_n(K)$ associée par la correspondance de Langlands locale, alors $\mathrm{Hom}_{\mathbf{GL}_n(\mathscr{O}_K)}(\mathrm{LL}(\alpha),\mathrm{LL}(\rho)) \neq 0$ si et seulement si la restriction de $\rho$ à $I_K$ est $\alpha$ et $N=0$.*

(Pour $n=2$, ce résultat a été prouvé par HENNIART dans un appendice à l'article de BREUIL, MÉZARD [36].)

La représentation $\mathrm{LL}(\alpha)\otimes M_H$ est irréductible ; elle possède un $\mathscr{O}_L$-réseau stable, et la réduction modulo $\mathfrak{m}_L$ de ce réseau est lisse. Sa semi-simplifiée s'écrit donc sous la forme $\oplus_W W^{\oplus m_W(H,\alpha)}$ (et les $m_W(H,\alpha)$ ne dépendent pas du choix du réseau).

Par ailleurs, d'après le th. 4.23, la dimension de $R_{\overline{\rho}}^{\square,H,\tau}[\frac{1}{p}]$ est $n^2+[K:\mathbf{Q}_p]\frac{n(n-1)}{2}$ (cf. (ii) de la rem. 4.24) ; comme $R_{\overline{\rho}}^{\square,H,\tau}$ est sans $p$-torsion, on en déduit que la fibre spéciale est de la même dimension, mais ses composantes irréductibles peuvent avoir une multiplicité. On note $Z(\overline{\mathbf{F}}_p\otimes R_{\overline{\rho}}^{\square,H,\tau})$ le cycle de dimension $n^2+[K:\mathbf{Q}_p]\frac{n(n-1)}{2}$ de $\mathrm{Spec}(\overline{\mathbf{F}}_p\otimes R_{\overline{\rho}}^{\square})$ correspondant.

CONJECTURE 4.26. — *Il existe des cycles $Z_{\overline{\rho},W}$ de dimension $n^2+[K:\mathbf{Q}_p]\frac{n(n-1)}{2}$ de $\mathrm{Spec}(\overline{\mathbf{F}}_p\otimes R_{\overline{\rho}}^{\square})$ tels que, pour tout $(H,\alpha)$, on ait $Z(\overline{\mathbf{F}}_p\otimes R_{\overline{\rho}}^{\square,H,\tau})=\sum_W m_W(H,\tau)Z_{\overline{\rho},W}$*

Les cycles $Z_{\overline{\rho},W}$, s'ils existent, sont surdéterminés : ils sont solutions d'une infinité d'équations en un nombre fini d'indéterminées.

Cette conjecture est connue si $n=2$ et $K=\mathbf{Q}_p$ (divers cas sont traités dans [154, 179] ; cf. [70] pour une version avec un arrière-goût catégorique). Les résultats les plus généraux sont ceux de FENG, LE, LE HUNG, LEVIN, MORRA [163, 99] ($n$ arbitraire mais $K$ non ramifiée et des restrictions sérieuses sur $\overline{\rho}, H, \alpha$ ; [99] considère aussi la conj. 4.29 ci-dessous).

*4.5.4. Le champ d'Emerton–Gee.* — Si $n\geq 1$, EMERTON, GEE [91] ont montré que les $(\varphi,\Gamma)$-modules de rang $n$ sur $\mathbf{A}_K$ sont paramétrés[77] par un champ formel[78] $\mathscr{X}_n$, noethérien. Le champ réduit sous-jacent $\mathscr{X}_{n,\mathrm{red}}$ est un champ algébrique, de dimension[79] $[K:\mathbf{Q}_p]\frac{n(n-1)}{2}$. Les points de $\mathscr{X}_{n,\mathrm{red}}(\overline{\mathbf{F}}_p)$ paramètrent les représentations continues $\overline{\rho}\colon \mathrm{Gal}_K\to \mathbf{GL}_n(\overline{\mathbf{F}}_p)$.

Si $\mathrm{End}\,\overline{\rho}=\overline{\mathbf{F}}_p$, le complété de $\mathscr{O}_{\mathscr{X}_n}$ en $\overline{\rho}$ est l'anneau $R_{\overline{\rho}}$ des déformations universelles de $\overline{\rho}$. Dans le cas général, on dispose d'un morphisme naturel $\mathrm{Spf}\,R_{\overline{\rho}}^{\square}\to \mathscr{X}_n$ qui est une surjection sur le tube de $\overline{\rho}$.

*Exemple 4.27.* — Si $n=1$, $\mathscr{X}_n$ est l'espace des caractères unitaires continus de $K^*$. On a $K^*\cong \boldsymbol{\mu}(K)\times \varpi^{\mathbf{Z}}\times \mathbf{Z}_p^{[K:\mathbf{Q}_p]}$, et $\boldsymbol{\mu}(K)=\boldsymbol{\mu}'(K)\times \boldsymbol{\mu}_{p^\infty}(K)$ ; posons $|\boldsymbol{\mu}_{p^\infty}(K)|=p^k$. Un caractère de $K^*$ est donc donné par sa valeur sur $\boldsymbol{\mu}'(K)$ (les caractères de $\boldsymbol{\mu}'(K)$ ont pour réduction modulo $p$ les poids de Serre pour $K,1$), sur un générateur de $\boldsymbol{\mu}_{p^\infty}(K)$ (une racine $p^k$-ième de l'unité), sur la base de $\mathbf{Z}_p^{[K:\mathbf{Q}_p]}$ (soit un élément de la boule unité ouverte de centre 1 pour chaque élément de la base) et sur $\varpi$ (une unité).

---

77. Si $\Lambda$ est une algèbre noethérienne $I$-adiquement complète, avec $p\in I$, un $(\varphi,\Gamma)$-module sur $\Lambda\widehat{\otimes}\mathbf{A}_K$ est équivalent à la donnée d'un morphisme $\mathrm{Spf}\,A\to \mathscr{X}_n$.

78. Mais pas formel $p$-adique : les $\mathscr{O}_{\mathscr{X}_n}$ ressemblent à $\mathbf{Z}_p[[x_1,\ldots,x_s]]\langle x_{s+1},\ldots,x_d\rangle$, i.e. certaines directions sont $p$-adiques et d'autres sont purement formelles.

79. Il s'agit de dimension champêtre : cette dimension est 1 de moins que ce à quoi on s'attendrait à cause du $\mathbf{G}_m$ contenu dans les automorphismes des objets considérés.

Il s'ensuit que $\mathscr{X}_1$ est une réunion disjointe de $|\boldsymbol{\mu}'(K)|$ copies de

$$\mathrm{Spf}\Big(\big(\mathbf{Z}_p[[x_0-1, x_1-1, \ldots, x_{[K:\mathbf{Q}_p]}-1]]/(x_0^{p^k}-1)\big)\langle y, y^{-1}\rangle\Big)$$

Le schéma $\mathscr{X}_{1,\mathrm{red}}$ s'obtient en tuant les variables topologiquement nilpotentes (les $x_i - 1$ et $p$), et est une réunion disjointe de $|\boldsymbol{\mu}'(K)|$ copies de [80] $\mathrm{Spec}\,\mathbf{F}_p[y, y^{-1}]$; ses $\overline{\mathbf{F}}_p$-points paramètrent les représentations continues $\bar\rho \colon \mathrm{Gal}_K \to \mathbf{GL}_1(\overline{\mathbf{F}}_p)$.

Par contraste, l'espace des représentations continues de dimension 1 de $\mathrm{Gal}_K$, est l'espace des caractères continus du complété profini $\widehat{K^*}$ de $K^*$; c'est la réunion des complétés de $\mathscr{X}_1$ en chacun des points de $\mathscr{X}_{1,\mathrm{red}}(\overline{\mathbf{F}}_p)$ et l'espace réduit sous jacent est la réunion des points de $\mathscr{X}_{1,\mathrm{red}}(\overline{\mathbf{F}}_p)$, sans géométrie.

Les composantes irréductibles de $\mathscr{X}_n$ sont en bijection naturelle avec les poids de Serre pour $(K, n)$; on note $\mathscr{X}_n^W$ la composante correspondant à $W$. L'image inverse dans $\mathrm{Spf}\,R_{\bar\rho}^{\square}$ de $\mathscr{X}_{n,\mathrm{red}}^W$ est la complétion d'un sous-schéma fermé $Z_{\bar\rho}^W$ de $\mathrm{Spec}\,R_{\bar\rho}^{\square}/p$ de dimension $n^2 + [K:\mathbf{Q}_p]\frac{n(n-1)}{2}$.

THÉORÈME 4.28. — [91, th. 1.7.1]. *Si $H$ est régulier, le support de $Z(\overline{\mathbf{F}}_p \otimes R_{\bar\rho}^{\square,H,\tau})$ est contenu dans la réunion des $Z_{\bar\rho}^W$.*

Ce résultat fournit une confirmation forte de la conj. 4.26 et suggère que les $Z_{\bar\rho,W}$ de cette conjecture sont des combinaisons linéaires à coefficients entiers $\geq 0$ des $Z_{\bar\rho}^{W'}$.

On peut aller encore plus loin dans la géométrisation de la conjecture de Breuil–Mézard : EMERTON, GEE [91, th. 4.8.12] construisent, pour tout type $(H, \alpha)$, un sous-champ fermé (et $p$-adique) $\mathscr{X}_n^{(H,\alpha)}$ paramétrant les $(\varphi, \Gamma)$-modules dont la représentation de $\mathrm{Gal}_K$ associée est potentiellement semi-stable de type $(H, \alpha)$; l'image inverse de $\mathscr{X}_n^{(H,\alpha)}$ dans $\mathscr{X}_{\bar\rho}^{\square}$ est l'espace $\mathscr{X}_{\bar\rho}^{\square,H,\alpha}$ ci-dessus.

CONJECTURE 4.29. — [91, conj. 8.2.2]. *Il existe des cycles $Z_W$, combinaisons à coefficients entiers $\geq 0$ des $\mathscr{X}_{n,\mathrm{red}}^{W'}$, tels que le cycle sous-jacent à la fibre spéciale de $\mathscr{X}_n^{(H,\alpha)}$ soit $\sum_W m_W(H,\alpha) Z_W$.*

*Remarque 4.30.* — (i) Des résultats assez généraux confirmant cette conjecture ont été annoncés par FENG, LE HUNG [99].

(ii) Une justification particulièrement séduisante de la conjecture est fournie par la catégorification (conjecturale) de la correspondance de Langlands locale $p$-adique pour $\mathbf{GL}_n(K)$ (dont l'existence est tout aussi conjecturale, sauf si $n = 2$ et $K = \mathbf{Q}_p$), cf. point (4) de [92, conj. 6.1.14] et [92, n° 6.1.37].

---

80. Il est de dimension 1 mais $\mathscr{X}_{1,\mathrm{red}}$ est de dimension 0 car il faut quotienter par $\mathbf{G}_m$ agissant trivialement.

# RÉFÉRENCES

Pierre Colmez

CNRS, IMJ-PRG, Sorbonne Université,
4 place Jussieu, 75010 Paris, France
*E-mail* : `pierre.colmez@imj-prg.fr`